\documentclass[11pt, leqno]{amsart}
\usepackage{latexsym, amsfonts,mathrsfs, amsmath, amssymb, amscd, epsfig}
\usepackage{mathrsfs}

\usepackage{amsthm}
\usepackage{array}
\usepackage{psfrag}
\usepackage{bm}
\usepackage{setspace}
\usepackage[usenames,dvipsnames]{xcolor}

\numberwithin{equation}{section}

\newtheorem{lemma}{Lemma}[section]
\newtheorem{theorem}[lemma]{Theorem}
\newtheorem{proposition}[lemma]{Proposition}
\newtheorem{corollary}[lemma]{Corollary}
\newtheorem{definition}[lemma]{Definition}
\newtheorem{remark}[lemma]{Remark}
\newtheorem{problem}[lemma]{Problem}

\theoremstyle{definition}

\def\beq#1\eeq{\begin{equation}#1\end{equation}}
\def\balign #1 #2 \ealign{\begin{aligned} #1 #2  \end{aligned} }

\def\Div{{\rm div}}

\def\bu{\mathbf{u}}

\def\bfF{\mathbf{F}}

\def\mE{\mathcal{E}}
\def\e{\mathfrak{e}}
\def\i{\mathfrak{i}}

\newcommand \alp{\alpha}
\newcommand \eps{\varepsilon}
\newcommand \vphi{\varphi}

\newcommand \Gam{\Gamma}
\newcommand \gam{\gamma}
\newcommand \om{\omega}
\newcommand \tx{\text}
\newcommand \R{\mathbb{R}}
\newcommand \til{\tilde}

\newcommand \der{\partial}
\newcommand \mcl{\mathcal}

\newcommand \ol{\overline}
\newcommand \Om{\Omega}

\newcommand \Gamen{\Gamma_0}
\newcommand \Gamex{\Gamma_L}

\newcommand \rx{{\rm{x}}}

\newcommand \bara{\bar{a}}

\newcommand \tpsi{\til{\psi}}
\newcommand \tPsi{\til{\Psi}}

\def \msB {\mathscr{B}}
\def \msK {\mathscr{K}}

\def\bfX{\mathbf{X}}

\def\Div{{\rm div}}
\def\Curl{{\rm curl}}

\newcommand \rhos{\rho_{s}}

\newcommand \ba{\bar a}
\newcommand \bb{\bar b}

\newcommand \baru{\bar u}

\newcommand \barE{\bar E}

\newcommand \mcW {\mathcal{W}}

\newcommand \hpsi{\hat{\psi}}
\newcommand \hPsi{\hat{\Psi}}

\newcommand \mfrak{\mathfrak}

\begin{document}
\title[supersonic solutions to Euler-Poisson system]
{Structural stability of Supersonic solutions to the Euler-Poisson system}

\author{Myoungjean Bae}
\address{Myoungjean Bae, Department of Mathematics\\
         POSTECH\\
          San 31, Hyojadong, Namgu, Pohang, Gyungbuk, Republic of Korea 37673;
         Korea Institute for Advanced Study
85 Hoegiro, Dongdaemun-gu,
Seoul 130-722,
Republic of Korea
}
\email{mjbae@postech.ac.kr}

\author{Ben Duan}
\address{B. Duan, School of Mathematical Sciences, Dalian University of Technology, Dalian, 116024, China}
\email{bduan@dlut.edu.cn}
\author{Jingjing Xiao}
\address{J. Xiao, Department of Mathematics, The Chinese University of Hong Kong, Shatin, Hong Kong
}
\email{jjxiao@math.cuhk.edu.hk}

\author{Chunjing Xie}
\address{C. Xie, School of Mathematical Sciences, Institute of Natural Sciences, Ministry of Education Key Laboratory of Scientific and Engineering Computing, and SHL-MAC, Shanghai Jiao Tong University, 800 Dongchuan Road, Shanghai, China
}
\email{cjxie@sjtu.edu.cn}

\begin{abstract}
The well-posedness for the supersonic solutions of the Euler-Poisson system for hydrodynamical model in semiconductor devices and plasmas is studied in this paper. We first reformulate the Euler-Poisson system in the supersonic region into a second order hyperbolic-elliptic coupled system together with several transport equations. One of the key ingredients of the analysis is to obtain  the well-posedness of the boundary value problem for the associated linearized hyperbolic-elliptic coupled system, which is achieved via a delicate choice of multiplier to gain energy estimate. The nonlinear structural stability of supersonic solution in the general situation is established by combining the iteration method with the estimate for hyperbolic-elliptic system and the transport equations together.

\end{abstract}


\thanks{Updated on \today}

\maketitle


\section{Introduction}

The hydrodynamical model of semiconductor devices or plasmas is described by the nonlinear system
\begin{equation}\label{UnsteadyEP}
\left\{
\begin{aligned}
& \rho_t+ \Div_{\rm x} (\rho \bu)=0, \\
& (\rho \bu)_t+\Div_{\rm x} (\rho \bu \otimes \bu) +\nabla_{\rm x} p=\rho \nabla_{\rm x} \Phi, \\
& (\rho \mE)_t +\Div_{\rm x}(\rho\mE \bu +p\bu)=\rho \bu\cdot \nabla_{\rm x}\Phi,\\
& \Delta_{\rm x} \Phi=\rho-b({\rm x}),%
\end{aligned}%
\right.
\end{equation}
called the {\emph{Euler-Poisson system}} (see \cite{MarkRSbook}). In the system above, $\bu,
\rho$, $p$, and $\mE$ represent the macroscopic particle velocity, electron density,
pressure, and the total energy density, respectively.  The electric potential $\Phi$ is  generated by the
Coulomb force of particles. The fixed positive function $b({\rx})>0$ represents the density of fixed,
positively charged background ions. In fact, the system \eqref{UnsteadyEP} can also be used to model the biological transport of ions in channel proteins \cite{Shu}.
The system \eqref{UnsteadyEP} is closed with the aid of definition of specific total energy and the equation of state
\begin{equation}
\mE=\frac{|\bu|^2}{2}+\e\quad \text{and}\quad p=p(\rho, \e),
\end{equation}
respectively,
where $\e$ is the internal energy.
In this paper, we consider the case for which the pressure $p$ and \emph{the enthalpy} $\i=\e+\frac{p}{\rho}$  are given by
\begin{equation}
\label{2d-1-b1}
p(\rho, S)=S {\rho^{\gam}}\quad \text{and}\quad
\i(\rho, S)= \frac{\gam}{\gam-1}S\rho^{\gam-1},
\end{equation}
respectively, where we follow the notations in gas dynamics to call the constant $\gamma>1$ the adiabatic constant and the quantity $\ln S$ entropy.
One of the interesting phenomenon for the system \eqref{UnsteadyEP} is that the electric field can provide more stabilizing effect. Mathematically speaking, when $b(x)\equiv b_0$ for some constant $b_0$, the associated linearized system around trivial steady state $(\rho,\bu, S)=(b_0, {\bf 0}, \bar S)$  where
$\bar S$ is a constant state, is a Klein-Gordon type system (equation) which has faster dispersive decay than the wave system (equation) which corresponds to the linearized Euler system. With the aid of this faster dispersive decay and the nice structure of the system \eqref{UnsteadyEP}, the global classical solutions of the system \eqref{UnsteadyEP} with small and smooth irrotational data were established in  \cite{Guo99, LW,IP, GHZ} in the cases with different spatial dimensions.

A natural problem is to see whether there are some other physically nontrivial steady states. If there are some nontrivial steady solutions, can we prove  the stability of these solutions?
The steady state of the system \eqref{UnsteadyEP}  is governed by  the following steady Euler-Poisson system
\begin{equation}\label{2-a2}
\begin{cases}
\Div_{\rm x} (\rho \bu)=0,\\
\Div_{\rm x} (\rho \bu \otimes \bu) +\nabla p=\rho \nabla_{\rm x} \Phi,\\
\Div_{\rm x}(\rho\bu \msB)=\rho {\bf u}\cdot \nabla_{\rm x}\Phi,\\
\Delta_{\rm x}\Phi=\rho-b(\rx),
\end{cases}
\end{equation}
where {\emph{Bernoulli's function}} $\msB$ is given by
\begin{equation}
\label{2d-1-a5}
\msB(\rho,|{\bf u}|, S)=\frac{|\bu|^2}{2}+\e+\frac{p}{\rho}
=
\frac{|\bu|^2}{2}+\frac{\gam}{\gamma-1}S\rho^{\gam-1} .
\end{equation}
The one dimensional solutions of the system \eqref{2-a2} were studied in \cite{MarkPhase, Gamba1d, GambaMorawetz, LuoXin, RosiniPhase}, and the structural and dynamical stability  of some transonic shock solutions in one dimensional setting were achieved in \cite{LRXX}. A natural question is to study the structural stability of nontrivial one dimensional steady solutions for the Euler-Poisson system under multidimensional steady perturbations of the boundary conditions.

The first main difficulty for the system \eqref{2-a2} is that it may change type as long as the Mach number $M_a=|\bu|/\sqrt{p_\rho(\rho, S)}$ of the flows varies, where $\sqrt{p_\rho(\rho, S)}$ is called the {\emph{sound speed}}. In the case $M_a<1$, i.e., the flow is subsonic, if, in addition, the current flux is sufficiently small, the existence of subsonic solutions was obtained in \cite{DeMark1d, DeMark3d, MarkZAMP,Yeh, Hattori, WEP}.
In \cite{BDX, BDX3}, the structural stability of subsonic solutions for the Euler-Poisson system has been achieved even when the background solutions have large variations.
More precisely, in \cite{BDX}, the subsonic potential flow model of the system \eqref{2-a2} was formulated as a second order quasilinear elliptic system for  the velocity potential and  the electric potential, respectively. Furthermore, it was interestingly discovered that this quasilinear elliptic system has a special structure to yield the well-posedness even when a nonlinear boundary condition, fixing the exit pressure, is prescribed. For general solutions  with nonzero vorticity, in \cite{BDX3}, two dimensional subsonic solutions to the system \eqref{2-a2}  have been studied with the aid of the Helmholtz decomposition. Hence a natural question is to understand the well-posedness for the supersonic solutions of the steady Euler-Poisson system \eqref{2-a2}.
It is our goal to provide an answer to this question in this work.

If the right hand sides of the second and the third equations in the system \eqref{2-a2} are zero and the Poisson equation in \eqref{2-a2} disappears,  then the system \eqref{2-a2} becomes the well-known steady Euler system which describes gas motion of compressible inviscid flow. If $M_a>1$, the steady Euler system is a quasilinear hyperbolic system. Therefore, the unique existence of smooth multidimensional supersonic solutions to steady Euler system in a nozzle of finite length
can be obtained by the standard theory for the initial boundary value problem for hyperbolic system \cite{Rauch}. However,  the system \eqref{2-a2} becomes a second order quasilinear  hyperbolic-elliptic coupled system even in the case of isentropic irrotational solutions as the Poisson equation (the fourth equation in \eqref{2-a2}) is included. Hence steady Euler-Poisson system for supersonic solutions is analytically different very much from both steady Euler system for supersonic solutions and steady Euler-Poisson system for subsonic solutions. Therefore, the question on the well-posedness theory for the multi-dimensional supersonic solutions to steady Euler-Poisson system becomes immediately nontrivial.

The main goal of the analysis is to establish structural stability of supersonic solutions to the Euler-Poisson system under multidimensional perturbations of boundary conditions.
One of the key ingredients of this paper is to establish the well-posedness of a boundary value problem of a class of quasilinear hyperbolic-elliptic coupled system which is the linearized problem for the Euler-Poisson system around one dimensional supersonic solutions. So far, there are very few known results about general hyperbolic-elliptic system of second order PDEs. In \cite{Anderson-Monc}, the local and strong well-posedness of a Riemann problem of a quasilinear hyperbolic-elliptic system is proved. But as far as we know, there has been no known result about the initial-boundary value problem for a general quasilinear hyperbolic-elliptic system  except for the study on the symmetric positive system in \cite{Fried}. On the other hand,
the study for supersonic solutions for the Euler-Poisson system not only is important for the supersonic solutions themselves but also plays a crucial role for the study on the transonic solutions for the Euler-Poisson system.
The study of subsonic and supersonic solutions to the system \eqref{2-a2} is the essential step to analyze multi-dimensional smooth or discontinuous transonic solutions. This is also one of our ultimate goal for the study on steady solutions of the Euler-Poisson system.

The rest of the paper is organized as follows. In Section \ref{secPb}, we propose the problem on the existence of multidimensional  supersonic solutions for the Euler-Poisson system and state the main results. In Section \ref{section-thm-pf-potential}, we introduce a linear hyperbolic-elliptic coupled system of second order and establish the well-posedness of the boundary value problem of this system. This is the core part of the work. By using the results obtained in  Section \ref{section-thm-pf-potential}, the nonlinear structural stability of supersonic solutions are achieved in Section \ref{secNL}  for both the irrotational case and the case with nonzero vorticity.  Finally,  the detailed analysis on the higher order estimate for the linearized problem of the Euler-Poisson system and the analysis on the transport equation for the entropy are included in the appendices.

\section{Problem and Main theorems}\label{secPb}
\subsection{One dimensional supersonic solutions for steady Euler-Poisson system}
The analysis on various types of one dimensional solutions to the steady Euler-Poisson system \eqref{2-a2} has been achieved in \cite{LuoXin}. Here, we recall and collect some important properties of supersonic solutions which are used in this work later.

Let $b_0>0$ be a fixed constant.
Assume that $(\rho, u, p, \Phi)$ solve the ODE system
\begin{equation}
\label{1-b1}
\begin{cases}
(\rho u)'=0,\\
(\rho u^2+p)'=\rho \Phi',\\
(\rho u \msB)'=\rho u \Phi',\\
\Phi''=\rho-b_0,
\end{cases}
\end{equation}
where $'$ denotes the derivative with respect to $x_1$ and $\msB=\msB(\rho, |u|, S)$ is defined by \eqref{2d-1-a5}. If $\rho>0$ and $u>0$, then the system \eqref{1-b1} is equivalent to
\begin{equation}
\label{ode-s}
\begin{cases}
\rho'=\frac{E\rho}{\gam {S_0}\rho^{\gam-1}-\frac{J_0^2}{\rho^2}},\\
E'=\rho-b_0
\end{cases}
\end{equation}
with
\begin{equation}\label{equS}
u=\frac{J_0}{\rho}\quad \text{and}\quad \frac{p}{\rho^{\gam}}={S_0},
\end{equation}
where the  constants $J_0>0$ and $S_0>0$ are determined by the values of $(\rho, u, p)$ at $x_1=0$. Here, $E=\Phi'$ represents the electric field in $x_1$-direction. The first equation in \eqref{ode-s} has a singularity at $\rho=\rhos$ which is given by
\begin{equation}
\label{def-rhos}
\rhos=\left(\frac{J_0^2}{\gam S_0}\right)^{\frac{1}{\gam+1}}.
\end{equation}
One can directly check from \eqref{equS} that the solution $(\rho, u, p, E)$ to \eqref{1-b1} is {\emph{supersonic}} if and only if $\rho<\rhos$, {\emph{subsonic}} if and only if $\rho>\rhos$, and {\emph{sonic}} if and only if $\rho=\rhos$.

\begin{figure}[htp]
\centering
\includegraphics[scale=0.5]{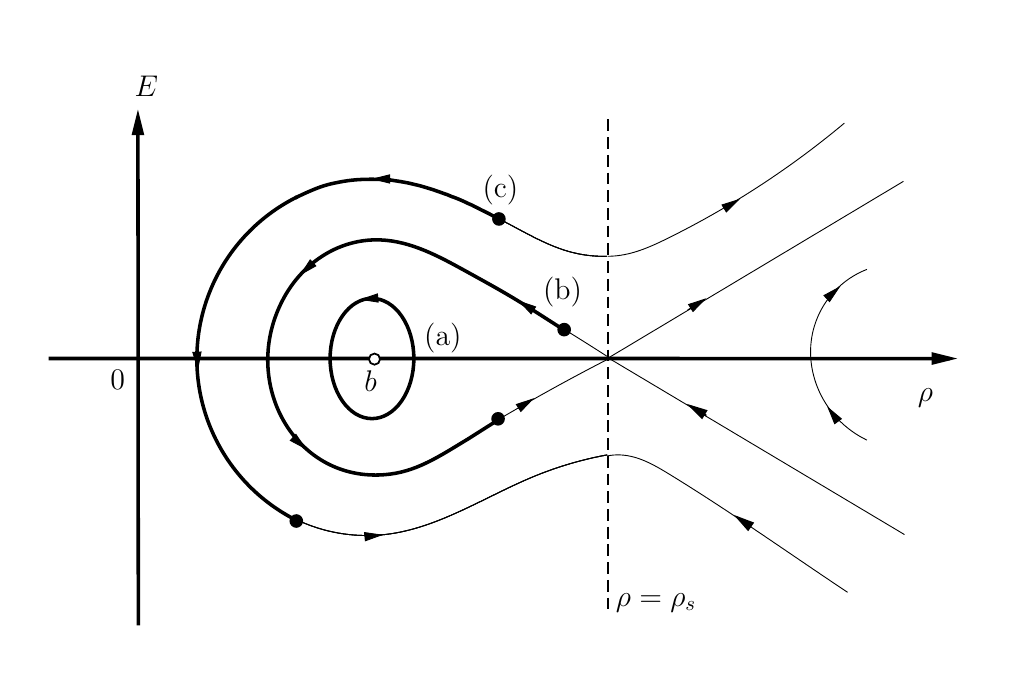}
\caption{Phase plane of the system \eqref{ode-s} }
\label{figure-1}
\end{figure}

It is easy to see that the system \eqref{ode-s} has a Hamiltonian. In fact, any $C^1$ solution $(\rho, E)$ to \eqref{ode-s} satisfies
\begin{equation}
\label{eqn-Ham}
\frac 12 E^2-H(\rho)={\rm constant}
\end{equation}
for
\begin{equation}
\label{def-H}
H(\rho)=\int_{\rhos}^\rho\frac{(t-b_0)}{t}\left(\gam S_0 t^{\gam-1}-\frac{J_0^2}{t^2}\right)\;dt.
\end{equation}
The phase plane of the system \eqref{ode-s} under the assumption of
\begin{equation}
\label{assumption}
0<b_0<\rhos
\end{equation}
is given in Figure \ref{figure-1} for fixed $(J_0, S_0, b_0)$.

We first study the solution of the ODEs \eqref{ode-s} with the initial data
\begin{equation}\label{ODE_IC}
(\rho, E)(0) =(\rho_0, E_0)
\end{equation}
satisfying $0<\rho_0<\rhos$.

Here we collect some results for supersonic solutions in \cite{LuoXin} and exploit more properties of the solutions of the problem \eqref{ode-s} and \eqref{ODE_IC} which we need for later analysis.
\begin{lemma}
\label{lemma-1dsol-new}
Fix $\gam>1$, $J_0>0$, and $S_0>0$. Suppose that $\rho_0\in(0, \rhos)$ and $E_0\in \R$.
\begin{itemize}
\item[(a)] If $\frac 12 E_0^2-H(\rho_0)<0$, then the initial value problem \eqref{ode-s} and \eqref{ODE_IC}  has a unique solution $(\rho, u)(x_1)\in C^{\infty}(\R)$. And, the solution is periodic. 
Furthermore, there exists a constant $\bar{\eps}>0$ depending on $(\gam, J_0, S_0, \rho_0, E_0)$ such that
    \begin{equation*}
       \bar{\eps}\le\rho(x_1)\le \rhos-\bar{\eps}\quad\tx{for all $x_1\in \R$}.
    \end{equation*}

\item[(b)] If $\frac 12 E_0^2-H(\rho_0)=0$, then there exist constants $T_{\min}$ and $T_{\max}$ with $-\infty<T_{\min}<0<T_{\max}<+\infty$ so that the initial value problem \eqref{ode-s} and \eqref{ODE_IC} has a unique solution $(\rho, u)(x_1)\in C^\infty([T_{\min}, T_{\max}])$ satisfying
    \begin{equation*}
      0<\rho(x_1)<\rhos\,\,\tx{on $(T_{\min}, T_{\max})$},\quad \rho(x_1)=\rhos\,\,\tx{at }x_1=T_{\min}\,\,\text{and}\,\, T_{\max}.
    \end{equation*}

\item[(c)] If $\frac 12 E_0^2-H(\rho_0)>0$, then there exist constants $T_{\min}$ and $T_{\max}$ with $-\infty<T_{\min}<0<T_{\max}<+\infty$ so that the initial value problem \eqref{ode-s} and \eqref{ODE_IC} has a unique solution $(\rho, u)(x_1)\in C^\infty(T_{\min}, T_{\max})$ satisfying
    \begin{equation*}
      0<\rho(x_1)<\rhos\,\,\tx{on $(T_{\min}, T_{\max})$},\quad \lim_{x_1\to T_{\min}+}\rho(x_1)=\lim_{x_1\to T_{\max}-}\rho(x_1)=\rhos,
      \end{equation*}
  and
  \begin{equation*}
      \lim_{x_1\to T_{\min}+}\rho'(x_1)=-\infty,\quad \lim_{x_1\to T_{\max}-}\rho'(x_1)=+\infty.
    \end{equation*}
\end{itemize}
\begin{proof}
It directly follows from the unique existence theorem of ODEs that the initial value problem \eqref{ode-s} and \eqref{ODE_IC} for $\rho_0\in (0, \rhos)$ has a unique smooth solution as long as $0<\rho(x_1)<\rhos$.
Denote
\begin{equation*}
\begin{split}
&\Upsilon_{(\rho_0, E_0)}:=\left\{(\rho, E)\in (0, \rhos]\times \R\Big|\frac 12 E^2-H(\rho)=\frac 12 E_0^2-H(\rho_0)\right\}.\\ 
\end{split}
\end{equation*}
First, we have
$\displaystyle{
\inf_{(\rho, E)\in \Upsilon_{(\rho_0, E_0)}}\rho>0}
$
for each $(\rho_0, E_0)\in (0,\rhos)\times \R$.\\
\quad\\
\emph{(i) Proof of (a)}: If $\frac 12 E_0^2-H(\rho_0)<0$, then we have $\displaystyle{
\sup_{(\rho, E)\in \Upsilon_{(\rho_0, E_0)}}\rho<\rhos}
$. Therefore,  the initial value problem \eqref{ode-s} and \eqref{ODE_IC} has a unique smooth solution on $\R$. Heuristically,
the periodicity can be observed from Fig. \ref{figure-1} (See the orbit (a) in Fig. \ref{figure-1}).
For further details, one can refer to \cite{LuoXin}. Hence (a) is proved. \\

\emph{(ii) Proof of (b)}: If $\frac 12 E_0^2-H(\rho_0)=0$, then we have (see the orbit (b) of Fig. \ref{figure-1})
\begin{equation}
\label{ineq-Upsilon-inf-sup}
0<\inf_{(\rho, E)\in \Upsilon_{(\rho_0, E_0)}}\rho
<\sup_{(\rho, E)\in \Upsilon_{(\rho_0, E_0)}}\rho=\rhos.
\end{equation}
The only equilibrium solution of \eqref{ode-s} is $(\rho, E)=(b_0, 0)$.
It follows from \eqref{def-H} and \eqref{assumption} that one has $H(b_0)>0=H(\rhos)$. Hence $(b_0,0)\not\in \Upsilon_{(\rho_0, E_0)}$ and $(\rhos,0)\in \Upsilon_{(\rho_0, E_0)}$. Therefore, there exist two finite constants $T_{\min}$ and $T_{\max}$ with $-\infty<T_{\min}<0<T_{\max}<\infty $ so that the solution $(\rho, E)(x_1)$ to the initial value problem \eqref{ode-s} and \eqref{ODE_IC} satisfy $0<\rho(x_1)<\rhos$ for $T_{\min}<x_1< T_{\max}$, and $\rho(T_{\min})=\rho(T_{\max})=\rhos$. Then, $(\rho, E)(x_1)$ is $C^{\infty}$ for $T_{\min}<x_1<T_{\max}$. Now it remains to prove that $(\rho, E)(x_1)$ is $C^{\infty}$ at $x_1=T_{\min}$
and $x_1= T_{\max}$.

Denote
\begin{equation*}
\begin{split}
&\bar\Upsilon:=\left\{(\rho, E)\Big|\frac 12 E^2-H(\rho)=0\right\},\\
&\bar\Upsilon_+:=\left\{(\rho, E)\in \bar\Upsilon\Big|-(\rho-\rhos)E\ge 0\right\},\quad
\bar\Upsilon_-:=\left\{(\rho, E)\in \bar\Upsilon\Big|-(\rho-\rhos)E\le 0\right\}.
\end{split}
\end{equation*}
Note that $(\rhos,0)\in \bar\Upsilon_+\cap \bar\Upsilon_-$.
By \eqref{def-rhos} and \eqref{def-H}, we have $H(\rhos)=H_{\rho}(\rhos)=0$. Therefore, any $(\rho, E)\in \bar\Upsilon$ satisfies
\begin{equation}\label{extra1}
E^2=2(\rho-\rhos)^2\int_0^1\int_0^1 t_1H_{\rho\rho}(t_1t_2\rho+(1-t_1t_2)\rhos)\,dt_1dt_2.
\end{equation}
And, a direct computation with using \eqref{def-H} yields that
\begin{equation}\label{extra2}
H_{\rho\rho}(\rho)=\frac{\gam S_0}{\rho^3}
\left((\rho^{\gam+1}-\rhos^{\gam+1})(1-\frac{3(\rho-b_0)}{\rho})+(\gam+1)(\rho-b_0)\rho^{\gam}\right).
\end{equation}
By \eqref{def-rhos}, we get
\begin{equation}\label{extra3}
  \gam S_0\rho^{\gam-1}-\frac{J_0^2}{\rho^2}
  =\frac{\gam S_0}{\rho^2}(\rho^{\gam+1}-\rhos^{\gam+1})
  =\frac{\gam(\gam+1) S_0}{\rho^2}(\rho-\rhos)\int_0^1(t\rho+(1-t)\rhos)^{\gam}\,dt.
\end{equation}
We use \eqref{extra1}--\eqref{extra3} to rewrite the differential equation for $\rho$ in \eqref{ode-s} as
\begin{equation*}
  \rho'=-F(\rho)\,\,\tx{on $\bar\Upsilon_+$},\quad \rho'=F(\rho)\,\,\tx{on $\bar\Upsilon_-$}
\end{equation*}
where
\begin{equation}
\label{ode-rho-new-ctt}
  F(\rho)=\frac{\rho^3\sqrt{2\int_0^1\int_0^1 t_1H_{\rho\rho}(t_1t_2\rho+(1-t_1t_2)\rhos)\,dt_1dt_2}}
  {\gam(\gam+1)S_0\int_0^1 (t\rho+(1-t)\rhos)^{\gam}\,dt}.
\end{equation}
We regard $(\rho(T_{\min}),0)\in \bar\Upsilon_-\cap \Upsilon_{(\rho_0, E_0)}$ and $(\rho(T_{\max}),0)\in \bar\Upsilon_+\cap \Upsilon_{(\rho_0, E_0)}$.
Since $F(\rho)$ is smooth with respect to $\rho$, we conclude from the differential equation \eqref{ode-rho-new-ctt} for $\rho$ that $\rho(x_1)$ is $C^{\infty}$ up to $x_1=T_{\min}$, and up to $x_1=T_{\max}$.\\
\quad\\
{\emph{(iii) Proof of (c)}}:
If $\frac 12 E_0^2-H(\rho_0)>0$ holds, then ${F}(\rho, E): =\frac{E\rho}{\gam {S_0}\rho^{\gam-1}-\frac{J_0^2}{\rho^2}}$ is continuous on the set
$\Upsilon_{(\rho_0,E_0)}\setminus \{(\rho, E): \rho=\rhos\}$ but discontinuous for $\rho=\rhos$ (see the orbit (c) in Fig. \ref{figure-1}).  Since the curve $\Upsilon_{(\rho_0, E_0)}$ passes through the points $(\rhos, \pm\sqrt{E_0^2-2H(\rho_0)})$, we have $|\rho'|=\infty$ if $\rho=\rhos$ and $\frac 12 E_0^2-H(\rho_0)>0$. Then, the statement (c) follows from the unique existence theorem of ODEs and the property $\displaystyle{\inf_{\Upsilon_{(\rho_0, E_0)}}|{F}(\rho, E)|>0}$.
\end{proof}
\end{lemma}

\begin{definition}\label{def1}
(i) If  $\frac 12 E_0^2-H(\rho_0)<0$, that is, the solution $(\rho, u)$ of \eqref{ode-s} and \eqref{ODE_IC} is periodic, we denote
\begin{equation}
\begin{split}
&T_{\max}=\sup \{x_1: \rho(t)< \sup_{(\rho, E)\in \Upsilon_{(\rho_0, E_0)}}\rho, \,\,\tx{for all $t\in (0, x_1)$}\},\\
&T_{\min}=\inf \{x_1: \rho(t)< \sup_{(\rho, E)\in \Upsilon_{(\rho_0, E_0)}}\rho, \,\,\tx{for all $t\in (x_1, 0)$}\}.
\end{split}
\end{equation}
(ii) If $E_0>0$, we define
\begin{equation}
T_*=\sup\{ x_1: E(t)>0, \,\,\text{for all $t\in (0,x_1)$}\}.
\end{equation}
\end{definition}

It is easy to see that $T_{\max}$ is the smallest positive value of $x_1$ at which $\rho(x_1)$ achieves its maximum, and that $T_*$ is the smallest positive value of $x_1$ at which the value of $E(x_1)$ becomes zero. Furthermore, the flow corresponding to the solution $(\rho, E)$ of \eqref{ode-s} and \eqref{ODE_IC} is accelerating on the interval $(T_{\min}, T_*)$ as long as $E_0>0$.

\subsection{Problem and Main Results}
\label{subsec-main-thm1}
For a constant $L>0$, denote
\begin{equation}\label{defineomega}
\Om_L:=\{(x_1,x_2): 0<x_1<L, \,\,-1<x_2<1\},
\end{equation}
and
\begin{equation}\label{defineomegabdy}
\Gam_0:=\der\Om_L\cap\{x_1=0\},\quad \Gamex:=\der\Om_L\cap\{x_1=L\},\quad \Lambda_L:=\der \Om_L\cap\{x_2=\pm 1\}.
\end{equation}
If $\rho>0$ in $\Om_L$, then it can be directly derived from the steady Euler-Poisson system \eqref{2-a2} that
\[
{\bf u} \cdot \nabla S={\bf u}\cdot \nabla \msK=0\quad\tx{in $\Om_L$}
\]
for $\msK:=\msB-\Phi$,  where $\msB$ is defined in \eqref{2d-1-a5}. We call $\msK$ {\emph{the pseudo-Bernoulli invariant}}.
Under the condition $\msK=0$ on $\Gamen$, $u_1>0$ in $\Omega_L$, and the slip boundary condition on $\Lambda_L$,
the transport equation ${\bf u}\cdot \nabla \msK=0$ yields $\msK\equiv 0$ in $\Om_L$. Since these two transport equations
for $S$ and $\msK$ can be dealt with in the same way, without of loss of generality, we assume $\msK\equiv 0$ in this paper for simplicity. 

Our main goal is to study the following well-posedness problem for the supersonic solutions of the steady Euler-Poisson system.
\begin{problem}
\label{problem-1}
Given functions $(b, u_{\rm en}, v_{\rm en}, S_{\rm en}, E_{\rm en}, \Phi_{\rm ex})$ with $\displaystyle{\min_{\ol{\Gamen}}u_{\rm en}>0 }$, find
$(\rho, {\bf u}, p, \Phi)\in [C^1(\Om_L)\cap C^0(\ol{\Om_L})]^4\times [C^2(\Om_L)\cap C^1(\ol{\Om_L})]$ with ${\bf u}=(u_1, u_2)$  so that
\begin{itemize}
\item[(i)] $(\rho, {\bf{u}}, p, \Phi)$ solve the system \eqref{2-a2} in $\Om_L$;
\item[(ii)]  $\rho>0$ and $u_1>0$ hold in $\ol{\Om_L}$;
\item[(iii)] $(\rho, {\bf u}, p, \Phi)$ satisfy boundary conditions:
\begin{equation}
\label{bc-1}
\begin{split}
(u_1, u_2, S, \Phi_{x_1}, \msK)=(u_{\rm en}, v_{\rm en}, S_{\rm en}, E_{\rm en},0)\quad&\tx{on $\Gamen$},\\
u_2=\Phi_{x_2}=0 \quad &\tx{on $\Lambda_L$},\\
\Phi=\Phi_{\rm ex} \quad &\tx{on $\Gamex$};
\end{split}
\end{equation}

\item[(iv)] The inequality $|{\bf u}|^2>\frac{\gam p}{\rho}$ holds in $\ol{\Om_L}$, i.e., the flow corresponding to $(\rho, {\bf u}, p, \Phi)$ is supersonic in $\ol{\Om_L}$.
\end{itemize}
\end{problem}
It follows from \cite{BDX2} that if $u_1> 0$ and $\msK\equiv 0$, then the two dimensional steady Euler-Poisson system is equivalent to
\begin{equation}\label{EP_1}
\left\{
\begin{aligned}
&\Div(\rho \bu)=0,\\
&\Delta \Phi =\rho -b,\\
& \omega =\frac{\rho^{\gamma-1}S_{x_2}}{(\gamma-1) u_1},\\
& \rho\bu \cdot \nabla S=0,\\
& \msB-\Phi\equiv 0,
\end{aligned}
\right.
\end{equation}
where $\omega=\Curl \bu =\partial_{x_1} u_2-\partial_{x_2}u_1$ is the vorticity of the flow. It follows from $\msB-\Phi= 0$ that
\begin{equation}\label{extra4}
\rho=\left(\frac{\gamma-1}{\gamma S} \left(\Phi-\frac{|\bu|^2}{2}\right)\right)^{\frac{1}{\gamma-1}}.
\end{equation}
Substituting \eqref{extra4} into \eqref{EP_1} implies that \eqref{EP_1} can be regarded as a nonlinear  PDE system for $({\bf u}, S, \Phi)$.

To solve the system \eqref{EP_1}, we introduce the Helmholtz decomposition for the velocity field as in \cite{BDX2}. Given a velocity field $\bu\in C^1(\Om_L)\cap C^0(\ol{\Om_L})$,
if $\phi\in C^2(\Om_L)\cap C^1(\ol{\Om_L})$ satisfies
\begin{equation}
\label{hdecomp-1}
\begin{cases}
\Delta \phi =\omega \qquad\tx{in}\,\,\Om_L,\\
\der_{x_1}\phi=0\qquad\tx{on}\;\;\Gamen,\qquad \phi=0\qquad\tx{on}\;\;\der\Om_L\setminus \Gamen,
\end{cases}
\end{equation}
then
$
\Curl (\bu -\nabla^\perp \phi)=0
$ holds in $\Om_L$
for $\nabla^\perp\phi :=(-\partial_{x_2}\phi, \partial_{x_1}\phi)$. Since $\Om_L$ is a simply connected domain, one can find a function $\varphi\in C^2(\Om_L)$ satisfying
\begin{equation*}
\nabla \varphi=\bu -\nabla^\perp \phi.
\end{equation*}
In other words, the velocity field $\bu$ can be expressed as
\begin{equation}
\label{hdecomp-2}
\bu=\nabla\vphi+\nabla^{\perp}\phi.
\end{equation}
Therefore, the system \eqref{EP_1} can be written into a nonlinear system for $(\vphi, \phi, \Phi, S, \msK)$ as follows,
\begin{equation}
\label{hdecomp-system}
\left\{
\begin{aligned}
&{\rm div}\left(\mcl{H}\left(S, \Phi-\frac 12|\nabla\vphi+\nabla^{\perp}\phi|^2\right)(\nabla\vphi+\nabla^{\perp}\phi)\right)=0,\\
&\Delta\Phi=\mcl{H}\left(S, \Phi-\frac 12|\nabla\vphi+\nabla^{\perp}\phi|^2\right)-b({\bf x}),
\end{aligned}
\right.
\end{equation}
together with
\begin{equation}
\label{eqn3-hdecomp}
\Delta\phi=\frac{ \mcl{H}^{\gamma-1}\left(S, \Phi-\frac 12|\nabla\vphi+\nabla^{\perp}\phi|^2\right)S_{x_2}}{(\gam-1)(\vphi_{x_1}+\phi_{x_2})},
\end{equation}
and
\begin{equation}
\label{eqn4-hdecomp}
\mcl{H}\left(S, \Phi-\frac 12|\nabla\vphi+\nabla^{\perp}\phi|^2\right)(\nabla\vphi+\nabla^{\perp}\phi)\cdot \nabla S=0,
\end{equation}
where
\begin{equation}
\label{rho-hdecomp}
\mcl{H}(S,\zeta)=
\left(\frac{\gam-1}{\gam}\frac{\zeta}{S}\right)^{\frac{1}{\gam-1}},
\,\,\,\,\text{provided that}\,\,\,\,\frac{\zeta}{S}>0.
\end{equation}
In terms of  $(\vphi, \phi, \Phi, S)$, the boundary conditions  \eqref{bc-1} can be formulated as
\begin{equation}\label{bc-vphi}
\vphi_{x_1}=u_{\rm en}(x_2)-\phi_{x_2},\,\, \vphi(0,x_2)=\varphi_{\rm en}(x_2)\,\,\tx{on}\,\, \Gamen,\quad
\vphi_{x_2}=0\,\, \tx{on}\,\, \Lambda_L,
\end{equation}
\begin{equation}
\label{bc-Phi}
\Phi_{x_1}=E_{\rm en}\,\,\tx{on}\,\, \Gamen,\quad
\Phi_{x_2}=0\,\,\tx{on}\,\, \Lambda_L,\quad
\Phi=\Phi_{\rm ex}\,\,\tx{on}\,\, \Gamex,
\end{equation}
\begin{equation}
\label{bc-phi}
\phi_{x_1}=0\,\, \tx{on}\,\, \Gamen,\quad \phi=0\,\, \tx{on}\,\, \der\Om_L\setminus \Gamen,
\end{equation}
\begin{equation}
\label{bc-T}
S =S_{\rm en}\,\,\tx{on}\,\, \Gamen,
\end{equation}
where
\begin{equation}
\label{def-ent-ftn}
\varphi_{\rm en}(x_2)=\int_{-1}^{x_2} v_{\rm en}(s)\;ds
\end{equation}
with $v_{\rm en}$ given in Problem \ref{problem-1}.

Suppose that
\begin{equation}
\label{bc-T-simplified}
  S_{\rm en}\equiv S_0\quad\tx{on $\Gamen$}
\end{equation}
for a constant $S_0>0$.  In the rest of the paper, we denote $\{{\bf e}_1, {\bf e}_2\}$ to be the canonical bases of $\mathbb{R}^2$.  If $\mcl{H}\left(S, \Phi-\frac 12|\nabla\vphi+\nabla^{\perp}\phi|^2\right)
(\nabla\vphi+\nabla^{\perp}\phi)\cdot{\bf e}_1>0$ in $\ol{\Om_L}$, then it follows from \eqref{eqn4-hdecomp}, \eqref{bc-T}, \eqref{bc-T-simplified} and the slip boundary condition for ${\bf u}\cdot {\bf e_2}=0$ on $\Lambda_L$, which is achieved by the boundary conditions $\vphi_{x_2}=\phi=0$ on $\Lambda_L$, that the characteristic method gives
\begin{equation*}
  S\equiv S_0\quad \tx{in $\ol{\Om_L}$}.
\end{equation*}
In this case, Eq. \eqref{eqn3-hdecomp} becomes $\Delta\phi=0$ in $\Om_L$, and this equation combined with \eqref{bc-phi} implies
$\phi\equiv 0$ in $\ol{\Om_L}$. Hence
the system of \eqref{hdecomp-system}--\eqref{eqn4-hdecomp} is simplified as
\begin{align}
\label{pfeqn-mass}
&{\rm div} (\mcl{H}_0(\Phi,\nabla\vphi)\nabla\vphi)=0,\\
\label{pfeqn-poisson}
&\Delta \Phi=\mcl{H}_0(\Phi, \nabla\vphi)-b(\rx),
\end{align}
where
\begin{equation}\label{def-rho}
\mcl{H}_0(z, {\bf q})=\left(\frac{\gam-1}{\gam S_0}\left(z-\frac{|{\bf q}|^2}{2}\right) \right)^{\frac{1}{\gam-1}}\qquad \text{for}\,\,z\in \mathbb{R},\,\, {\bf q}\in \mathbb{R}^2\,\,\tx{with $z-\frac{|{\bf q}|^2}{2}>0$}.
\end{equation}
The system \eqref{pfeqn-mass}--\eqref{pfeqn-poisson} describes exactly the irrotational solutions for the Euler-Poisson system.

For given constants $J_0>0$, $S_0>0$, $b_0\in(0, \rhos)$, $\rho_0\in(0, \rhos)$ and $E_0\in \R$, let $(\bar{\rho}, \bar{E})(x_1)$ be the solution to \eqref{ode-s} and \eqref{ODE_IC}. As in \eqref{equS}, we define
\begin{equation*}
  \bar{u}(x_1):=\frac{J_0}{\bar{\rho}(x_1)}.
\end{equation*}

Let $T_{\max}$ be from Lemma \ref{lemma-1dsol-new} and Definition \ref{def1}.
For each $T\in(0, T_{\max})$, there exists a constant $\eps_0\in(0, \rhos)$ depending on $(\gam, J_0, S_0, \rho_0, E_0, T)$ to satisfy
\begin{equation}
\label{1dsol-supersonicity}
 \eps_0 \le \bar{\rho}(x_1)\le \rhos-\eps_0\quad\tx{for $0\le x_1\le T$}.
\end{equation}
Let us set
\begin{equation}\label{definition-m0}
  m_0:=J_0S_0^{\frac{1}{\gam-1}}.
\end{equation}
A direct computation using \eqref{1dsol-supersonicity} shows that there exists a constant $\mu_0\in(0,1)$ depending only on $(\gam, J_0, S_0, \rho_0, E_0, T)$ to satisfy
    \begin{equation}
    \label{lwrbd-supersonicity}
      \mu_0\le \bar{u}^2(x_1)-\frac{\gam m_0^{\gam-1}}{\bar u^{\gam-1}(x_1)}\le \frac{1}{\mu_0}\quad \tx{for all $0\le x_1\le T$}.
    \end{equation}
For the given background solution $(\bar u, \bar E)$, define functions $(\vphi_0, \Phi_0)$ by
\begin{equation}
\label{definition-background-potential}
  \begin{split}
  &\vphi_0(x_1,x_2):=\int_0^{x_1} \bar{u}(t)\,dt,\\
  &\Phi_0(x_1,x_2):=\int_0^{x_1} \bar{E}(t)\,dt+\msB_0\quad\tx{with $\msB_0:=\frac 12 \left(\frac{J_0}{\rho_0}\right)^2+\frac{\gam S_0}{\gam-1}\rho_0^{\gam-1}$}.
  \end{split}
\end{equation}
Note that $\Phi_0$ is defined so that $\Phi_0-\msB_0=0$ holds on $\Gamen$. This guarantees that $\msK\equiv 0$ for the background solutions.

The functions $(\vphi_0, \Phi_0)$ satisfy the equations \eqref{pfeqn-mass} and \eqref{pfeqn-poisson} in $\Om_L$ provided that $L\in (0, T_{\max})$,
and satisfies the following boundary conditions:
\begin{equation*}
  \begin{split}
  &
  \vphi_0=0,\,\,
  \der_{x_1}\vphi_0=u_0,\,\, \der_{x_1}\Phi_0=E_0, \quad \tx{on $\Gamen$,}\\
  &\der_{x_2}\vphi_0=\der_{x_2}\Phi_0=0 \quad \tx{on $\Lambda_L$},
\end{split}
\end{equation*}
where $u_0:=\frac{J_0}{\rho_0}$.

\begin{definition}\label{defbg}
We call $(\vphi_0, \Phi_0)$ defined in \eqref{definition-background-potential}  and associated $(\bar\rho, \bar u, \bar E)$ the background solutions of the Euler-Poisson system \eqref{2-a2}.
\end{definition}

A special case of Problem \ref{problem-1} is to find supersonic irrotational solutions for the Euler-Poisson system.
\begin{problem}
  \label{problem-2-potential}
Given functions $(u_{\rm en}, E_{\rm en}, b, \Phi_{\rm ex})$ sufficiently close to $(u_0, E_0, b_0, \Phi_0(L,x_2))$, find a solution $(\vphi, \Phi)$ to the system of \eqref{pfeqn-mass}--\eqref{pfeqn-poisson} in $\Om_L=(0,L)\times (-1,1)$ with the following boundary conditions
\begin{equation}
\label{BCs-potential}
  \begin{split}
  \vphi=0,\,\,\vphi_{x_1}=u_{\rm en},\,\, \Phi_{x_1}=E_{\rm en}\quad &\tx{on $\Gamen$},\\
  \Phi=\Phi_{\rm ex}\quad &\tx{on $\Gam_{L}$},\\
  \vphi_{x_2}=\Phi_{x_2}=0\quad &\tx{on $\Lambda_{L}$}.
  \end{split}
\end{equation}
\end{problem}

Our first result gives the structural stability of irrotational supersonic solutions to \eqref{2-a2}. Namely, we prove the well-posedness of
Problem \ref{problem-2-potential} when the boundary data in \eqref{BCs-potential} are given as
small perturbations of the background solution on the associated part of the boundary.
\begin{theorem}
\label{theorem-2-potential}
For fixed constants $J_0>0$, $S_0>0$, $b_0\in(0, \rhos)$, $\rho_0\in(0, \rhos)$ and $E_0\in \R$, let $(\bar{\rho}, \bar{E})(x_1)$ and $(\vphi_0, \Phi_0)$ be the associate background solutions defined in Definition \ref{defbg} and let $T_{\max}>0$ be given in Lemma \ref{lemma-1dsol-new} and Definition \ref{def1}.

For each $\epsilon_0\in(0, T_{\rm{max}})$, there exist constants $\bar{L}\in(0, T_{\max}-\epsilon_0]$, $\bar \sigma_{\rm bd}>0$, $\bar C >0$ and $\kappa_0>0$ depending only on $(J_0, S_0, b_0, \rho_0, E_0, \epsilon_0)$
so that, whenever $L\in(0, \bar{L}]$, if given functions $u_{\rm en}\in C^3(\ol{\Gamen})$, $E_{\rm en}\in C^4(\ol{\Gamen})$, $\Phi_{\rm ex}\in C^4(\ol{\Gamex})$ and $b\in C^2(\ol{\Om_L})$
satisfy
\begin{equation}
\label{estimate-bc}
  \|u_{\rm en}-u_0\|_{C^3({\ol{\Gamen}})}+\|E_{\rm en}-E_0\|_{C^4(\ol{\Gamen})}+\|\Phi_{\rm ex}-\Phi_0(L,\cdot)\|_{C^4(\ol{\Gamex})}+\|b-b_0\|_{C^2({\ol{\Om_L}})}\le \bar\sigma_{\rm bd},
\end{equation}
and the compatibility conditions:
\begin{equation}  \label{compatibility-condition-BC1}
  \der_{x_2}b=0\,\,\tx{on $\Lambda_L$},\quad \frac{du_{\rm en}}{dx_2}(\pm 1)=\frac{d^kE_{\rm en}}{dx_2^k}(\pm 1)= \frac{d^k\Phi_{\rm ex}}{dx_2^k}(\pm 1)=0\quad\tx{for $k=1,3$ },
\end{equation}
then Problem \ref{problem-2-potential} has a unique solution $(\vphi, \Phi)\in [H^4(\Om_L)]^2$ that satisfies
\begin{equation}\label{apriori-nlbvp-potential}
\begin{split}
  &\|(\vphi-\vphi_0, \Phi-\Phi_0)\|_{H^4(\Om_L)} \le \bar C \Big(
   \|u_{\rm en}-u_0\|_{C^3({\ol{\Gamen}})}+\|E_{\rm en}-E_0\|_{C^4(\ol{\Gamen})}\\
   &\qquad \qquad +\|\Phi_{\rm ex}-\Phi_0(L,\cdot)\|_{C^4({\ol{\Gamex}})}+\|b-b_0\|_{C^2({\ol{\Om_L}})}\Big),
   \end{split}
\end{equation}
and
\begin{equation}\label{supersonicity-potential}
 c^2(\Phi, \nabla\vphi)-|\nabla\vphi|^2\le - \kappa_0 \quad\tx{in $\ol{\Om_L}$},
\end{equation}
where  $c(z, {\bf q})$ is the sound speed given by
\begin{equation}\label{defeqsound}
c(z, {\bf q})=\sqrt{(\gamma-1)\left(z-\frac{1}{2}|{\bf q}|^2\right)}\,\,.
\end{equation}
Furthermore, the solution $(\vphi, \Phi)$ satisfies
\begin{equation}
\label{comp-cond-nlbvp-potential}
  \der_{x_2}^k\vphi=\der_{x_2}^k\Phi=0\quad\tx{on $\Lambda_L$ for $k=1,3$ in the sense of trace}.
\end{equation}
Furthermore, for each $\alp\in(0,1)$, it follows from \eqref{apriori-nlbvp-potential} that
\begin{equation}\label{apriori-nlbvp-potential-holder}
\begin{split}
  &\|(\vphi-\vphi_0, \Phi-\Phi_0)\|_{C^{2,\alp}(\ol{\Om_L})}
  \le \bar C^*\Big(
   \|u_{\rm en}-u_0\|_{C^3(\ol{\Gamen})}+\|E_{\rm en}-E_0\|_{C^4(\ol{\Gamen})}\\
   &\qquad \qquad +\|\Phi_{\rm ex}-\Phi_0(L,\cdot)\|_{C^4(\ol{\Gamex})}+\|b-b_0\|_{C^2(\ol{\Om_L})}\Big),
   \end{split}
\end{equation}
for some constant $\bar C^*>0$ depending only on $(J_0, S_0, b_0, \rho_0, E_0, \epsilon_0, L, \alp)$.

\end{theorem}

When the boundary data $S_{\rm en}$ in \eqref{bc-T} is not a constant, we have the following results on the flows with non-zero vorticity.
\begin{theorem}
\label{theorem-1}
For given functions $b\in C^2(\ol{\Om_L})$, $u_{\rm en}\in C^3(\ol{\Gamen})$, $(v_{\rm en},S_{\rm en},E_{\rm en})\in C^4(\ol{\Gamen})$ and $\Phi_{\rm ex}\in C^4(\ol{\Gamex})$, let us set
\begin{equation}
\label{definition-sigma}
\begin{split}
  \sigma(b,u_{\rm en}, v_{\rm en}, E_{\rm en},\Phi_{\rm ex}, S_{\rm en}):=&\|b-b_0\|_{C^2(\ol{\Om_L})}+\|u_{\rm en}-u_0\|_{C^3(\ol{\Gamen})}+\|(v_{\rm en}, E_{\rm en}-E_0)\|_{C^4(\ol{\Gamen})}\\
  &+\|\Phi_{\rm ex}-\Phi_0(L,\cdot)\|_{C^4(\ol{\Gamex})}+\|S_{\rm en}-S_0\|_{C^4(\ol{\Gamen})}.
  \end{split}
\end{equation}
For fixed constants $J_0>0$, $S_0>0$, $b_0\in(0, \rhos)$, $\rho_0\in(0, \rhos)$ and $E_0\in \R$, let $(\bar{\rho}, \bar{E})(x_1)$ and $(\vphi_0, \Phi_0)$ be the associate background solutions defined in Definition \ref{defbg}, and let $T_{\max}>0$ be given in Lemma \ref{lemma-1dsol-new} and Definition \ref{def1}.

For each $\epsilon_0\in(0, T_{\rm{max}})$, there exist constants $\bar{L}\in(0, T_{\max}-\epsilon_0]$, $\hat \sigma_{\rm bd}>0$, $\hat C >0$ and $\hat\kappa_0>0$ depending only on $(J_0, S_0, b_0, \rho_0, E_0, \epsilon_0)$
so that, whenever $L\in(0, \bar{L}]$, if given functions $(b, u_{\rm en}, v_{\rm en}, E_{\rm en}, \Phi_{\rm ex}, S_{\rm en})$ satisfy
\begin{equation}
\label{estimate-bc-full}
 \sigma(b,u_{\rm en}, v_{\rm en}, E_{\rm en},\Phi_{\rm ex}, S_{\rm en}) \le \hat\sigma_{\rm bd},
\end{equation}
and the compatibility conditions \eqref{compatibility-condition-BC1} and
\begin{equation}  \label{compatibility-condition-full}
\begin{split}
  \frac{d^k S_{\rm en}}{dx_2^k}(\pm 1)=\frac{d^{k-1}v_{\rm en}}{dx_2^{k-1}}(\pm 1)=0\quad\tx{for $k=1,3$},
  \end{split}
\end{equation}
then the nonlinear boundary value problem \eqref{hdecomp-system}--\eqref{eqn4-hdecomp} with boundary conditions \eqref{bc-vphi}--\eqref{bc-T} has a unique solution $(\vphi, \Phi, \phi, S)\in [H^4(\Om_L)]^2\times H^5(\Om_L)\times H^4(\Om_L)$ which satisfies
\begin{equation}\label{apriori-nlbvp-full}
\begin{split}
  &\|(\vphi-\vphi_0, \Phi-\Phi_0)\|_{H^4(\Om_L)}\le \hat C \sigma(b,u_{\rm en}, v_{\rm en}, E_{\rm en},\Phi_{\rm ex}, S_{\rm en}),
  \\
  &\|\phi\|_{H^5(\Om_L)}+\|S-S_0\|_{H^4(\Om_L)}\le \hat C\|S_{\rm en}-S_0\|_{C^4(\ol{\Om_L})},
   \end{split}
\end{equation}
and
\begin{equation}\label{supersonicity-full}
 c^2(\Phi, \nabla\vphi+\nabla^{\perp}\phi)-|\nabla\vphi+\nabla^{\perp}\phi|^2\le -\hat \kappa_0\quad\tx{in $\ol{\Om_L}$}
\end{equation}
with the sound speed $c(z, {\bf q})$ given by \eqref{defeqsound}.
Furthermore, the solution $(\vphi, \Phi, \phi, S)$ satisfies
\begin{equation}
\label{comp-cond-nlbvp-full}
\der_{x_2}^k\vphi=\der_{x_2}^k\Phi=\der_{x_2}^kS =\der_{x_2}^{k+1}\phi =0\quad\tx{on $\Lambda_L$ for $k=1,3$}
\end{equation}
in the sense of trace.
Hence for each $\alp\in(0,1)$, it follows from \eqref{apriori-nlbvp-full} that
\begin{equation}\label{apriori-nlbvp-full-holder}
\begin{split}
  &\|(\vphi-\vphi_0, \Phi-\Phi_0)\|_{C^{2,\alp}(\ol{\Om_L})}
  \le \hat C^{*}\sigma(b,u_{\rm en}, v_{\rm en}, E_{\rm en},\Phi_{\rm ex}, S_{\rm en}),\\
  &\|\phi\|_{C^{3,\alp}(\ol{\Om_L})}+\|S-S_0\|_{C^{2,\alp}(\ol{\Om_L})}
  \le \hat C^{*} \|S_{\rm en}-S_0\|_{C^4(\ol{\Om_L})}
   \end{split}
\end{equation}
for some constant $\hat C^{*}>0$ depending only on $(J_0, S_0, b_0, \rho_0, E_0, \epsilon_0, L, \alp)$.


\end{theorem}
We have the following remark on Theorems \ref{theorem-2-potential} and \ref{theorem-1}.
\begin{remark}
\label{remark-Lbar}
According to Theorems \ref{theorem-2-potential} and \ref{theorem-1}, one dimensional supersonic solutions of the Euler-Poisson system are structurally stable under small perturbations of boundary data provided that the length of the nozzle $\Om_L$ is less than a critical length $\bar{L}$. The significance of these theorems is that $\bar{L}$ in general is not a small constant. As we shall see later in Proposition \ref{lemma-H1-estimate}, $\bar{L}$ is chosen so that a boundary value problem of second order hyperbolic-elliptic coupled system is well posed in $\Om_L$ whenever $L\le \bar{L}$ (See \eqref{definition-bar-L1} and \eqref{definition-bar-L2}). Furthermore, such $\bar{L}$ can be precisely computed depending on background one dimensional supersonic solutions of the Euler-Poisson system.
\end{remark}

Once we prove Theorem \ref{theorem-1}, the well-posedness of Problem \ref{problem-1} is given as a corollary of Theorem \ref{theorem-1}.

\begin{corollary}[Well-posedness of Problem \ref{problem-1}]
\label{corollary-wp-problem-1}
Under the same conditions on the given functions $(b, u_{\rm en}, v_{\rm en}, E_{\rm en}, \Phi_{\rm ex}, S_{\rm en})$ as in Theorem \ref{theorem-1}, if the condition \eqref{estimate-bc-full} holds, then Problem \ref{problem-1} has a unique solution $({\bf u}, \Phi, S)\in [H^3(\Om_L)]^2\times [H^4(\Om_L)]^2$ which satisfies
\begin{equation}\label{apriori-nlbvp-full2}
\begin{split}
&\|{\bf u}-\frac{J_0}{\bar{\rho}}{\bf e}_1\|_{H^3(\Om_L)}
 +\|\Phi-\Phi_0\|_{H^4(\Om_L)}\le C \sigma(b,u_{\rm en}, v_{\rm en}, E_{\rm en},\Phi_{\rm ex}, S_{\rm en}),
  \\
  &\|S-S_0\|_{H^4(\Om_L)}\le C\|S_{\rm en}-S_0\|_{C^4(\ol{\Om_L})},
   \end{split}
\end{equation}
for $C>0$ depending only on $(J_0, S_0, b_0, \rho_0, E_0, \epsilon_0,L )$.
\end{corollary}
A detailed proof of Corollary \ref{corollary-wp-problem-1} can be easily given by adjusting \cite[Proof of Theorem 2.7]{BDX3}, so we skip proving it in this paper.

\begin{remark}[Regularity improvement of the vorticity $\nabla\times {\bf u}$]
According to \cite[Proof of Theorem 2.7]{BDX3}, if  $({\bf u}, \Phi, S)\in [H^3(\Om_L)]^2\times [H^4(\Om_L)]^2$ solves Problem \ref{problem-1}, then, \eqref{hdecomp-1} and \eqref{hdecomp-2} yield a unique $(\vphi, \phi)\in H^4(\Om_L)\times H^5(\Om_L)$ so that $(\vphi, \Phi, \phi, S)$ solves the nonlinear boundary value problem \eqref{hdecomp-system}--\eqref{eqn4-hdecomp} with boundary conditions \eqref{bc-vphi}--\eqref{bc-T}. Furthermore, it is interesting to see from \eqref{hdecomp-2} that
\begin{equation*}
\nabla\times {\bf u}=-\Delta \phi\in H^3(\Om_L).
\end{equation*}
Therefore,  the regularity of the vorticity $\nabla\times {\bf u}$ is same as the one of ${\bf u}$. In other words, the regularity of the vorticity is improved.
\end{remark}

\section{Linearized problems for the irrotational flows}
\label{section-thm-pf-potential}
In this section, we establish the unique solvability of the boundary value problem for a second order linear hyperbolic-elliptic coupled system which corresponds to the linearized problem for Problem \ref{problem-2-potential}. The analysis for this linearized problem plays a crucial role in proving  nonlinear stability for not only irrotational flows but also the flows with non-zero vorticity.

For $(\vphi_0, \Phi_0)$ given by \eqref{definition-background-potential}, define
\begin{equation*}
  \psi:=\vphi-\vphi_0\quad\text{and}\quad  \Psi:=\Phi-\Phi_0.
\end{equation*}
Suppose that $(\vphi, \Phi)$ solves Problem \ref{problem-2-potential} and satisfies
\begin{equation}\label{supercond}
  c^2(\Phi, \nabla\vphi)-\vphi_{x_1}^2< 0\quad\tx{in $\ol{\Om_L}$}
\end{equation}
for the sound speed $c(z, {\bf q})$ given by \eqref{defeqsound}. Then, Eq. \eqref{pfeqn-mass} is equivalent to
\begin{equation}
\label{pfeqn-mass-modi}
\vphi_{x_1x_1}+2A_{12}(\Phi, \nabla\vphi)\vphi_{x_1x_2}-A_{22}(\Phi, \nabla\vphi)\vphi_{x_2x_2}+B(\Phi, \nabla\Phi, \nabla\vphi )
=0\quad\tx{in}\quad \Om_L,
\end{equation}
where $A_{12}$, $A_{22}$, and $B$ are defined by
\begin{equation}
\label{coefficients-vphi-eqn}
\begin{split}
&A_{12}(z, {\bf q})=\frac{-q_1q_2}{c^2(z, |{\bf q}|^2)-q_1^2},
\quad A_{22}(z, {\bf q})=\frac{q_2^2-c^2(z,|{\bf q}|^2)}{c^2(z, |{\bf q}|^2)-q_1^2},\\
 &B(z, {\bf p}, {\bf q} )=\frac{{\bf p}\cdot{\bf q}}{c^2(z, |{\bf q}|^2)-q_1^2}
\end{split}
\end{equation}
for $z\in \R$, ${\bf p}=(p_1, p_2)\in \R^2$, ${\bf q}=(q_1, q_2)\in \R^2$ with $z-\frac 12|{\bf q}|^2>0$.

Note that $(\varphi_0, \Phi_0)$ defined in \eqref{definition-background-potential} satisfies
\begin{equation}\label{eqvphi0}
\begin{cases}
\der_{x_1x_1}\vphi_0+\bar B=0, \\
\Delta\Phi_0=\mcl{H}_0(\Phi_0, \nabla \varphi_0)-b_0
\end{cases}
\end{equation}
where $\bar B=B(\Phi_0, \nabla \Phi_0, \nabla\vphi_0)$ and  $\mcl{H}_0$ is defined in \eqref{def-rho}.
Let
\begin{equation}
\label{def-as}
\begin{split}
&a_{j2}(\rx, z, {\bf q})=A_{j2}(\Phi_0+z, \nabla\vphi_0+{\bf q})\quad \text{for}\,\, j=1, 2,\\
& \bara_{22}(x_1) =A_{22}(\Phi_0, \nabla\vphi_0),\quad Q(\rx, z, {\bf p}, {\bf q})=B(z, {\bf p}, {\bf q})-\bar{B},\\
&(\bara_1, \bb_1, \bb_2)=(\der_{q_1}, \der_{p_1}, \der_z)B(\rx, z, {\bf p}, {\bf q})|_{(z, {\bf p}, {\bf q})=(\Phi_0, \nabla\Phi_0, \nabla \varphi_0)},
\end{split}
\end{equation}
where $(\vphi_0, \Phi_0)$ are evaluated at ${\rx}=(x_1,x_2)$ in the above.
Subtracting the first equation in \eqref{eqvphi0} from Eq. \eqref{pfeqn-mass-modi} yields
\begin{equation}\label{eqn-N1}
\mcl{N}_1(\psi, \Psi)=\mathfrak{f}_1(\rx, \Psi, \nabla\Psi, \nabla\psi)\quad\tx{in $\Om_L$},
\end{equation}
where $\mcl{N}_1(\psi, \Psi)$ and $\mathfrak{f}_1(\rx, z, {\bf p}, {\bf q})$ are defined by
\begin{equation}
\label{defN1}
\begin{aligned}
&\mcl{N}_1(\psi, \Psi):=
\psi_{x_1x_1}+2a_{12}(\rx, \Psi, \nabla\psi)\psi_{x_1x_2}
-a_{22}(\rx, \Psi, \nabla\psi)\psi_{x_2x_2}\\
&\phantom{\mcl{N}_1(\psi, \Psi):=}+\bara_1({\rx})\psi_{x_1}
+\bb_1({\rx})\Psi_{x_1}+\bb_2({\rx})\Psi
\end{aligned}
\end{equation}
and
\begin{equation}
\label{def-f1}
\mathfrak{f}_1(\rx, z, {\bf p}, {\bf q})
=-Q(\rx, \Phi_0+z, \nabla\Phi_0+{\bf p}, \nabla\vphi_0+{\bf q})
+\bara_1(x_1)q_1+\bb_1(x_1) p_1+\bb_2(x_1) z,
\end{equation}
respectively.
With the aid of \eqref{ode-s} and \eqref{equS}, the direct computations  give
\begin{equation}
\label{def-coefficients1}
\begin{split}
&\ba_{22}({\rx})=\frac{1}{\frac{\bar u^{\gam+1}}{\gam m_0^{\gam-1}}-1},
\quad
\ba_1({\rx})=\frac{\barE\left(\gam \baru^{2}+\frac{\gam m_0^{\gam-1}}{\baru^{\gam-1}}\right)}{\left(\baru^2-\frac{\gam m_0^{\gam-1}}{\baru^{\gam-1}}\right)^2},\\
&\bb_1({\rx})=\frac{\baru}{\frac{\gam m_0^{\gam-1}}{\baru^{\gam-1}}-\baru^2},\quad \bb_2({\rx})=\frac{-(\gam-1)\barE\baru}{(\baru^2-\frac{\gam m_0^{\gam-1}}{\baru^{\gam-1}})^2}
\end{split}
\end{equation}
for the constant $m_0>0$ given by \eqref{definition-m0}.

\begin{lemma}
\label{lemma-bgdcoeff}
Let $T_{\max}>0$ be given by Lemma \ref{lemma-1dsol-new} and Definition \ref{def1}. For a constant $T\in (0, T_{\max})$, we set
\begin{equation*}
  \Om_T:=\{{\rx}\in \R^2: 0<x_1<T,\,\,-1<x_2<1\}.
\end{equation*}
Then, the coefficients $(\ba_{22}, \ba_1, \bb_1, \bb_2)$ given in \eqref{def-coefficients1} satisfy the following properties:
\begin{itemize}
\item[(a)]  there exists a constant $\mu'_0\in(0,1)$ so that $\bar{a}_{22}(\rx)$ satisfies
    \begin{equation*}
     \mu'_0 \le \bar{a}_{22}(\rx)\le \frac{1}{\mu'_0}\quad\tx{for ${\rx}\in \overline{\Om_{T}}$.}
    \end{equation*}

\item[(b)]  $(\bar{a}_{22}, \bar{a}_1, \bar{b}_1, \bar{b}_2)$ are smooth in $\overline{\Om_{T}}$. More precisely, for each $k\in \mathbb{Z}^+$, there exists a constant $\bar C_k>0$ to satisfy
    \begin{equation*}
      \|(\bar{a}_{22}, \bar{a}_1, \bar{b}_1, \bar{b}_2)\|_{C^k(\overline{\Om_{T}})}\le \bar C_k.
    \end{equation*}
     \end{itemize}
In the statements above, the constant $\mu_0'$ depends only on $(\gam, J_0, S_0, \rho_0, E_0, T)$, and the constant $\bar C_k$ depends on $(\gam, J_0, S_0, \rho_0, E_0, T,k)$.
    \begin{proof}
    The statement (a) can be directly obtained from \eqref{lwrbd-supersonicity}. Lemma \ref{lemma-1dsol-new} and \eqref{lwrbd-supersonicity} yield the statement (b).
    \end{proof}
\end{lemma}
Given positive constants $T$ and $\delta_0$, define the set
\begin{equation}
\begin{aligned}
 U^T_{\delta_0}=\{({\rx}, z, {\bf q})\in &\Om_T\times\R\times \R^2:
\,\, -\delta_0<z< \delta_0,\,\,|{\bf q}|<\delta_0\ \}.
\end{aligned}
\end{equation}
The following lemma is a direct consequence of  \eqref{def-as} and Lemma \ref{lemma-bgdcoeff}.
\begin{lemma}
\label{lemma-coeff-reg}
For each $\epsilon\in(0, T_{\max})$,
there exist positive constants $(\delta_0, \kappa_0, \kappa_1)$ depending only on $(\gam, J_0, S_0, \rho_0, E_0, \epsilon)$ so that, whenever $T\in(0, T_{\max}-\epsilon]$,  the following properties hold:
\begin{itemize}
\item[(a)]
$c^2(\Phi_0(\rx)+z, \nabla\vphi_0(\rx)+{\bf q} )-(\der_{x_1}\vphi_0(\rx)+q_1)^2\le -\kappa_0$ holds in $\ol{U^T_{\delta_0}}$;

\item[(b)] $\displaystyle{a_{22}({\rx},z, {\bf q})}\ge \kappa_1$ holds in $\ol{U^T_{\delta_0}}$;

\item[(c)] for each $k\in \mathbb{Z}^+$ and $j=1,2$, there exists a positive  constant $C>0$ depending only on $(\gam, J_0, S_0, \rho_0, E_0, \epsilon, k)$ such that
    \begin{equation}
\label{est-coeff-gen-3}
|D^k_{(\rx,z, {\bf q})}a_{j2}(\rx, z, {\bf q})|\le C \quad\tx{in $\ol{U_{ \delta_0}^T}$}.
\end{equation}
\end{itemize}
\end{lemma}

Suppose that $\mcl{H}_0(\Phi, \nabla\vphi)>0$ in $\ol{\Om_L}$. Subtracting the second  equation in \eqref{eqvphi0}
from \eqref{pfeqn-poisson} yields
\begin{equation}
\label{pfeqn-poisson-modi}
{\mcl{L}}_2(\psi, \Psi):=\Delta \Psi-\bar h_1(x_1)\Psi-\bar h_2(x_1)\psi_{x_1}
={\mathfrak{f}}_2(\rx, \Psi, \nabla\psi)\quad\tx{in}\quad \Om_L,
\end{equation}
where
\begin{equation}
\label{def-coefficients3}
(\bar h_1, \bar h_2)(x_1)=(\partial_z \mcl{H}_0,\partial_{q_1}\mcl{H}_0)(z, {\bf q})|_{\small (z,{\bf q})=(\Phi_0, \nabla \varphi_0)}=
\frac{\baru^{\gam-2}}{\gam m_0^{\gam-2}S_0^{\frac{1}{\gam-1}}}\left(1,-\bar u\right),
\end{equation}
and
\begin{equation}
\label{def-f2}
\begin{split}
{\mathfrak{f}}_2(\rx, z, {\bf q})=
\mcl{H}_0(\Phi_0+z, \nabla\vphi_0+{\bf q})-\mcl{H}_0(\Phi_0, \nabla\vphi_0)
-\bar h_1z-\bar h_2q_1-(b-b_0).
\end{split}
\end{equation}
  For each $T\in(0, T_{\max})$, the constant $\mu_{1, T}$ given by
\begin{equation}
\label{lwrbd-bar-h1}
\bar\mu_{1, T}:=\inf_{x_1\in [0, T]}\bar{h}_1(x_1)
\end{equation}
is strictly positive depending on $(\gam, J_0, S_0, \rho_0, E_0, T)$.

If the conditions
\begin{equation}
\label{condition-for-equiv}
0<L \le T_{\max}-\epsilon,\quad \|\Psi\|_{C^1(\ol{\Om_L})}<\delta_0,\quad\text{and}\quad  \|\psi\|_{C^1(\ol{\Om_L})}<\delta_0
\end{equation}
are satisfied for the constant $\delta_0$ from Lemma \ref{lemma-coeff-reg}, then the nonlinear boundary value problem \eqref{pfeqn-mass}, \eqref{pfeqn-poisson}, and \eqref{BCs-potential} for $(\vphi, \Phi)$ is equivalent to the following nonlinear system
\begin{equation}
\label{system-of-perturbations}
\begin{cases}
\mcl{N}_1(\psi, \Psi)=\mathfrak{f}_1(\rx, \Psi, \nabla\Psi, \nabla\psi)\\
{\mcl{L}}_2(\psi, \Psi)=\mathfrak{f}_2(\rx, \Psi, \nabla\Psi, \nabla\psi)
\end{cases}\quad\tx{in $\Om_L$},
\end{equation}
with boundary conditions
\begin{equation}
\label{BCs-of-perturbations}
\begin{split}
(\psi, \psi_{x_1})=(0,u_{en}-u_0)=:(0,g_1),\quad \Psi_{x_1}=E_{en}-E_0=:g_2 \quad &\tx{on $\Gamen$},\\
\Psi=\Phi_{ex}-\Phi_0(L,\cdot)=:\Psi_{ex} \quad &\tx{on $\Gam_{L}$},\\
\psi_{x_2}=\Psi_{x_2}=0 \quad &\tx{on $\Lambda_L$}.
\end{split}
\end{equation}
Given constants $\delta>0$ and $L\in(0, T_{\max})$ to be fixed later, we define an iteration set:
\begin{equation}
\label{def-iter}
\begin{split}
\mcl{J}_{\delta,L}:= \Bigl\{(\psi, \Psi)\in [H^4(\Om_L)]^2: \,\, \|(\psi, \Psi)\|_{H^4(\Om_L)}\le \delta,\,\, \der_{x_2}^k\psi=\der_{x_2}^k\Psi=0\,\,\tx{ on $\Lambda_{L}$ }&\\
\tx{in the sense of trace for $k=1,3$}&\Bigr\}.
\end{split}
\end{equation}
If $b$, $u_{\rm en}$, $E_{\rm en}$, and $\Phi_{\rm ex}$ satisfy the compatibility conditions \eqref{compatibility-condition-BC1}, then
 $g_1$ and $g_2$ satisfy the  compatibility conditions
\begin{equation}
\label{comp-cond1}
\begin{split}
\frac{d g_1}{dx_2}(\pm 1)=0,\quad \frac{d^k g_2}{dx_2^k}(\pm 1)=\frac{d^k \Psi_{\rm ex}}{dx_2^k}(\pm 1)=0\,\,\tx{for $k=1$ and $3$},\quad \der_{x_2} b(\rx)=0\,\,\tx{on $\Lambda_L$}.
\end{split}
\end{equation}

For a fixed $(\tpsi, \tPsi)\in \mcl{J}_{\delta,L}$, denote
\begin{equation}
\label{definition-til-as}
\begin{split}
\til{a}_{j2}(\rx):=a_{j2}(\rx, \tPsi, \nabla\tpsi)\quad\tx{for $j=1$ and $2$},
\end{split}
\end{equation}
where $(a_{12}, a_{22})(\rx, z, {\bf q})$ are given by  \eqref{def-as}.
The following lemma is a direct consequence of \eqref{coefficients-vphi-eqn}, \eqref{def-as}, \eqref{def-iter}, and Lemma \ref{lemma-coeff-reg}.
\begin{lemma}
\label{lemma-coeff-regularity2}
Let $\delta_0>0$ be from Lemma \ref{lemma-coeff-reg}.
For each $\epsilon\in(0, T_{\max})$, there exist positive constants $\delta_1\in(0, \delta_0]$ and $C$ depending only on $(\gam, S_0, J_0, \rho_0, E_0, \epsilon)$ so that whenever $\delta\in(0, \delta_1]$ and $L\in(0, T_{\max}-\epsilon]$, the coefficients $(\til{a}_{12}, \til{a}_{22})$ defined by \eqref{definition-til-as} for $(\tpsi, \tPsi)\in \mcl{J}_{\delta, L}$ with $\mcl{J}_{\delta, L}$ given by \eqref{def-iter}, satisfy the following properties:
\begin{itemize}
\item[(a)] $\displaystyle{\|(\til{a}_{12},\til{a}_{22}-\ba_{22})\|_{H^3(\Om_L)}\le C\delta}$,
\item[(b)] $\displaystyle{\til{a}_{22}\ge \kappa_1}$ in $\ol{\Om_L}$ for the constant $\kappa_1>0$ from Lemma \ref{lemma-coeff-reg}(b),
\item[(c)] $\displaystyle{\til{a}_{12}=0}$ on $\Lambda_L$.
\end{itemize}
\end{lemma}

For a fixed $(\tpsi, \tPsi)\in \mcl{J}_{\delta,T}$, define a  linear operator $\mcl{L}_1^{(\tpsi, \tPsi)}$ by
\begin{equation}
\label{def-L1}
\mcl{L}_1^{(\tpsi, \tPsi)}(\psi, \Psi)=
\psi_{x_1x_1}+2\til{a}_{12}\psi_{x_1x_2}
-\til{a}_{22}\psi_{x_2x_2}+\bara_1\psi_{x_1}
+\bb_1\Psi_{x_1}+\bb_2\Psi
\end{equation}
with $(\bara_1, \bb_1, \bb_2)$ given in \eqref{def-as}.

For a fixed $\epsilon \in(0, T_{\max})$, let $\delta_1(\epsilon)$ represent the constant $\delta_1$ from Lemma \ref{lemma-coeff-regularity2}, and let $(g_1, g_2, \Psi_{ex})$ be given by \eqref{BCs-of-perturbations}. For fixed $L\in(0, T_{\max}-\epsilon]$ and $\delta\in(0, \delta_1(\epsilon)]$, we consider the following linear boundary value problem associated with $(\tpsi, \tPsi)\in \mcl{J}_{\delta,L}$:
\begin{equation}
\label{lbvp-iter}
\begin{split}
&\begin{cases}
\mcl{L}_1^{(\tpsi, \tPsi)}(\psi, \Psi)=\mathfrak{f}_1(\rx, \tPsi, \nabla\tPsi, \nabla\tpsi)\\
\mcl{L}_2 (\psi, \Psi)=\mathfrak{f}_2(\rx, \tPsi, \nabla\tpsi)
\end{cases}\quad\tx{in $\Om_L$},\\
&(\psi, \psi_{x_1})=(0,g_1),\quad \Psi_{x_1}=g_2
\quad \tx{on $\Gamen$},\\
&\Psi=\Psi_{ex} \quad \tx{on $\Gam_{L}$},\\
&\psi_{x_2}=\Psi_{x_2}=0 \quad \tx{on $\Lambda_{L}$}.
\end{split}
\end{equation}
To simplify the boundary conditions of $\Psi$ in \eqref{lbvp-iter}, we define a function $\Psi_{\rm bd}$ by
\begin{equation}
\label{definition-bds}
 \Psi_{\rm bd}(\rx)=(x_1-L)g_2(x_2)+\Psi_{ex}(x_2)\quad\tx{for $\rx=(x_1,x_2)\in\ol{\Om_L}$},
 \end{equation}
and set
\begin{equation}
\label{definition-hat-Psi}
 \hat{\Psi}:=\Psi-\Psi_{\rm bd}.
\end{equation}
Hence $(\psi, \Psi)$ solves \eqref{lbvp-iter} if and only if $({\psi}, \hat{\Psi})$ satisfies
\begin{equation}
\label{lbvp-iter2}
\begin{split}
&\begin{cases}
\mcl{L}_1^{(\tpsi, \tPsi)}({\psi}, \hat{\Psi})=\hat{\mathfrak{f}}_1(\rx, \tPsi, \nabla\tPsi, \nabla\tpsi)\\
\mcl{L}_2 ({\psi}, \hat{\Psi})=\hat{\mathfrak{f}}_2(\rx, \tPsi, \nabla\tpsi)
\end{cases}\quad\tx{in $\Om_L$},\\
&{\psi}=\hat{\Psi}_{x_1}=0,\quad {\psi}_{x_1}=g_1
\quad \tx{on $\Gamen$},\\
&\hat{\Psi}=0 \quad \tx{on $\Gam_{L}$},\\
&{\psi}_{x_2}=\hat{\Psi}_{x_2}=0 \quad \tx{on $\Lambda_{L}$},
\end{split}
\end{equation}
where $(\hat{\mathfrak{f}}_1, \hat{\mathfrak{f}}_2)$ is given by
\begin{equation}
\label{definition-f-hats}
  \begin{split}
  &\hat{\mathfrak{f}}_1(\rx, \tPsi, \nabla\tPsi, \nabla\tpsi)=\mathfrak{f}_1(\rx, \tPsi, \nabla\tPsi, \nabla\tpsi)-
  \mcl{L}_1^{(\tpsi, \tPsi)}(0, \Psi_{\rm bd}),\\
  &\hat{\mathfrak{f}}_2(\rx, \tPsi, \nabla\tpsi)={\mathfrak{f}}_2(\rx, \tPsi, \nabla\tpsi)-\mcl{L}_2 (0, \Psi_{\rm bd}).
  \end{split}
\end{equation}
In the following two propositions, we establish a priori estimates of $(\psi, \hat{\Psi})$, and prove the well-posedness of \eqref{lbvp-iter2}, from which the well-posedness of \eqref{lbvp-iter} easily follows.

The first proposition provides the key estimate to prove the well-posedness of \eqref{lbvp-iter2}, or equivalently the well-posedness of \eqref{lbvp-iter}.
\begin{proposition}[A priori $H^1$ estimate]
\label{lemma-H1-estimate}
For each $\epsilon_0\in(0, \frac{T_{\max}}{10}]$,   let $\delta_1(\epsilon_0)$ represent the constant $\delta_1$ from Lemma \ref{lemma-coeff-regularity2}. Then, there exists a constant $\bar{L}\in(0, T_{\max}-\epsilon_0]$ depending only on $(\gam, S_0, J_0, \rho_0, E_0, \epsilon_0)$ so that if the iteration set $\mcl{J}_{\delta, L}$ is given by \eqref{def-iter} for $L\in(0, \bar{L})$ and $\delta\in(0, \delta_1(\epsilon_0)]$, then the following statement holds:
for a fixed $(\tpsi, \tPsi)\in \mcl{J}_{\delta, L}$ , if  $(\psi, \hat{\Psi})\in [C^2(\ol{\Om_L})]^2$ solves the linear boundary value problem \eqref{lbvp-iter2} associated with $(\tpsi, \tPsi)\in \mcl{J}_{\delta, L}$, then it satisfies the estimate
    \begin{equation}
      \label{apriori-H1-for-smooth-sol}
      \begin{split}
      \|\psi\|_{H^1(\Om_L)}+\|\Psi\|_{H^1(\Om_L)}\le C\Bigl(\|\hat{\mathfrak{f}}_1\|_{L^2(\Om_L)}
      +\|\hat{\mathfrak{f}}_2\|_{L^2(\Om_L)}
           +\|g_1\|_{C^0(\ol{\Gamen})}\Bigr),
      \end{split}
    \end{equation}
    for $\hat{\mfrak{f}}_1=\mathfrak{f}_1(\rx, \tPsi, \nabla\tPsi, \nabla\tpsi)$ and $\hat{\mfrak{f}}_2=\mathfrak{f}_2(\rx, \tPsi, \nabla\tpsi)$. Here, the constant $C$ depends only on $\gam$, $J_0$, $S_0$, $\rho_0$, $E_0$, $\epsilon_0$,  and $L$.
\end{proposition}

\begin{proof}
The proof is divided into three steps.

{\textbf{Step 1. Weighted energy equality.}} Fix a small constant $\epsilon_0\in(0, \frac{T_{\max}}{10}]$ and fix $(\tpsi, \tPsi)\in \mcl{J}_{\delta, L}$, and suppose that $(\psi, \hat{\Psi})\in [C^2(\ol{\Om_L})]^2$ solve \eqref{lbvp-iter2}.

By Lemma \ref{lemma-coeff-regularity2}(b), if
  $\delta\in(0, \delta_1(\epsilon_0)]$  and $L\in(0, T_{\max}-\epsilon_0]$, then
the coefficient $\til{a}_{22}$ satisfies
\begin{equation}
\label{lwrbd-til-a22}
  \til{a}_{22}\ge \kappa_1(\epsilon_0)\quad\tx{in $\ol{\Om_L}$},
\end{equation}
where the positive constants $(\delta_1, \kappa_1)(\epsilon_0)$ are chosen depending on $(\gam, S_0, J_0, \rho_0, E_0, \epsilon_0)$.

For a smooth function $\mcW(x_1)$ to be specified later, let us define a functional $\mcl{I}_L({\psi}, \hat{\Psi}, \mcW)$ as follows
\begin{equation}
\label{definition-IL}
 \mcl{I}_L({\psi}, \hat{\Psi}, \mcW):= \int_{\Om_L} \mcW(x_1){\psi}_{x_1} \mcl{L}_1^{(\tpsi, \tPsi)}({\psi}, \hat{\Psi})-\hat{\Psi} \mcl{L}_2 ({\psi}, \hat{\Psi})\,d\rx.
\end{equation}
Since $(\psi, \hat{\Psi})$ solves \eqref{lbvp-iter2}, one has
\begin{equation}
\label{equation-for-IL}
   \mcl{I}_L({\psi}, \hat{\Psi}, \mcW)=
   \int_{\Om_L} \mcW(x_1){\psi}_{x_1}\hat{\mathfrak{f}}_1(\rx, \tPsi, \nabla\tPsi, \nabla\tpsi)-\hat{\Psi}\hat{\mathfrak{f}}_2(\rx, \tPsi, \nabla\tpsi)\,d\rx.
\end{equation}
Integrating by parts gives
\begin{equation}
\label{IL-1}
  \mcl{I}_L(\psi, \hat{\Psi}, \mcW)
  =J_1(\psi, \hat{\Psi}, \mcW)+J_2(\hpsi, \hat{\Psi}, \mcW)
  +J_3(\psi, \mcW)
\end{equation}
with
\begin{equation*}
  \begin{split}
  &J_1(\psi, \hat{\Psi}, \mcW)=\int_{\Om_L}
  q_1(\mcW,\rx)\frac{\psi_{x_1}^2}{2}+q_2(\mcW, \rx)\frac{\psi_{x_2}^2}{2}+|\nabla\hat{\Psi}|^2
  +\bar{h}_1\hat{\Psi}^2\,d\rx,\\
  &J_2(\psi, \hat{\Psi}, \mcW)=\int_{\Om_L} \mcW\psi_{x_1}(\bb_1\hat{\Psi}_{x_1}+\bb_2\hat{\Psi})
  +\bar{h}_2\psi_{x_1}\hat{\Psi} -\partial_{x_2} \til{a}_{22} \mcW\psi_{x_1}\psi_{x_2}\,d\rx,\\
  &J_3(\psi, \mcW)=\int_{\Gam_{L}}\frac{\mcW}{2}(\psi_{x_1}^2
  +\til{a}_{22}\psi_{x_2}^2)\,dx_2-\int_{\Gamen} \frac{\mcW}{2}g_1^2\,dx_2,
  \end{split}
\end{equation*}
where $(q_1, q_2)(\mcW,\rx)$ is defined by
\begin{equation}
\label{definition-q1-q2}
q_1(\mcW, \rx)=-\mcW_{x_1}+2(\ba_1-\der_{x_2}\til{a}_{12})\mcW\quad\text{and}\quad
q_2(\mcW, \rx)=-\til{a}_{22}\left(\mcW_{x_1}+\frac{\der_{x_1} \til{a}_{22}}{\til{a}_{22}}\mcW\right).
\end{equation}

{\textbf{Step 2. A priori estimate of $(\psi, \hat{\Psi})$.}}
It follows from Cauchy-Schwarz inequality and Lemma \ref{lemma-coeff-regularity2}(a) that
\begin{equation}
\label{IL-2}
  |J_2(\psi, \hat{\Psi}, \mcW)|
  \le \frac 18 \hat{\Psi}_{x_1}^2+\frac{\bar\mu_{1,L}}{4} \hat{\Psi}^2+q_3(\mcW,\rx)\frac{\psi_{x_1}^2}{2} +C_*\delta \psi_{x_2}^2,
\end{equation}
for a constant $C_*>0$ chosen depending only on $(\gam, S_0, J_0, \rho_0, E_0, \epsilon_0)$,
where the constant $\bar\mu_{1,L}>0$ is given by \eqref{lwrbd-bar-h1} and
\begin{equation}
\label{definition-q3}
  q_3(\mcW,\rx)=
2\left(\Bigl(\bb_1^2
  +\frac{\bb_2^2}{\bar\mu_{1,L}}+\delta\Bigr)\mcW^2
  +\frac{\bar{h}_2^2}{{\bar\mu_{1,L}}}\right).
\end{equation}
By using \eqref{lwrbd-til-a22} and \eqref{IL-2}, we can estimate $\mcl{I}_L(\psi, \hat{\Psi}, \mcW)$ given by \eqref{IL-1} as
\begin{equation}
\label{lwrbd-IL}
\begin{split}
  \mcl{I}_L(\psi, \hat{\Psi}, \mcW)\ge&
  \int_{\Om_L}(q_1(\mcW,\rx)-q_3(W,\rx))\frac{\psi_{x_1}^2}{2}+(q_2(
  \mcW, \rx)-2C_*\delta)\frac{\psi_{x_2}^2}{2}\,d\rx\\
  &+\int_{\Om_L}\frac 78|\nabla \hat{\Psi}|^2+\frac{3\bar\mu_{1,L}}{4}\hat{\Psi}^2\,d\rx\\
  &+\min\{1,\kappa_1(\epsilon_0)\}
  \int_{\Gam_{L}}\frac{\mcW}{2}(\psi_{x_1}^2+\psi_{x_2}^2)\,dx_2-\int_{\Gamen} \frac{\mcW}{2}g_1^2\,dx_2
  \end{split}
\end{equation}
provided that the function $\mcW$ satisfies
\begin{equation}
\label{W-at-exit}
  \mcW(L)\ge 0.
\end{equation}
Furthermore, if $(\mcW, L, \delta)$ are chosen to satisfy
\begin{equation}\label{W1}
\mcW(x_1)>0\quad  \text{for} \,\,0\le x_1\le L,
\end{equation}
\begin{equation}\label{W2}
 q_1(\mcW, \rx)-q_3(\mcW, \rx)\ge \lambda_0\quad  \text{in} \,\,\ol{\Om_L},
 \end{equation}
 and
 \begin{equation}\label{W3}
 \frac{1}{\til{a}_{22}}\left(q_2(\mcW, \rx)-2C_*\delta\right) \ge \lambda_0\quad \text{in}\,\, \ol{\Om_L}
 \end{equation}
for some constant $\lambda_0>0$, then we obtain from \eqref{lwrbd-til-a22} and \eqref{lwrbd-IL} that
\begin{equation}
\label{coerceivity}
\begin{split}
&\mcl{I}_L(\psi, {\hPsi},\mcW)\\
\ge \,\,& \frac{\lambda_0}{2}\int_{\Om_L}\psi_{x_1}^2+\kappa_1(\epsilon_0)\psi_{x_2}^2\,d\rx
+\frac 78\int_{\Om_L}|\nabla\hat{\Psi}|^2\,d\rx+\frac{3\bar\mu_{1,L}}{4} \int_{\Om_L}\hat{\Psi}^2\,d\rx-\int_{\Gamen} \frac{\mcW}{2}g_1^2\,dx_2.
\end{split}
\end{equation}
Applying Cauchy-Schwarz inequality to the right-hand side of \eqref{equation-for-IL} yields
\begin{equation*}
  \mcl{I}_L(\psi, {\hPsi},\mcW)\le
  \int_{\Om_L} \frac{\lambda_0}{4}\psi_{x_1}^2+\frac{\bar{\mu}_{1,L}}{4}\hat{\Psi}^2
  +\frac{1}{\lambda_0}|\mcW\hat{\mathfrak{f}}_1|^2
  +\frac{1}{\bar{\mu}_{1,L}}|\hat{\mathfrak{f}}_2|^2
  \,d\rx.
\end{equation*}
This estimate, together with \eqref{coerceivity}, gives
\begin{equation*}
  \begin{split}
  &\frac{\lambda_0}{4}\min\{1,\kappa_1(\epsilon_0)\}\int_{\Om_L}|\nabla\psi|^2\,d\rx
  +\frac 78 \int_{\Om_L}|\nabla\hat{\Psi}|^2\,d\rx
  +\frac{\bar\mu_{1,L}}{2}\int_{\Om_L}\hat{\Psi}^2\,d\rx \\
  \le \,\, & \frac{1}{\lambda_0}\int_{\Om_L}|\mcW\hat{\mathfrak{f}}_1|^2\,d\rx
  +\frac{1}{\bar\mu_{1,L}}\int_{\Om_L} |\hat{\mathfrak{f}}_2|^2\,d\rx+\int_{\Gamen} \frac{\mcW}{2}g_1^2\,dx_2.
  \end{split}
\end{equation*}
Since $\psi=0$ on $\Gamen$, we apply Poincar\'{e} inequality to derive from the estimate right above that
\begin{equation}\label{H1-estimate-1}
  \|\psi\|_{H^1(\Om_L)}+\|\hat{\Psi}\|_{H^1(\Om_L)}
  \le \mathfrak{C}_*\left (\|\hat{\mathfrak{f}}_1\|_{L^2(\Om_L)}+\|\hat{\mathfrak{f}}_2\|_{L^2(\Om_L)}
  +\|g_1\|_{C^0(\ol{\Gamen})}\right)
\end{equation}
for a constant $\mathfrak{C}_*>0$ chosen depending on $(\gam, J_0, S_0, \rho_0, E_0, \epsilon_0, L)$ and $(\lambda_0,\max_{x_1\in[0,L]}\mcW(x_1))$.
This proves Proposition \ref{lemma-H1-estimate} once we find $(\mcW,L)$ so that they satisfy the conditions \eqref{W1}--\eqref{W3}, where $(\lambda_0,\max_{x_1\in[0,L]}\mcW(x_1))$ depend only on $(\gam, J_0, S_0, \rho_0, E_0, \epsilon_0, L)$ when $L$ is chosen appropriately.

{\textbf{Step 3. Construction of $\mcW$.}}
To complete the proof, it remains to find $(\mcW, L)$ to satisfy the conditions \eqref{W1}--\eqref{W3}.

Let us set
\begin{equation*}
  \bar{\mu}_1:=\bar{\mu}_{1, T_*}\quad\tx{for $T_*=T_{\max}-\epsilon_0$}
\end{equation*}
for $\bar\mu_{1,T}$ given by \eqref{lwrbd-bar-h1}.
Then, we have
\begin{equation}
\label{inequality-mu1}
 \bar{\mu}_{1,L}\ge  \bar{\mu}_1>0 \,\,\tx{for any $L\in(0,T_{\max}-\epsilon_0]$.}
\end{equation}
By \eqref{lwrbd-til-a22}, \eqref{definition-q1-q2}, \eqref{definition-q3} and \eqref{inequality-mu1}, for any $\delta\in (0, \delta_1(\epsilon_0)]$, $L\in(0, T_{\max}-\epsilon_0]$ and ${\rx}\in \ol{\Om_L}$, we have
\begin{equation}
\label{inequality-Ws}
  \begin{split}
  q_1(\mcW,\rx)-q_3(\mcW,\rx) &\ge -\mcW_{x_1}+2(\ba_1-\der_{x_2}\til{a}_{12})\mcW-
  2\Bigl(\bb_1^2
  +\frac{\bb_2^2}{\bar\mu_{1}}+\delta\Bigr)\mcW^2
  -\frac{2\bar{h}_2^2}{{\bar\mu_{1}}},\\
  \frac{1}{\til{a}_{22}}\left(q_2(\mcW, \rx)-2C_*\delta\right)
  &\ge -\mcW_{x_1}-2\frac{\der_{x_1} \til{a}_{22}}{2\til{a}_{22}}\mcW-\frac{2C_*\delta}{\kappa_1(\epsilon_0)}.
  \end{split}
\end{equation}
Set
\begin{equation}
\label{definition-coefficnets-mfrakas-bgd}
  \begin{split}
  &\mathfrak{a}^{\epsilon_0}_0:=\frac{2C_*\delta_1(\epsilon_0)}{\kappa_1(\epsilon_0)}+\max_{[0, T_{\max}-\epsilon_0]}\frac{2\bar{h}_2^2}{\bar\mu_{1}}, \quad \mathfrak{a}^{\epsilon_0}_2:=
  \max_{[0, T_{\max}-\epsilon_0]}2\left(\bb_1^2
  +\frac{\bb_2^2}{\bar\mu_{1}}+\delta_1(\epsilon_0)\right),\\
  &\mathfrak{a}^{\epsilon_0}_{1,1}:=\min_{[0, T_{\max}-\epsilon_0]}\ba_1,\quad
  \mathfrak{a}^{\epsilon_0}_{1,2}:=\min_{[0, T_{\max}-\epsilon_0]}\left(\frac{-\der_{x_1}\bar{a}_{22}}{2\bar{a}_{22}}\right),\quad
  \mathfrak{a}^{\epsilon_0}_1:=\min\{\mathfrak{a}^{\epsilon_0}_{1,1}, \mathfrak{a}^{\epsilon_0}_{1,2}\}.
  \end{split}
\end{equation}
By Lemma \ref{lemma-bgdcoeff} and the positivity of $\bar\mu_{1}$ depending only on $(\gam, J_0, S_0, \rho_0, E_0, \epsilon_0)$, there exists a constant $C>0$ depending only on $(\gam, J_0, S_0, \rho_0, E_0, \epsilon_0)$ to satisfy
\begin{equation}
\label{positivity-a0-a2}
  0<\mathfrak{a}^{\epsilon_0}_0,\,\,\mathfrak{a}^{\epsilon_0}_2\le C,\quad\tx{and}\quad |\mathfrak{a}^{\epsilon_0}_1|\le C.
\end{equation}
For the rest of the proof, all estimate constants are chosen depending only on $(\gam, J_0, S_0, \rho_0, E_0, \epsilon_0)$ unless otherwise specified.

For each $(\tpsi, \tPsi)\in \mcl{J}_{\delta, L}$ with $(\delta,L)\in (0, \delta_1(\epsilon_0)]\times (0, T_{\max}-\epsilon_0]$, let us set
\begin{equation}
\label{definition-coefficnets-mfrakas-iterset}
  \til{\mathfrak{a}}_{1,1}:=\min_{\overline{\Omega_L}}(\ba_1-\der_{x_2}\til{a}_{12}),\quad
  \til{\mathfrak{a}}_{1,2}:=\min_{\overline{\Omega_L}}\left(\frac{-\der_{x_1}\til{a}_{22}}{2\til{a}_{22}}\right),\quad
  \til{\mathfrak{a}}_1:=\min\{\til{\mathfrak{a}}_{1,1}, \til{\mathfrak{a}}_{1,2}\}.
\end{equation}
By Lemma \ref{lemma-coeff-regularity2} and Morrey's inequality, there exists a constant $C^{\flat}>0$ satisfying
\begin{equation}
\label{lwrbd-til-a1}
  \til{\mathfrak{a}}_1\ge \mathfrak{a}_1^{\epsilon_0}-C^{\flat}\delta_1(\epsilon_0)=:
  \hat{\mathfrak{a}}_1^{\epsilon_0}.
\end{equation}
Namely, $\hat{\mathfrak{a}}_1^{\epsilon_0}$ is a constant fixed depending only on $(\gam, J_0, S_0, \rho_0, E_0, \epsilon_0)$.
To simplify notations, let $(\mathfrak{a}_0, \mathfrak{a}_1, \mathfrak{a}_2)$ denote $(\mathfrak{a}_0^{\epsilon_0}, \hat{\mathfrak{a}}_1^{\epsilon_0}, \mathfrak{a}_2^{\epsilon_0})$
for the rest of the proof.

It follows from \eqref{inequality-Ws} and \eqref{positivity-a0-a2}--\eqref{lwrbd-til-a1} that $\mcW$ satisfies all the conditions \eqref{W1}--\eqref{W3} as long as it satisfies \eqref{W1} and solves the ODE
\begin{equation}\label{Weq1}
  -\mcW_{x_1}-\mathfrak{a}_2\mcW^2+2\mathfrak{a}_1\mcW
  -\mathfrak{a}_0= \lambda_0\quad\tx{for $0\le x_1\le L$}
\end{equation}
for some constant $\lambda_0>0$.

Since $\mathfrak{a}_2>0$, Eq. \eqref{Weq1} is equivalent to
\begin{equation}
\label{ODE-for-W}
  -\frac{\mcW_{x_1}}{\mathfrak{a}_2}
  =\left(\mcW-\frac{\mathfrak{a}_1}{\mathfrak{a}_2}\right)^2
  +\frac{1}{\mathfrak{a}_2}
  \left(\frac{\mathfrak{a}_0\mathfrak{a}_2-\mathfrak{a}_1^2}{\mathfrak{a}_2}
  +\lambda_0\right).
\end{equation}
Let us set
\begin{equation}\label{defY}
  Y(x_1):=\mcW(x_1)-\frac{\mathfrak{a}_1}{\mathfrak{a}_2},
\end{equation}
to rewrite \eqref{ODE-for-W} as
\begin{equation}
\label{ODE-for-Y}
  -\frac{Y'}{\mathfrak{a}_2}=Y^2
  +\frac{1}{\mathfrak{a}_2}
  \left(\frac{\mathfrak{a}_0\mathfrak{a}_2-\mathfrak{a}_1^2}{\mathfrak{a}_2}
  +\lambda_0\right).
\end{equation}
Note that $\lambda_0>0$ is a free parameter. We solve $Y$ by considering two cases
separately according to the sign of $\mathfrak{a}_0\mathfrak{a}_2-\mathfrak{a}_1^2$.

\underline{\emph{Case 1.  $\mathfrak{a}_0\mathfrak{a}_2-\mathfrak{a}_1^2>0$}}. Set
\begin{equation}
\label{definition-nu-lambda0}
  \nu(\lambda_0):=\sqrt{\frac{1}{\mathfrak{a}_2}
  \left(\frac{\mathfrak{a}_0\mathfrak{a}_2-\mathfrak{a}_1^2}{\mathfrak{a}_2}
  +\lambda_0\right)}.
\end{equation}
Regarding $\nu(\lambda_0)$ as a function of $\lambda_0>0$, a general solution to \eqref{ODE-for-Y} for $\lambda_0>0$ is given by
\begin{equation}
\label{general-Y}
  Y(x_1;\lambda_0)=\nu(\lambda_0)\tan(C-\mathfrak{a}_2\nu(\lambda_0)x_1)
\end{equation}
for some constant $C\in \R$.
Then it follows from \eqref{defY} that a general solution $\mcW(\cdot, \lambda_0)$ to \eqref{Weq1} for $\lambda_0>0$ is given by
\begin{equation}
\label{definition-W-case1}
  \mcW(x_1;\lambda_0):=Y(x_1;\lambda_0)+\frac{\mathfrak{a}_1}{\mathfrak{a}_2}
  =\nu(\lambda_0)\tan(C-\mathfrak{a}_2\nu(\lambda_0)x_1)+\frac{\mathfrak{a}_1}{\mathfrak{a}_2}.
\end{equation}
First, we choose $C$ as
\begin{equation}
\label{choice-C-for-W}
C=\frac{\pi}{2}-\lambda_0,
\end{equation}
so that we have
\begin{equation*}
\mcW(0;\lambda_0)=\frac{\mathfrak{a}_1}{\mathfrak{a}_2}+\nu(\lambda_0)\cot \lambda_0> \frac{\mathfrak{a}_1}{\mathfrak{a}_2}
+\sqrt{\frac{\lambda_0}{\mathfrak{a}_2}}\cot \lambda_0.
\end{equation*}

 It follows from the positivity of $\mathfrak{a}_2$ in \eqref{positivity-a0-a2} and L'H\^{o}pital's rule that
 \begin{equation}\label{limit1}
\lim_{\lambda_0\to 0+} \sqrt{\frac{\lambda_0}{\mathfrak{a}_2}}\cot \lambda_0 =\infty.
\end{equation}
Hence one can fix a small constant $\lambda_1^*\in(0, \frac{\pi}{4})$ depending on $(\gam, J_0, S_0, \rho_0, E_0, \epsilon_0)$ so that whenever $\lambda_0\in(0,\lambda^*_1]$, we have $\mcW(0;\lambda_0)>0$. Since $\mcW$ solves \eqref{Weq1}, it follows from the observations made in \eqref{inequality-Ws} that if $\mcW(x_1;\lambda_0)>0$ holds for all $x_1\in[0, L]$, then $\mcW$ satisfies all the conditions \eqref{W1}--\eqref{W3} so that the energy estimate \eqref{H1-estimate-1} is validated.

By \eqref{definition-nu-lambda0}--\eqref{choice-C-for-W}, $\mcW(x_1;\lambda_0)>0$ holds for all $x_1\in[0, L]$, if $(\lambda_0, L)\in (0, \lambda^*_1]\times (0, T_{\max}-\epsilon_0]$ satisfies
\begin{equation}
\label{definition-L1-ld}
  L\le
 \frac{1}{\mathfrak{a}_2\nu(\lambda_0)}\left(
 \frac{\pi}{2}-\lambda_0-
 \arctan\left(\frac{|\mathfrak{a}_1|}{\mathfrak{a}_2\nu(\lambda_0)}\right)
 \right).
\end{equation}
Let us set
\begin{equation}
\label{definition-bar-L1}
  \bar{L}:=\min\left\{T_{\max}-\epsilon_0, \sup_{\lambda\in(0, \lambda^*]} \frac{1}{\mathfrak{a}_2\nu(\lambda_0)}\left(
 \frac{\pi}{2}-\lambda_0-
 \arctan\left(\frac{|\mathfrak{a}_1|}{\mathfrak{a}_2\nu(\lambda_0)}\right)
 \right)\right\}.
\end{equation}
Then, for any $L\in(0, \bar{L})$, one can find a constant $\lambda_0\in(0, \lambda^*_1]$ depending on $(\gam, J_0, S_0, \rho_0, E_0, \epsilon_0, L)$ to satisfy \eqref{definition-L1-ld} so that $W(x_1;\lambda_0)$ given by \eqref{definition-W-case1} satisfy \eqref{W1}--\eqref{W3} in the case of $\mathfrak{a}_0\mathfrak{a}_2-\mathfrak{a}_1^2>0$.

According to \eqref{definition-bar-L1}, the choice of $\bar{L}$ depends on $(\gam, J_0, S_0, \rho_0, E_0, \epsilon_0)$,
and the choice of $\lambda_0>0$ depends on $(\gam, J_0, S_0, \rho_0, E_0, \epsilon_0, L)$. Therefore, $\max_{x_1\in[0,L]}\mcW(x_1)$ depends on
$(\gam, J_0, S_0, \rho_0, E_0, \epsilon_0, L)$ as well.

\underline{\emph{Case 2.  $\mathfrak{a}_0\mathfrak{a}_2-\mathfrak{a}_1^2\le 0$}}.
We first denote
\begin{equation*}
  \mathfrak{a}_*:={\frac{\mathfrak{a}_1^2-\mathfrak{a}_0\mathfrak{a}_2}{\mathfrak{a}_2}}\ge 0
  \quad \tx{and} \quad \beta(\lambda_0):={\frac{\lambda_0-\mathfrak{a}_*}{\mathfrak{a}_2}}.
\end{equation*}
We choose $\lambda_0>\mathfrak{a}_*$ so that $\beta(\lambda_0)>0$, and Eq. \eqref{ODE-for-Y} can be rewritten as
\begin{equation*}
  -\frac{Y'}{\mathfrak{a}_2}=Y^2+\sqrt{\beta(\lambda_0)}^2.
\end{equation*}
Therefore,  the general solution is given by \eqref{general-Y} except that we replace $\nu(\lambda_0)$ by $\sqrt{\beta(\lambda_0)}$. We choose $C$
as $C=\frac{\pi}{2}-\beta(\lambda_0)$ to obtain from \eqref{defY} that
\begin{equation}
\label{definition-W-case2}
  \mcW(x_1;\lambda_0)=\sqrt{\beta(\lambda_0)}\cot\left(\beta(\lambda_0)
  +\mathfrak{a}_2\sqrt{\beta(\lambda_0)}x_1 \right)+\frac{\mathfrak{a}_1}{\mathfrak{a}_2}.
\end{equation}
With the aid of L'H\^{o}pital's rule, the direct computations  yield that
\begin{equation*}
  \lim_{\lambda_0\rightarrow \mathfrak{a}_*+}\mcW(0, \lambda_0)=\infty,\quad\tx{and}
  \quad \lim_{\lambda_0\rightarrow \mathfrak{a}_*+}\mcW(L, \lambda_0)=\frac{1}{\mathfrak{a}_2}\left(\frac{1}{L}+\mathfrak{a}_1\right).
\end{equation*}
Therefore, if
\begin{equation}
\label{definition-bar-L2}
\frac 1L+\mathfrak{a}_1>0\quad\tx{and} \quad L<T_{\max}-\epsilon_0,
\end{equation}
then
there exists a constant $\lambda_0>\mathfrak{a}_*$ depending only on $(\gam, J_0, S_0, \rho_0, E_0, \epsilon_0,L)$ so that $\mcW(x_1;\lambda_0)$ satisfies
  \begin{equation*}
  \mcW(x_1;\lambda_0)>0\quad \tx{for all $x_1\in[0,L]$.}
  \end{equation*}
Furthermore, it directly follows from the choice of $\lambda_0$ and \eqref{definition-W-case2} that  $\max_{x_1\in[0,L]}\mcW(x_1)$ depends on $(\gam, J_0, S_0, \rho_0, E_0, \epsilon_0, L)$. This completes the proof of  Proposition \ref{lemma-H1-estimate}.
\end{proof}

\begin{remark}[Accelerating supersonic flow]
\label{remark-H1-lemma-accel}
In \eqref{ODE_IC}, assume that $E_0>0$, and let $T_*\in(0, T_{\max})$ be given by Definition \ref{def1}(ii).
A direct computation using \eqref{ode-s}, \eqref{equS} and \eqref{def-coefficients1} shows that
  $\til{\mathfrak a}_1$ given by \eqref{definition-coefficnets-mfrakas-iterset} satisfy
\begin{equation}
\label{positivity-tila}
\til{\mathfrak a}_1\ge \frac{\bar{\mathfrak{a}}_1^*(\eps_0)}{2}>0
\end{equation}
whenever $L\in(0, T_{\max}-\eps_0]$ and $\delta\in(0, \delta^*]$. The estimate \eqref{H1-estimate-1} holds for any $L$ satisfying
\begin{equation}
\label{relaxed-relation-L}
  L <\min\left\{T_*-\eps_0, \sup_{\lambda\in(0, \bar{\lambda}]}\frac{1}{\bar{\mathfrak{a}}_2^*(\eps_0)\nu(\lambda_0)}  \left(\frac{\pi}{2}-\lambda_0+\arctan \left(\frac{\bar{\mathfrak{a}}_1^*(\eps_0)}{2\bar{\mathfrak{a}}_2^*(\eps_0)\nu(\lambda_0)}\right)\right)  \right\},
\end{equation}
and
\begin{equation}
\label{relaxed-relation-L2}
0<L< T_{\max}-\eps_0
\end{equation}
when the coefficients satisfy
$(\bar{\mathfrak{a}}^{*}_0 \bar{\mathfrak{a}}^{*}_2)(\eps_0)-(\til{ \mathfrak{a}}_1)^2 >0$ and $(\bar{\mathfrak{a}}^{*}_0 \bar{\mathfrak{a}}^{*}_2)(\eps_0)-(\til{ \mathfrak{a}}_1)^2  \le 0$, respectively.
Therefore, \eqref{relaxed-relation-L2} is more relaxed compared to \eqref{definition-bar-L2} for accelerating supersonic flow.
From this perspective,  accelerating supersonic flow of the Euler-Poisson system is relatively more stable.
\end{remark}

With the aid of Proposition \ref{lemma-H1-estimate}, we next prove the well-posedness of  the linear boundary value problem \eqref{lbvp-iter2}.
\begin{proposition}
\label{proposition-wp-lbvp}
 Fix $\epsilon_0\in(0, T_{\max})$, and let $\bar{L}\in(0, T_{\max}-\epsilon_0]$ be from Proposition \ref{lemma-H1-estimate}. For $L\in(0, \bar{L}]$ and $\delta\in(0, \frac{\delta_1(\epsilon_0)}{2}]$, let the iteration set $\mcl{J}_{\delta, L}$ be given by \eqref{def-iter}. Then, for each $(\tpsi, \tPsi)\in \mcl{J}_{\delta, L}$, the associated linear boundary value problem \eqref{lbvp-iter2} has a unique solution $(\psi, \hat{\Psi})\in [H^4(\Om_L)]^2$ that satisfies
    \begin{equation}
    \label{estimate-sol-lbvp-H4}
      \begin{split}
      &\|\hat{\Psi}\|_{H^4(\Om_L)}\le
  C\left(\|\hat{\mathfrak{f}}_1\|_{H^2(\Om_L)}+\|\hat{\mathfrak{f}}_2\|_{H^2(\Om_L)}
+\|g_1\|_{C^2(\ol{\Gamen})} \right),\\
&\|\psi\|_{H^4(\Om_L)}\le C\left(\|\hat{\mathfrak{f}}_1\|_{H^3(\Om_L)}+
\|\hat{\mathfrak{f}}_2\|_{H^2(\Om_L)}
+\|g_1\|_{C^3(\ol{\Gamen})} \right)
      \end{split}
    \end{equation}
 for a constant $C>0$ depending only on $(\gam, J_0, S_0, \rho_0, E_0, \epsilon_0, L)$. Furthermore, the solution $(\psi, \hat{\Psi})$ satisfies the compatibility conditions
 \begin{equation}
 \label{comp-cond-sol-lbvp}
   \der_{x_2}^k\psi=\der_{x_2}^k\hat{\Psi}=0\;
\tx{for $k=1,3$ on $\Lambda_{L}$ in the sense of trace}.
 \end{equation}

\begin{proof}
{\textbf{Step 1. Approximation of \eqref{lbvp-iter2} by a problem with smooth coefficients.}}

{\emph{Claim 1. For any given constant $\eps>0$, and any given function $\phi \in H^4(\Om_L)$ satisfying
\begin{equation*}
\label{compatibility-cond-fr-smoothapx}
  \der_{x_2}\phi=\der_{x_2}^3\phi=0\quad\tx{on $\Lambda_{L}$ in the sense of trace},
\end{equation*}
there exists a function $\phi^{\eps}\in C^{\infty}(\ol{\Om_L})$ satisfying
\begin{equation*}
\label{smallness-aprx}
\|\phi-\phi^{\eps}\|_{H^4(\Om_L)}\le \eps\quad\tx{and}\quad
 \der_{x_2}\phi^{\eps}=\der_{x_2}^3\phi^{\eps}=0\quad\tx{on $\Lambda_{L}$}.
\end{equation*}
}}
Given $\phi\in H^4(\Om_L)$ satisfying
\eqref{compatibility-cond-fr-smoothapx}, denote its extension $\tilde{\phi}$ onto $\widetilde{\Om}_L:=(0,L)\times (-\frac 32, \frac 32)$ by
\begin{equation*}
  \tilde{\phi}(x_1,x_2)=\begin{cases}
  \phi(x_1,-x_2+2)\quad&\mbox{for $x_2>1$},\\
  \phi(x_1,x_2)\quad&\mbox{for $-1\le x_2\le 1$},\\
  \phi(x_1, -x_2-2)\quad&\mbox{for $x_2<-1$}.
  \end{cases}
\end{equation*}
It is easy to check that $\tilde{\phi}\in H^4(\widetilde{\Om}_L)$.  We define $\phi^{\eps}$ by the convolution of $\tilde{\phi}$ with a standard mollifier $\eta_{r}(\rx)$(which depends only on $|\rx|$) with a compact support in a disk of radius $r>0$. One can choose $r>0$ sufficiently small depending on $\eps$ so that $\phi^{\eps}:=\til{\phi}*\eta_r$ satisfies $\|\phi-\phi^{\eps}\|_{H^4(\Om_L)}\le \eps$. Furthermore, it can be directly checked from the definition that $\phi^{\eps}$ satisfies
$\der_{x_2}\phi^{\eps}=\der_{x_2}^3\phi^{\eps}=0$ on $\Lambda_{L}$.  This is because $\til{\phi}$ is an even function about $x_2=\pm 1$, and $\eta_r$ is radially symmetric. This verifies the claim.

By {\emph{Claim 1}}, we can introduce a sequence $\{(\tpsi_n, \tPsi_n)\}_{n=1}^\infty \subset C^{\infty}(\ol{\Om_L})$ such that
\begin{equation}
\label{distance-approx-origina}
\|\tpsi_n-\tpsi\|_{H^4(\Om_L)}+\|\tPsi_n-\tPsi\|_{H^4(\Om_L)}\le \frac{\delta}{2n},\quad
\der_{x_2}^k\tpsi_n=\der_{x_2}^k\tPsi_n=0\quad\tx{on  $\Lambda_L$ for $k=1,3$}.
\end{equation}

For each $n\in \mathbb{N}$, define
\begin{equation}
\label{smooth-approx-in-lbvp}
  \begin{split}
  &{a}_{j2}^{(n)}(\rx):=a_{j2}(\rx,\tPsi_n, \nabla\tpsi_n)\quad\tx{for $j=1,2$}, \\
  &\mcl{L}_1^{(n)}(\psi, \Psi):=
  \psi_{x_1x_1}+2{a}_{12}^{(n)}\psi_{x_1x_2}
  -{a}_{22}^{(n)}\psi_{x_2x_2}+\ba_1\psi_{x_1}+\bb_1\Psi_{x_1}+\bb_2\Psi,\\
  &\hat{\mathfrak{f}}_1^{(n)}(\rx)
  :=\mathfrak{f}_1(\rx, \tPsi_n, \nabla\tPsi_n, \nabla\tpsi_n)
  -\mcl{L}_1^{(n)}(0, \Psi_{\rm bd}),\\
  & \hat{\mathfrak{f}}_2^{(n)}(\rx):=\mathfrak{f}_2(\rx, \tPsi_n, \nabla\tpsi_n)-\mcl{L}_2(0, \Psi_{\rm bd}),
  \end{split}
\end{equation}
where $a_{j2}(\rx, z, {\bf q})$, $\mathfrak{f}_1(\rx, z, {\bf p}, {\bf q})$ and $\mathfrak{f}_2(\rx, z, {\bf q})$ are given by \eqref{def-as}, \eqref{def-f1} and \eqref{def-f2}, respectively. It is easy to see that $\til{a}_{j2}^{(n)}\in C^{\infty}(\ol{\Om}_L)$, $\hat{\mfrak{f}}_k^{(n)}\in C^2(\ol{\Om_L})$, and $\mcl{L}_1^{(n)}$ is a linear differential operator with smooth coefficients.

Now we have the following lemma on the well-posedness of the problem with smooth coefficients.
\begin{lemma}
\label{lemma-wp-lbvp-approx}
For each $n\in \mathbb{N}$, the linear boundary value problem:
\begin{equation}
\label{lbvp-iter-approx}
\begin{split}
&\begin{cases}
\mcl{L}_1^{(n)}(v, w)=\hat{\mathfrak{f}}_1^{(n)}\\
\mcl{L}_2 (v, w)=\hat{\mathfrak{f}}_2^{(n)}
\end{cases}\quad\tx{in $\Om_L$},\\
&v=w_{x_1}=0,\quad v_{x_1}=g_1
\quad \tx{on $\Gamen$},\\
&w=0 \quad \tx{on $\Gam_{L}$},\\
&v_{x_2}=w_{x_2}=0 \quad \tx{on $\Lambda_{L}$}.
\end{split}
\end{equation}
has a unique solution $(v, w)\in [H^4(\Om_L)]^2$, and the solution satisfies
 \begin{equation}
    \label{estimate-sol-approx}
      \begin{split}
      &\|v\|_{H^4(\Om_L)}\le C\left(\|\hat{\mathfrak{f}}^{(n)}_1\|_{H^3(\Om_L)}+
\|\hat{\mathfrak{f}}^{(n)}_2\|_{H^2(\Om_L)}+\|g_1\|_{C^3(\ol{\Gamen})}\right),\\
      &\|w\|_{H^4(\Om_L)}\le
  C\left(\|\hat{\mathfrak{f}}^{(n)}_1\|_{H^2(\Om_L)}+\|\hat{\mathfrak{f}}_2^{(n)}\|_{H^2(\Om_L)}+\|g_1\|_{C^2(\ol{\Gamen})}\right)
      \end{split}
    \end{equation}
 for a constant $C>0$ depending only on $(\gam, J_0, S_0, \rho_0, E_0, \epsilon_0, L)$. Also, the solution $(v, w)$ satisfies
 \begin{equation}
 \label{comp-cond-sol-lbvp-approx}
   \der_{x_2}^kv=\der_{x_2}^k w=0\;
\tx{on $\Lambda_{L}$ in the sense of trace, for $k=1,3$ }.
 \end{equation}
 \end{lemma}
We first prove the well-posedness of \eqref{lbvp-iter2} and the estimate \eqref{estimate-sol-lbvp-H4} by assuming Lemma \ref{lemma-wp-lbvp-approx}. We  will give the proof of Lemma \ref{lemma-wp-lbvp-approx} later.

 {\bf Step 2. The well-posedness of \eqref{lbvp-iter2}.}
Assume that Lemma \ref{lemma-wp-lbvp-approx} holds.
For each $n\in \mathbb{N}$, let $(v_n, w_n)\in [H^4(\Om_L)]^2$ be the solution to the problem \eqref{lbvp-iter-approx}. It follows from \eqref{estimate-sol-approx}, \eqref{comp-cond-sol-lbvp-approx} and Morrey's inequality that the sequence $\{(v_n, w_n)\}_{n\in \mathbb{N}}$ is uniformly bounded in $C^{2,\frac 23}(\ol{\Om_L})$. By Arzel\`{a}-Ascoli theorem, the weak compactness property of $H^4(\Om_L)$, and by \eqref{comp-cond-sol-lbvp-approx}, there exists a subsequence $\{(v_{n_k}, w_{n_k})\}$ and the functions $(\psi, \hat{\Psi})$ so that
\begin{align*}
(v_{n_k}, w_{n_k})\to (\psi, \hat{\Psi})\quad \text{in}\,\,\,\, [C^{2,\frac 12}(\ol{\Om_L})]^{2}\quad
\tx{and}\quad
(v_{n_k}, w_{n_k})\rightharpoonup(\psi, \hat{\Psi})\quad \text{in}\,\,\,\, [H^4({\Om_L})]^2.
\end{align*}
Furthermore, $\displaystyle{\der_{x_2}^k\psi=\der_{x_2}^k \hat{\Psi}=0}$ holds  on $\Lambda_L$ in the sense of trace for $k=1$ and $3$.
These properties of $(\psi,\hat{\Psi})$ combined with Lemma \ref{lemma-coeff-reg} and \eqref{distance-approx-origina} imply that $(\psi,\hat{\Psi})$ is a classical solution to \eqref{lbvp-iter2}, and satisfies the estimate
\begin{equation}
\label{estimate-linear-H4-sol-new2018}
      \begin{split}
      &\|\psi\|_{H^4(\Om_L)}\le C\left(\|\hat{\mathfrak{f}}_1\|_{H^3(\Om_L)}+
\|\hat{\mathfrak{f}}_2\|_{H^2(\Om_L)}+\|g_1\|_{H^3({\Gamen})}\right),\\
      &\|\hat{\Psi}\|_{H^4(\Om_L)}\le
  C\left(\|\hat{\mathfrak{f}}_1\|_{H^2(\Om_L)}+
\|\hat{\mathfrak{f}}_2\|_{H^2(\Om_L)}+\|g_1\|_{H^2({\Gamen})}\right)
      \end{split}
    \end{equation}
where the estimate constant $C$ is from \eqref{estimate-sol-approx}. This finishes the proof of Proposition \ref{proposition-wp-lbvp}.
\end{proof}
\end{proposition}

To make the proof of Proposition \ref{proposition-wp-lbvp} complete, we now prove Lemma \ref{lemma-wp-lbvp-approx}.
 \begin{proof}[Proof of  Lemma \ref{lemma-wp-lbvp-approx}]
 We prove this lemma by the Galerkin method with the aid of the a priori $H^1$ estimate established in Propositions \ref{lemma-H1-estimate}. The proof is divided into 4 steps.

{\textbf{Step 1. Approximate Problems.}} Let us denote the open interval $\Gam:=(-1,1)$, and define the standard inner product  $\langle \cdot, \cdot\rangle$  in $L^2(\Gam)$ by
\begin{equation*}
\langle\zeta_1,\zeta_2\rangle=\int_{-1}^1 \zeta_1(x_2)\zeta_2(x_2)\;dx_2.
\end{equation*}
For each $k\in \mathbb{Z}^+$, we define a function $\eta_k$ by
\begin{equation}
\label{eigenfunction}
  \eta_k(x_2):=\cos (k\pi x_2).
\end{equation}
The set $\mathfrak{E}:=\{\eta_k\}_{k=0}^{\infty}$ is the collection of all eigenfunctions to the eigenvalue problem
\[
-\eta''=\om \eta \,\, \text{on}\,\, \Gam,\quad \eta'(\pm 1)=0,
 \]
 and the corresponding eigenvalues are $\{\om_k: \om_k=(k\pi)^2,\,\,k\in\mathbb{Z}_+\}$. It can be directly cheked that the set $\mathfrak{E}$  forms both an orthonormal basis in $L^2(\Gam)$, and an orthogonal basis in $H^1(\Gam)$.

We fix $n\in \mathbb{N}$, and consider the linear boundary value problem \eqref{lbvp-iter-approx} for each fixed $n$.
 For each $m=1,2,\cdots$, let $(V_m, W_m)$ be of the form
\begin{equation}
\label{def-projection-n}
V_m(x_1,x_2)=\sum_{j=0}^m \vartheta_j(x_1)\eta_j(x_2), \,\,  W_m(x_1,x_2)=\sum_{j=0}^m \Theta_j(x_1)\eta_j(x_2)\quad\tx{for ${\rx}=(x_1,x_2)\in \Om_L$}.
\end{equation}
For each $m\in \mathbb{N}$, we first find a solution $(V_m, W_m)$ satisfying
\begin{equation}\label{galerkin-he-system}
\langle \mcl{L}^{(n)}_1(V_m, W_m), \eta_k\rangle=\langle \hat{\mathfrak{f}}_1^{(n)}, \eta_k\rangle\,\,\tx{and}\,\,
\langle \mcl{L}_2(V_m, W_m), \eta_k\rangle = \langle \hat{\mathfrak{f}}_2^{(n)}, \eta_k\rangle
\quad\tx{for $0<x_1<L$}
\end{equation}
for all $k=0,1,\cdots, m$,
with boundary conditions
\begin{equation}
\label{galerkin-bc-others}
V_m=\der_{x_1}W_m=0 \,\,\tx{on}\,\, \Gamen,\quad \der_{x_1}V_m=\sum_{j=0}^m \langle g_1, \eta_j\rangle \eta_j\,\,\tx{on}\,\,\Gamen,
\quad W_m=0\,\,\tx{on}\,\, \Gam_L.
\end{equation}
It follows from $\eta_j'(\pm 1)=0$ that  $(V_m,
W_m)$ automatically satisfy the slip boundary conditions
\begin{equation}
\label{galerkin-slip-bc}
\der_{x_2}V_m=\der_{x_2}W_m=0\quad\tx{on}\,\,\Lambda_{L}.
\end{equation}
One can directly check from \eqref{distance-approx-origina} that
\begin{equation}
\label{galerkin-coeff}
  {a}_{12}^{(n)}=0\,\,\,\,\tx{on}\,\,\,\, \Lambda_L\quad \text{and}\quad
  \|\tpsi_n\|_{H^4(\Om_L)}+\|\tPsi_n\|_{H^4(\Om_L)}\le \delta_1(\epsilon_0).
\end{equation}
If $(V_m, W_m)$ are smooth, and satisfy \eqref{galerkin-he-system}--\eqref{galerkin-slip-bc}, then for any function $\mcl{U}(x_1)$, one has
\begin{equation}
\label{weak-formulation-galerkin}
\begin{split}
 \mcl{I}^{(n)}_L(V_m ,W_m, \mcl{U})
 &:=\int_{\Om_L} \mcl{W}(x_1)\der_{x_1}V_m \mcl{L}_1^{(n)}(v_m, w_m)-W_m\mcl{L}_2(V_m, W_m)\,d\rx\\
 &=\int_{0}^L\sum_{k=0}^m \mcl{U}(x_1)\langle \mcl{L}^{(n)}_1(V_m ,W_m), \eta_k\rangle  \vartheta_k'- \langle \mcl{L}_2(V_m ,W_m), \eta_k\rangle \Theta_k\,
 dx_1\\
 &=\int_0^L\sum_{k=0}^m \mcl{U}(x_1)\langle \hat{\mathfrak{f}}_1^{(n)}, \eta_k\rangle  \vartheta_k' - \langle \hat{\mathfrak{f}}_2, \eta_k\rangle \Theta_k\,dx_1\\
 &=\int_{\Om_L} \mcl{U}(x_1)\hat{\mathfrak{f}}_1^{(n)} \der_{x_1}V_m -W _m\hat{\mathfrak{f}}_2^{(n)}\,d\rx.
 \end{split}
\end{equation}
Substituting $\mcl{U}$ in \eqref{weak-formulation-galerkin} by the function $\mcl{W}(x_1)$ constructed in Step 3 of the proof of Proposition \ref{lemma-H1-estimate}
yields
\begin{equation}
\label{estimate-H1-galerkin}
  \|V_m\|_{H^1(\Om_L)}+\|W_m\|_{H^1(\Om_L)}
  \le C\left(
   \|\hat{\mathfrak{f}}_1^{(n)}\|_{L^2(\Om_L)}+\|\hat{\mathfrak{f}}_2^{(n)}\|_{L^2(\Om_L)}+\|g_1\|_{C^2(\ol{\Gamen})}
  \right)
\end{equation}
for a constant $C>0$ depending only on $(\gam, J_0, S_0, \rho_0, E_0, \epsilon_0, L)$.

{\textbf{Step 2. Analysis on ODE problems.}}
For each $k\in\mathbb{Z}^+$, set
\begin{equation*}
  f_{1,k}:= \langle \hat{\mathfrak{f}}_1^{(n)}, \eta_k\rangle,\quad
  f_{2,k}:= \langle \hat{\mathfrak{f}}_2^{(n)}, \eta_k\rangle,\quad
  g_{1,k}:= \langle g_1, \eta_k\rangle.
\end{equation*}
By the orthonormality of the set $\mathfrak{E}$ in $L^2(\Gam)$, for each $m\in \mathbb{N}$, \eqref{galerkin-he-system}--\eqref{galerkin-slip-bc} for $(V_m, W_m)$ is written as a boundary value problem for $\{(\vartheta_k,\Theta_k)\}_{k=0}^m$ as follows:
\begin{equation}
\label{galerking-ibvp-eqns}
\begin{cases}
\vartheta_k''+\sum_{j=0}^m (2 a_{12}^{(jk)}+a_1^{(jk)})\vartheta_j'
+j^2a_{22}^{(jk)}\vartheta_j+b_1^{(jk)}\Theta_j'+b_2^{(jk)}\Theta_j
=f_{1,k}\\
\Theta_k''-(k\pi)^2\Theta_k-\sum_{j=0}^m h_1^{(jk)}\Theta_j+h_2^{(jk)}\vartheta_j'=f_{2,k}
\end{cases}
\,\,\tx{for $0<x_1<L$},
\end{equation}
\begin{equation}
\label{galerkin-bc-n}
\vartheta_k(0)=\Theta_k'(0)=0,\quad \vartheta'_k(0)=g_{1,k}\quad \tx{and}\quad \Theta_k(L)=0,
\end{equation}
where $k=0, \cdots, m$, and $(a_{12}^{(jk)}, a_{22}^{(jk)}, a_1^{(jk)}, b_1^{(jk)}, b_2^{(jk)}, h_1^{(jk)}, h_2^{(jk)})$ are defined by
\begin{equation*}
a_{12}^{(jk)}=\langle a_{12}^{(n)}\eta_j', \eta_k\rangle,\quad (a_{22}^{(jk)}, a_1^{(jk)}, b_1^{(jk)}, b_2^{(jk)}, h_1^{(jk)}, h_2^{(jk)})=\langle (a_{22}^{(n)}, \ba_1, \bb_1, \bb_2, \bar{h}_1, \bar{h}_2)\eta_j, \eta_k \rangle.
\end{equation*}

{\emph{Claim 2. For each $m\in \mathbb{N}$, the boundary value problem \eqref{galerking-ibvp-eqns}--\eqref{galerkin-bc-n} for $\{(\vartheta_k,\Theta_k)\}_{k=0}^m$ has a unique smooth solution on the interval $[0,L]$.}}

Let us denote
\begin{equation*}
{\bf X}:=(X_1,\cdots,X_{4(m+1)})^T:=(\vartheta_0,\cdots, \vartheta_{m}, \vartheta_0',\cdots, \vartheta_m', \Theta_0',\cdots, \Theta_m', \Theta_0,\cdots, \Theta_{m})^T.
\end{equation*}
Now the system \eqref{galerking-ibvp-eqns} can be written  as a first order ODE system for ${\bf X}$ in the following form:
\begin{equation}
\label{ODE-n}
{\bf X}'={\mathbb A}{\bf X} +\bfF, \quad x_1\in [0, L],
\end{equation}
where  ${\mathbb A}: [0, L]\to \R^{4(m+1)\times 4(m+1)}$ and $\bfF:[0, L]\to \R^{4(m+1)}$ are smooth functions represented in terms of $(a_{l2}^{(jk)}, b_l^{(jk)}, h_l^{(jk)}, f_l^{(jk)})$ for $l=1,2$ and $j,k=0,1,\cdots, m$.

Define projection mapping $\Pi_i$ ($i=1$, $2$) by
\begin{equation*}
\Pi_1{\bf X}=(X_1,\cdots, X_{3(m+1)}, 0,\cdots,0)^T\quad \text{and}\quad
\Pi_2:={\bf{Id}}-\Pi_1
\end{equation*}
where ${\bf Id}$ is an identity map from $\R^{4(m+1)}$ to $\R^{4(m+1)}$. Then, we have
\begin{equation}
\label{ODE-bcs-n}
  (\Pi_1{\bf X})(0)=\sum_{j=0}^m g_{1j}\hat{\bf e}_{j+m+1}=:{\bf P}_0,\quad
  (\Pi_2{\bf X})(L)={\bf 0},
\end{equation}
where we set $\hat{\bf e}_l:=(\delta_{il})_{i=1}^{4(m+1)}$, i.e., $\hat{\bf e}_l$ is the unit vector in the direction of $X_l$ for $l=1,\cdots, 4(m+1)$.

To verify {\emph{Claim 2}},  it suffice to show that there exists a $C^1$ function $\bfX: [0, L]\to \R^{4(m+1)}$ satisfying
\begin{equation}
\label{integral-eqn-n}
\begin{split}
&\Pi_1\bfX(x_1)-\int_0^{x_1}\Pi_1{\mathbb A}\bfX(t)\;dt={\bf P}_0+\int_0^{x_1} \Pi_1\bfF(t)\;dt,\\
&\Pi_2\bfX(x_1)-\int_L^{x_1}\Pi_2{\mathbb A}\bfX(t)\;dt=\int_L^{x_1} \Pi_2\bfF(t)\;dt.
\end{split}
\end{equation}
Set $\mathfrak{B}:=C^1([0,L];\R^{4(m+1)})$. Define a linear operator
$\mathfrak{K}:\mathfrak{B}\to \mathfrak{B}$ by
\begin{equation*}
\mathfrak{K} \bfX (x_1)= \Pi_1\int_0^{x_1} {\mathbb A}\bfX(t)\;dt+\Pi_2\int_L^{x_1} {\mathbb A}\bfX(t)\;dt.
\end{equation*}
Let us rewrite \eqref{integral-eqn-n} as
\begin{equation}
\label{integral-formulation}
({\bf Id}-\mathfrak{K})\bfX(x_1)={\bf P}_0+\Pi_1\int_0^{x_1} \bfF(t)\;dt+\Pi_2\int_L^{x_1} \bfF(t)\;dt,
\end{equation}
where ${\bf P}_0$ is given in \eqref{ODE-bcs-n}.
Since ${\mathbb A}: [0, L]\rightarrow \R^{4(m+1)\times 4(m+1)}$ is smooth, there exists a constant $k_0>0$  satisfying
\begin{equation*}
\|\mathfrak{K}\bfX\|_{C^2([0, L])}\le k_0\|\bfX\|_{C^1([0, L])}\quad\tx{for all}\,\,\bfX\in \mathfrak{B},
\end{equation*}
and this implies that the linear mapping $\mathfrak{K}:\mathfrak{B}\rightarrow \mathfrak{B}$ is compact by Arzel\`{a}-Ascoli theorem.

Suppose that $({\bf I}-\mathfrak{K})\bfX_*=0$. Then $\bfX_*$ solves \eqref{ODE-n}--\eqref{ODE-bcs-n} with ${\bf F}={\bf 0}$ and ${\bf P}_0=0$, and the corresponding $(V_m, W_m)$ given by \eqref{def-projection-n} satisfy \eqref{galerkin-he-system}--\eqref{galerkin-slip-bc} with $\hat{\mathfrak{f}}_1^{(n)}=\hat{\mathfrak{f}}_2^{(n)}=g_1=0$. Then the estimate \eqref{estimate-H1-galerkin} implies that $V_m=W_m=0$ in $\Om_L$, from which we obtain that $\bfX_*={\bf 0}$ on $[0,L]$. Then, we apply Fredholm alternative theorem to conclude that \eqref{integral-formulation} has a unique solution $\bfX\in \mathfrak{B}$.

Furthermore, the bootstrap argument with using \eqref{integral-eqn-n} implies that the solution ${\bf X}$ is smooth on $[0,L]$. This verifies {\emph{Claim 2}}.

{\textbf{Step 3. Estimate for $(V_m, W_m)$.}}
In order to prove the estimate \eqref{estimate-linear-H4-sol-new2018} in Lemma \ref{lemma-wp-lbvp-approx}, we will estimate $(\|V_m\|_{H^4(\Om_L)}, \|W_m\|_{H^4(\Om_L)})$, then pass to the limit $m\to \infty$. The estimate of $(\|V_m\|_{H^4(\Om_L)}, \|W_m\|_{H^4(\Om_L)})$ is given in the following lemma.
\begin{lemma}\label{lemma-H4-galerkin}
Fix $n\in \mathbb{N}$. For each $m\in \mathbb{N}$, let $(V_m, W_m)$ given in the form \eqref{def-projection-n} be the solution to \eqref{galerkin-he-system}--\eqref{galerkin-bc-others}. Then, there exists a constant $C>0$ depending only on $(\gam, J_0, S_0, \rho_0, E_0, \epsilon_0, L)$ such that
\begin{equation}\label{estimate-H4-galerkin-new}
\begin{split}
&\|V_m\|_{H^4(\Om_L)}\le C\Bigl(\|\hat{\mathfrak{f}}_1^{(n)}\|_{H^3(\Om_L)}
+\|\hat{\mathfrak{f}}_2^{(n)}\|_{H^2(\Om_L)}
+\|g_1\|_{C^3(\ol{\Gamen})}\Bigr),\\
&\|W_m\|_{H^4(\Om_L)}\le
C\Bigl(\|\hat{\mathfrak{f}}_1^{(n)}\|_{H^2(\Om_L)}+\| \hat{\mathfrak{f}}_2^{(n)}\|_{H^2(\Om_L)}
+\|g_1\|_{C^2(\ol{\Gamen})}\Bigr).
\end{split}
\end{equation}
And, the estimate constant $C$ in the above is independent of $n,m\in \mathbb{N}$ and $(\tpsi, \tPsi)\in \mcl{J}_{\delta, L}$.
\end{lemma}
This lemma can be proved by applying standard theory for elliptic and hyperbolic equations, using \eqref{estimate-H1-galerkin}, and by the bootstrap argument with additional careful computations to treat corner points on the boundary of $\der \Om_L$.
A detailed proof of Lemma \ref{lemma-H4-galerkin} is given in Appendix \ref{appendix-1}.

{\textbf{Step 4. Taking the limit.}} Now, we are ready to complete the proof of Lemma \ref{lemma-wp-lbvp-approx}. For a fixed $n\in \mathbb{N}$, let $(V_m, W_m)$ given in the form \eqref{def-projection-n} be the solution to \eqref{galerkin-he-system}--\eqref{galerkin-bc-others}. By Lemma \ref{lemma-H4-galerkin}, the sequence $\{(V_m, W_m)\}_{m=1}^{\infty}$ is bounded in $[H^4(\Om_L)]^2$. Due to $\eta_k'(\pm 1)=\eta_k'''(\pm 1)=0$ for $k=0,1,2\cdots$,  each $(V_m, W_m)$ satisfies the compatibility condition $\der_{x_2}^jV_m=\der_{x_2}^jW_m=0$ on $\Lambda_L$ for $j=1$ and $3$. Then it follows from Morrey's inequality that the sequence $\{(V_m, W_m)\}_{m\in \mathbb{N}}$ is bounded in $[C^{2,\frac 23}(\ol{\Om_L})]^2$. Then, by adjusting the limiting argument given in Step 2 in this proof, we can extract a subsequence $\{(V_{m_k}, W_{m_k})\}$ so that the limit $(v, w)\in [H^4(\Om_L)\cap C^{2,\frac 12}(\ol{\Om_L})]^2$ of the subsequence is a classical solution to \eqref{lbvp-iter-approx}, and satisfies the compatibility boundary conditions \eqref{comp-cond-sol-lbvp-approx}. Furthermore, the estimate \eqref{estimate-sol-approx} given in Lemma \ref{lemma-wp-lbvp-approx} is directly obtained from \eqref{estimate-H4-galerkin-new} combined with the weak $H^4$ convergence of  $\{(V_{m_k}, W_{m_k})\}$. This completes the proof of Lemma \ref{lemma-wp-lbvp-approx}.
\end{proof}

From Proposition \ref{proposition-wp-lbvp}, the well-posedness of \eqref{lbvp-iter} directly follows.
\begin{corollary}
\label{corollary-wp-original-lbvp}
Fix $\epsilon_0\in(0, T_{\max})$, and let $\bar{L}\in(0, T_{\max}-\epsilon_0]$ be from Proposition \ref{lemma-H1-estimate}. For $L\in(0, \bar{L}]$ and $\delta\in(0, \frac{\delta_1(\epsilon_0)}{2}]$, let the iteration set $\mcl{J}_{\delta, L}$ be given by \eqref{def-iter}. Then, for each $(\tpsi, \tPsi)\in \mcl{J}_{\delta, L}$, the associated linear boundary value problem \eqref{lbvp-iter} has a unique solution $(\psi, {\Psi})\in [H^4(\Om_L)]^2$ that satisfies
    \begin{equation}
    \label{estimate-sol-lbvp-H4-original}
      \begin{split}
      &\|{\Psi}\|_{H^4(\Om_L)}\le
  C\left(\|{\mathfrak{f}}_1\|_{H^2(\Om_L)}+\|{\mathfrak{f}}_2\|_{H^2(\Om_L)}
+\|g_1\|_{C^2(\ol{\Gamen})}+\|g_2\|_{C^4(\ol{\Gamen})}+\|\Psi_{ex}\|_{C^4(\ol{\Gamen})} \right),\\
&\|\psi\|_{H^4(\Om_L)}\le C\left(\|{\mathfrak{f}}_1\|_{H^3(\Om_L)}+
\|{\mathfrak{f}}_2\|_{H^2(\Om_L)}
+\|g_1\|_{C^3(\ol{\Gamen})}+\|g_2\|_{C^4(\ol{\Gamen})}+\|\Psi_{ex}\|_{C^4(\ol{\Gamen})} \right)
      \end{split}
    \end{equation}
 for a constant $C>0$ depending only on $(\gam, J_0, S_0, \rho_0, E_0, \epsilon_0, L)$. Furthermore, the solution $(\psi, {\Psi})$ satisfies the compatibility conditions
 \begin{equation}
 \label{comp-cond-sol-lbvp-original}
   \der_{x_2}^k\psi=\der_{x_2}^k\Psi=0\;
\tx{for $k=1,3$ on $\Lambda_{L}$ in the sense of trace}.
 \end{equation}
\end{corollary}
\begin{proof}
By Proposition \ref{proposition-wp-lbvp}, for each $(\tpsi, \tPsi)\in \mcl{J}_{\delta, L}$, Since the associated linear boundary value problem \eqref{lbvp-iter2} has a unique solution $(\psi, \hat{\Psi})\in [H^4(\Om_L)]^2$, and it satisfies the estimate \eqref{estimate-sol-lbvp-H4}. Then, $(\psi, \Psi)$ given by
$$
(\psi, \Psi)=(\psi ,\hat{\Psi})+(0, \Psi_{\rm bd})
$$
solves \eqref{lbvp-iter}. And, the estimate \eqref{estimate-sol-lbvp-H4-original} can be easily checked by using \eqref{definition-bds}, \eqref{definition-f-hats} and \eqref{estimate-sol-lbvp-H4}. Finally, the compatibility condition \eqref{comp-cond-sol-lbvp-original} holds for $(\psi, \Psi)$ due to \eqref{comp-cond1} and \eqref{comp-cond-sol-lbvp}.
\end{proof}


\section{Nonlinear stability of supersonic solutions}\label{secNL}
In this section, we apply the results obtained in Section \ref{section-thm-pf-potential}
to achieve nonlinear structural stability of supersonic solutions for the Euler-Poisson system. We first prove Theorem \ref{theorem-2-potential}  in Section \ref{subsection-thm-pf1}, which gives the structural stability of irrotational supersonic solutions. The iteration scheme for the supersonic solutions with non-zero vorticity is presented in detail in Section \ref{seciteration}. The structural stability of supersonic solutions with non-zero vorticity (Theorem \ref{theorem-1}) is proved in Section \ref{subsec-proof-theorem1-a}.

\subsection{Irrotational flows (Proof of Theorem \ref{theorem-2-potential})}
\label{subsection-thm-pf1}

We prove Theorem \ref{theorem-2-potential} in two steps.

{\bf Step 1. } In this step, we prove the existence of a solution to Problem \ref{problem-2-potential}.

Fix $\epsilon_0\in(0, T_{\max})$, and let $\bar{L}\in (0, T_{\max}-\epsilon_0]$ and $\delta_1(\epsilon_0)$ be from Proposition \ref{lemma-H1-estimate}. We fix $L\in(0, \bar{L}]$. For a constant $\delta\in(0, \frac{\delta_1(\epsilon_0)}{2}]$ to be specified later, let the iteration set $\mcl{J}_{\delta, L}$ be given by \eqref{def-iter}. We define a mapping $\mathfrak{F}: \mcl{J}_{\delta, L}\rightarrow [H^4(\Om_L)]^2$ by
\begin{equation*}
  \mathfrak{F}(\til{\psi}, \til{\Psi}):=(\psi, \Psi),
\end{equation*}
where $(\psi, \Psi)$ is the solution to the linear boundary value problem \eqref{lbvp-iter} associated with $(\til{\psi}, \til{\Psi})\in \mcl{J}_{\delta, L}$.
By Corollary \ref{corollary-wp-original-lbvp}, the mapping $\mathfrak{F}$ is well defined. Furthermore, if
$\mathfrak{F}(\psi_*, \Psi_*)=(\psi_*, \Psi_*)$ holds, then $(\vphi, \Phi):=(\vphi_0, \Phi_0)+(\psi_*, \Psi_*)$ solves Problem \ref{problem-2-potential} where $(\vphi_0, \Phi_0)$ is given by \eqref{definition-background-potential}.

Let us define
\begin{equation*}
\sigma_p=\sigma(b, u_{\rm en}, 0, E_{\rm en}, \Phi_{\rm ex}, S_0),
\end{equation*}
where $\sigma(b, u_{\rm en}, 0, E_{\rm en}, \Phi_{\rm ex}, S_0)$ is defined by \eqref{definition-sigma}.

By using \eqref{def-rho}, \eqref{coefficients-vphi-eqn}, \eqref{def-as}, \eqref{def-f1}, \eqref{def-coefficients3},\eqref{def-f2}, and Lemma \ref{lemma-coeff-reg}, it can be directly checked that there exists a constant $C>0$ depending only on $(\gam, J_0, S_0, \rho_0, E_0, \epsilon_0, L)$ to satisfy
\begin{align}
\label{estimate-f1}
&\|\mathfrak{f}_1(\cdot, \tPsi, \nabla\tPsi, \nabla\tpsi)\|_{H^3(\Om_L)}\le
C(\|\til\psi\|_{H^4(\Om_L)}^2
+\|\til\Psi\|_{H^4(\Om_L)}^2),\\
\label{estimate-f2}
&\|\mathfrak{f}_2(\cdot, \nabla\tPsi, \nabla\tpsi)\|_{H^2(\Om_L)}\le C\left(\|\til\psi\|_{H^4(\Om_L)}^2
+\|\til\Psi\|_{H^4(\Om_L)}^2 +\sigma_p\right)
\end{align}
for all $(\tPsi, \tpsi)\in \mcl{J}_{\delta, L}$.
We combine \eqref{estimate-f1}--\eqref{estimate-f2} with Corollary \ref{corollary-wp-original-lbvp} to get \begin{equation}\label{estimate-image-F}
\begin{split}
 \|\mathfrak{F}(\tpsi, \tPsi)\|_{H^4(\Om_L)}
 \le C_*\Bigl(\delta^2+\sigma_{p}\Bigr).
  \end{split}
\end{equation}
for a constant $C_*>0$ depending only on $(\gam, J_0, S_0, \rho_0, E_0, \epsilon_0, L)$. Therefore, if $(\delta,\sigma_p)$ satisfy
\begin{equation}
\label{parameter-condition1}
  \delta\le \frac 12 \min\{\frac{1}{C_*}, \delta_1(\epsilon_0)\}\quad\tx{and}\quad
  \sigma_{p}\le \frac{\delta}{4C_*},
\end{equation}
then $\mathfrak{F}$ maps $\mcl{J}_{\delta,L}$ into itself.

By the Rellich's theorem, the iteration set $\mcl{J}_{\delta,L}$ is compact in $[H^3(\Om_L)]^2$. So if $\mathfrak{F}$ is continuous in $[H^3(\Om_L)]^2$, then it follows from the Schauder fixed point theorem that $\mathfrak{F}$ has a fixed point in $\mcl{J}_{\delta, L}$.
Suppose that a sequence $\{(\tpsi_j, \tPsi_j)\}_{j\in \mathbb{N}}$ in $\mcl{J}_{\delta, L}$ converges to $(\tpsi_{\infty}, \tPsi_{\infty})\in \mcl{J}_{\delta, L}$ in $[H^3(\Om_L)]^2$. For each $j\in \mathbb{N}\cup\{\infty\}$, we set
\begin{equation*}
  (\psi_j, \Psi_j):=\mathfrak{F}(\tpsi_j, \tPsi_j).
\end{equation*}
Since $\{(\psi_j, \Psi_j)\}_{j\in \mathbb{N}}$ is bounded in $[H^4(\Om_L)\cap C^{2,\frac 23}(\ol{\Om_L})]^2$, by the Rellich's theorem and the Arzel\'{a}-Ascoli theorem, one can take a subsequence $\{(\psi_{j_k}, \Psi_{j_k})\}$ which converges to some $(\psi_{\infty}^*, \Psi_{\infty}^*)$ in $[H^3(\Om_L)\cap C^{2,\frac 12}(\ol{\Om_L})]^2$.
Due to Lemma \ref{lemma-coeff-reg} and the $H^3$-convergence of $\{(\tpsi_j, \tPsi_j)\}$ to $(\tpsi_{\infty}, \tPsi_{\infty})$, $(\psi_{\infty}^*, \Psi_{\infty}^*)$  is a classical solution to the linear boundary value problem \eqref{lbvp-iter} associated with $(\tpsi_{\infty}, \tPsi_{\infty})$. Then the uniqueness of the solution to \eqref{lbvp-iter} given by Corollary \ref{corollary-wp-original-lbvp} implies that $(\psi^*_{\infty}, \Psi^*_{\infty})=(\psi_{\infty}, \Psi_{\infty})$. And, this shows that $\mathfrak{F}: \mcl{J}_{\delta,L}\rightarrow \mcl{J}_{\delta, L}$ is continuous in $[H^3(\Om_L)]^2$. Therefore, $\mathfrak{F}$ has a fixed point in $\mcl{J}_{\delta,L}$ provided that \eqref{parameter-condition1} holds.

Let $(\psi, \Psi)\in \mcl{J}_{\delta,L}$ be a fixed point of $\mathfrak{F}$. Then, the functions $(\vphi, \Phi):=(\vphi_0, \Phi_0)+(\psi, \Psi)$ solve Problem \ref{problem-2-potential} and satisfy \eqref{supersonicity-potential} and \eqref{comp-cond-nlbvp-potential} with  the constant $\kappa_0>0$ from Lemma \ref{lemma-coeff-reg}(a). Furthermore, it directly follows from \eqref{estimate-sol-lbvp-H4-original} and \eqref{estimate-f1}--\eqref{parameter-condition1} that $(\vphi, \Phi)$ satisfies
\begin{equation}
\label{pre-esti-nonl-irrotational}
\|(\vphi-\vphi_0, \Phi-\Phi_0)\|_{H^4(\Om_L)}\le C_{\sharp}(\delta \|(\vphi-\vphi_0, \Phi-\Phi_0)\|_{H^4(\Om_L)}+\sigma_p)
\end{equation}
for a constant $C_{\sharp}>0$ depending only on $(\gam, J_0, S_0, \rho_0, E_0, \epsilon_0, L)$. Therefore, if $(\delta,\sigma_p)$ satisfy
\begin{equation}
\label{final-cond-exist}
\delta\le \frac 12 \min\{\frac{1}{C_*}, \frac{1}{C_{\sharp}},  \delta_1(\epsilon_0)\}\quad\tx{and}\quad
  \sigma_{p}\le \frac{\delta}{4C_*},
\end{equation}
then
the estimate \eqref{apriori-nlbvp-potential} can be obtained from \eqref{pre-esti-nonl-irrotational}.

{\bf Step 2.} To complete the proof of Theorem \ref{theorem-2-potential}, it remains to prove the uniqueness of a solution to Problem \ref{problem-2-potential}.
Let $(\vphi^{(k)}, \Phi^{(k)})$ ($k=1$, $2$) be two solutions to Problem \ref{problem-2-potential} that satisfy \eqref{apriori-nlbvp-potential}, \eqref{supersonicity-potential} and \eqref{comp-cond-nlbvp-potential}. For each $k=1,2$, set
\begin{equation*}
\begin{split}
  &(\psi^{(k)}, \Psi^{(k)}):=(\vphi^{(k)}, \Phi^{(k)})-(\vphi_0, \Phi_0),\quad \mathfrak{f}_1^{(k)}:=\mathfrak{f}_1(\rx, \Psi_k, \nabla\Psi_k, \nabla\psi_k),\\
&(a^{(k)}_{12}, a^{(k)}_{22},\mathfrak{f}_2^{(k)}):=(a_{12}, a_{22}, \mathfrak{f}_2)(\rx, \Psi_k, \nabla\psi_k),
\end{split}
\end{equation*}
where $(a_{12}, a_{22})(\rx,z,{\bf q})$, $\mathfrak{f}_1(\rx, z, {\bf p}, {\bf q})$ and $\mathfrak{f}_2(\rx, z, {\bf q})$ are defined by \eqref{def-f1}, \eqref{def-as} and \eqref{def-f2}, respectively.
Then, $(\check{\psi},\check{\Psi}):=(\psi_1, \Psi_1)-(\psi_2, \Psi_2)$ satisfies
\begin{equation*}
\begin{cases}
\mcl{L}_1^{(\psi_1, \Psi_1)}(\check{\psi},\check{\Psi})=
F_1\\
\mcl{L}_2(\check{\psi},\check{\Psi})=F_2
\end{cases}\;\;\tx{in}\;\;\Om_L,
\end{equation*}
where
\begin{equation*}
F_1= (\mathfrak{f}_1^{(1)}-\mathfrak{f}_1^{(2)})+
(a_{22}^{(1)}-a_{22}^{(2)})\der_{x_2x_2}\psi^{(2)}-2(a_{12}^{(1)}-a_{12}^{(2)})\der_{x_1x_2}\psi^{(2)},\,\, \quad F_2=\mathfrak{f}_2^{(1)}-\mathfrak{f}_2^{(2)},
\end{equation*}
and the linear differential operators $\mcl{L}_1^{(\psi_1, \Psi_1)}$ and $\mcl{L}_2$ are given by \eqref{def-L1} and \eqref{pfeqn-poisson-modi}, respectively.
Furthermore, $(\check{\psi},\check{\Psi})$ satisfies the homogeneous boundary conditions:
\begin{equation*}
\begin{aligned}
&(\check{\psi}, \check{\psi}_{x_1}, \check{\Psi}_{x_1})=(0,0,0)\,\,\tx{on}\;\;\Gamen,\quad &\check{\Psi}=0\,\,\tx{on}\;\;\Gamex,\quad  (\check{\psi}_{x_2}, \check{\Psi}_{x_2})=(0,0)\,\,\tx{on}\,\, \Lambda_L.
\end{aligned}
\end{equation*}

With the aid of \eqref{apriori-nlbvp-potential}, \eqref{def-f1}, \eqref{def-as}, \eqref{def-f2}, and Proposition \ref{lemma-H1-estimate}, the direct computations  yield
\begin{equation}
\label{est-H1-contraction-potential}
\|(\check{\psi},\check{\Psi})\|_{H^1(\Om_L)}\le C_{\flat}\sigma_p
\|(\check{\psi}, \check{\Psi})\|_{H^1(\Om_L)}
\end{equation}
for a constant $C_{\flat}>0$ depending only on $(\gam, J_0, S_0, \rho_0, E_0, \epsilon_0, L)$.
Finally, we choose $\delta$ as
\begin{equation*}
\delta= \frac 12 \min\{\frac{1}{C_*}, \frac{1}{C_{\sharp}},  \delta_1(\epsilon_0)\}.
 \end{equation*}
And, we choose $\bar{\sigma}_{\rm bd}$ as
\begin{equation*}
\bar{\sigma}_{\rm bd}=\min\left\{\frac{\delta}{4C_*}, \frac{1}{2C_{\flat}}\right\}.
\end{equation*}
Therefore,  if $\sigma_p\le \bar{\sigma}_{\rm bd}$, then the inequality \eqref{est-H1-contraction-potential} implies $\check{\psi}=\check{\Psi}=0$ so that  Problem \ref{problem-2-potential} has a unique solution that satisfies \eqref{apriori-nlbvp-potential}, \eqref{supersonicity-potential} and \eqref{comp-cond-nlbvp-potential}.
This finishes the proof of Theorem \ref{theorem-2-potential}. \hfill $\square$

\subsection{Flows with non-zero vorticity (Proof of Theorem \ref{theorem-1})} Theorem \ref{theorem-1} is also proved by the method of iteration. But the iteration scheme is more complicated than the case of irrotational flow due to the additional equations \eqref{eqn3-hdecomp}--\eqref{eqn4-hdecomp} coupled with \eqref{hdecomp-system}. In Section \ref{seciteration}, we explain the iteration scheme applied to prove Theorem \ref{theorem-1}, then we complete the proof of Theorem \ref{theorem-1} in Section \ref{subsec-proof-theorem1-a}.

\subsubsection{Iteration scheme for the flows with non-zero vorticity}\label{seciteration}
Suppose that $(\vphi, \phi, \Phi, S)$ is a solution to \eqref{hdecomp-system}--\eqref{eqn4-hdecomp} and \eqref{bc-vphi}--\eqref{bc-T}. In addition, suppose that the solution satisfies
\begin{equation}\label{nondegeneracy-vel-x1}
  \mcl{H}\left(S, \Phi-\frac 12|\nabla\vphi+\nabla^{\perp}\phi|^2\right) (\nabla\vphi+\nabla^{\perp}\phi)\cdot{\bf e}_1 \ge \eps\quad\tx{in $\ol{\Om_L}$}
\end{equation}
for some constant $\eps>0$.
For $(\vphi_0, \Phi_0)$ given by \eqref{definition-background-potential}, we set
\begin{equation*}
  (\psi, \Psi, Y):=(\vphi, \Phi, S)-(\vphi_0, \Phi_0, S_0).
\end{equation*}
Assume that $(\psi, \Psi, \phi)$ satisfies
\begin{equation}
\label{supersonicity1-full}
  c^2(\Phi_0+\Psi, \nabla\vphi+\nabla^{\perp}\phi)-|\nabla\vphi+\nabla^{\perp}\phi|^2<0
  \quad\tx{in $\ol{\Om_L}$}
\end{equation}
with $c(z, {\bf q})$ defined by \eqref{defeqsound}. To present a boundary value problem for the iteration scheme, we first introduce the following notations.
\begin{definition}
\label{definition-nlbvp-pert-full}
(i) For a given ${\bf r}=(r_1,r_2)\in \R^2$, we set ${\bf r}^{\perp}:=(r_2, -r_1)$. For $A_{j2}(z, {\bf q})(j=1,2)$ given by \eqref{coefficients-vphi-eqn}, we define
\begin{equation*}
  a_{j2}(\rx, z, {\bf q}, {\bf r})=A_{j2}(\Phi_0(\rx)+z, \nabla\vphi_0(\rx)+{\bf q}+{\bf r}^{\perp}).
\end{equation*}
Then $\tilde{\mcl{N}}_1(\psi, \Psi, \phi)$ is defined by
\begin{equation*}
  \begin{split}
\tilde{\mcl{N}}_1(\psi, \Psi, \phi):=
&\psi_{x_1x_1}+
2a_{12}(\rx, \Psi, \nabla\psi, \nabla\phi)\psi_{x_1x_2}
-a_{22}(\rx, \Psi, \nabla\psi, \nabla\phi)\psi_{x_2x_2}\\
&+\bara_1({\rx})\psi_{x_1}
+\bb_1({\rx})\Psi_{x_1}+\bb_2({\rx})\Psi
  \end{split}
\end{equation*}
where $(\ba_1, \bb_1, \bb_2)$ are given by \eqref{def-as}. Note that the definition of the nonlinear differential operator $\tilde{\mcl{N}}_1$ is almost the same as \eqref{defN1} except that the coefficients $a_{j2}$ with $j=1,2$ additionally depend on $\nabla\phi$.

(ii) For $(\bar h_1, \bar h_2)$ given by \eqref{def-coefficients3}, we define
\begin{equation*}
{\mcl{L}}_2(\psi, \Psi):=\Delta \Psi-\bar h_1(x_1)\Psi-\bar h_2(x_1)\psi_{x_1}
\end{equation*}
Note that the definition of $\mcl{L}_2$ is the same as \eqref{pfeqn-poisson-modi}.

(iii) For $\mcl{H}$ defined by \eqref{rho-hdecomp}, we set
\begin{equation*}
\begin{split}
  &{\bf M}(\rx, Y, \Psi, \nabla\psi, \nabla\phi)\\
  &:=
  \mcl{H}\left(S_0+Y,\Phi_0(\rx)+\Psi-\frac 12|\nabla\vphi_0(\rx)+\nabla\psi+\nabla^{\perp}\phi|^2\right)
  \left(\nabla\vphi_0(\rx)+\nabla\psi+\nabla^{\perp}\phi\right).
  \end{split}
\end{equation*}

(iv) For $\mathbb{M}=(m_{ij})_{i,j=1}^2\in \R^{2\times 2}$ and ${\bf w}=(w_1,w_1)\in \R^2$, we set $\mathbb{M}[{\bf w}, {\bf w}]:=\sum_{i,j=1}^2 m_{ij}w_iw_j$.

For $\rx\in \Om_L$, $\xi, z\in\R$, ${\bf p}, {\bf q}, {\bf r}\in \R^2$, and $\mathbb{M}\in \R^{2\times 2}$, we define
\begin{align*}
&F_1(\rx, \xi,{\bm \zeta}, z, {\bf p}, {\bf q}, {\bf r}, \mathbb{M}):=
  (\gam-1)\bar{u}'(x_1)\left(\bar u(x_1)\phi_{x_2}+\frac 12|{\bf q}+{\bf r}^{\perp}|^2\right)-{\bf p}\cdot{\bf q}-(\nabla\Phi_0(\rx)+{\bf p})\cdot {\bf r}^{\perp}\\
  &\phantom{F_1(\rx, \xi,{\bm \zeta}, z, {\bf p}, {\bf q}, {\bf r}, \mathbb{M}):=}
  +\frac{\Phi_0(\rx)+z-\frac 12 |\nabla\vphi_0(\rx)+{\bf q}+{\bf r}^{\perp}|^2}{S_0+Y}\left(\nabla\vphi_0(\rx)+{\bf q}+{\bf r}^{\perp}\right)\cdot {\bm\zeta}\\
  &\phantom{F_1(\rx, \xi,{\bm \zeta}, z, {\bf p}, {\bf q}, {\bf r}, \mathbb{M}):=}
  +\mathbb{M}[\nabla\vphi_0(\rx)+{\bf q}+{\bf r}^{\perp},\nabla\vphi_0(\rx)+{\bf q}+{\bf r}^{\perp}],\\
&F_2(\rx, z, {\bf p}, {\bf q}):=\bar{u}(x_1)p_1+\bar{E}(\rx)q_1
  +(\gam-1)\bar{u}'(x_1)(z-\bar{u}(\rx)q_1),\\
&\beta(\rx, z, {\bf q}, {\bf r}):=
(\gam-1)\left(\Phi_0(\rx)+z-\frac 12 |\nabla\vphi_0(\rx)+{\bf q}+{\bf r}^{\perp}|^2\right)-(\bar{u}(x_1)+q_1+r_2)^2.
\end{align*}
Finally, we define $f_1$ by
\begin{equation*}
  f_1(\rx, \xi,{\bm \zeta}, z, {\bf p}, {\bf q}, {\bf r}, \mathbb{M})
  :=\frac{F_1(\rx, \xi,{\bm \zeta}, z, {\bf p}, {\bf q}, {\bf r}, \mathbb{M})}{\beta(\rx, z, {\bf q}, {\bf r})}-\left(\frac{1}{\beta(\rx, z, {\bf q}, {\bf r})}-
  \frac{1}{\beta(\rx, 0, {\bf 0}, {\bf 0})}\right)F_2(\rx, z, {\bf p}, {\bf q})
\end{equation*}
provided that $\beta(\rx, z, {\bf q}, {\bf r})\neq 0$.

(v) For $(\bar{h}_1, \bar{h}_2)$ given by \eqref{def-coefficients3}, define $f_2$ by
\begin{equation*}
\begin{split}
  f_2(\rx, \xi, z, {\bf q}, {\bf r})
  :=&\mcl{H}(S_0+\xi, \Phi_0(\rx)+z-\frac 12|\nabla\vphi_0(\rx)+{\bf q}+{\bf r}^{\perp}|^2)-
  \mcl{H}(S_0, \Phi_0(\rx)-\frac 12|\nabla\vphi_0(\rx)|^2)\\
  &-\bar{h}_1(\rx)z-\bar{h}_2(\rx)q_1+b(\rx)-b_0.
  \end{split}
\end{equation*}

(vi) For $\rx\in \Om_L$, $\xi, \eta, z\in \R$ and ${\bf q}, {\bf r}\in \R^2$ with satisfying $\bar{u}(\rx)+q_1+r_2\neq 0$, define $f_3$ by
\begin{equation*}
  f_3(\rx, \xi, \eta, z, {\bf q}, {\bf r})
  =\frac{-\eta \mcl{H}^{\gam-1}(S_0+\xi, \Phi_0(\rx)+z-\frac 12|\nabla\vphi_0(\rx)+{\bf q}+{\bf r}^{\perp}|^2)}{(\gam-1)(\bar u(x_1)+q_1+r_2)}.
\end{equation*}
\end{definition}

For $a_{j2}(\rm x, z, {\bf q}, {\bf r})$ given by Definition \ref{definition-nlbvp-pert-full}(i), we have an analogy of Lemma \ref{lemma-coeff-reg} (c) as follows:
\begin{lemma}
\label{lemma-coeff-full-new2018}
For $T>0$ and $R>0$, set
\begin{equation*}
\mathfrak{U}^T_R:=\{({\rx}, z, {\bf q}, {\bf r})\in \Om_L\times\R\times \R^2\times \R^2:
\,\, -R<z< R,\,\,|{\bf q}|<R,\,\,|{\bf r}|<R\ \}.
\end{equation*}
For each $\epsilon\in(0, T_{\max})$, there exists an $R_{\epsilon}>0$ depending only on $(\gam, J_0, S_0, \rho_0, E_0, \epsilon)$ so that, whenever $T\in(0, T_{\rm max}-\epsilon]$, the following property holds:
for each $k\in \mathbb{Z}^+$ and $j=1,2$,
\begin{equation}
\label{est-coeff-gen-3-full}
|D^k_{(\rx,z, {\bf q}, {\bf r})}a_{j2}(\rx, z, {\bf q}, {\bf r})|\le C_{k}\quad\tx{in $\ol{\mathfrak{U}_{ R_{\epsilon}}^T}$},
\end{equation}
for $C_k>0$ depending only on $(\gam, J_0, S_0, \rho_0, E_0, \epsilon, k)$
\end{lemma}

Suppose that $\beta(\rx, \Psi, \nabla\psi, \nabla\phi)\neq 0$ and $\bar{u}(\rx)+\psi_{x_1}+\phi_{x_2}\neq 0$ hold so that $\tilde{\mcl{N}}_1$, $\mcl{L}_2$, ${\bf M}$, $(f_1, f_2, f_3)$ are well defined by Definition \ref{definition-nlbvp-pert-full} with
\begin{equation*}
(\xi, \eta, z, {\bf p}, {\bf q}, {\bf r}, \mathbb{M})=(Y,  Y_{x_2}, \Psi, \nabla\Psi, \nabla\psi, \nabla\phi, \nabla(\nabla^{\perp}\phi)).
\end{equation*}
Under this assumption, $(\vphi, \Phi, \phi, S):=(\vphi_0, \Phi_0, 0, S_0)+(\psi, \Psi, \phi, Y)$ solves the nonlinear boundary value problem \eqref{hdecomp-system}--\eqref{eqn4-hdecomp} and \eqref{bc-vphi}--\eqref{bc-T} if and only if $(\psi, \Psi, \phi, Y)$ solves the following nonlinear boundary value problem:
\begin{align}
\label{nlbvp-perturbation-full}
&\begin{cases}
&\tilde{\mcl{N}}_1(\psi, \Psi, \phi)=f_1(\rx,Y, \nabla Y,  \Psi, \nabla\Psi, \nabla\psi, \nabla \phi, \nabla(\nabla^{\perp}\phi))\\
&\mcl{L}_2(\psi, \Psi)=f_2(\rx, Y, \Psi, \nabla\psi, \nabla \phi)\\
&\Delta\phi=f_3(\rx, Y, Y_{x_2}, \Psi, \nabla\psi, \nabla \phi)
\end{cases} \qquad\tx{in\quad $\Om_L$},\\
\label{transport-eqn-full}
&{\bf M}(\rx, Y, \Psi, \nabla\psi, \nabla\phi)\cdot \nabla Y=0 \qquad\tx{in\quad $\Om_L$},
\end{align}
\begin{align}
\label{bc1-psi-pert}
&\psi_{x_1}=u_{\rm en}-u_0-\phi_{x_2}\,\,\text{and}\,\, \psi=\varphi_{\rm en}
\,\,\tx{on}\,\, \Gamen,\quad
\psi_{x_2}=0\,\, \tx{on}\,\, \Lambda_L,\\
\label{bc2-Psi-pert}
&\Psi_{x_1}=E_{\rm en}-E_0\,\,\tx{on}\,\, \Gamen,\quad
\Psi_{x_2}=0\,\,\tx{on}\,\, \Lambda_L,\quad
\Psi=\Phi_{\rm ex}-\Phi_0(L,\cdot)\,\,\tx{on}\,\, \Gamex,\\
\label{bc3-phi-pert}
&\phi_{x_1}=0\,\, \tx{on}\,\, \Gamen,\quad \phi=0\,\, \tx{on}\,\, \der\Om_L\setminus \Gamen,\\
\label{bc4-Y-pert}
&Y=S_{\rm en}-S_0\,\,\tx{on}\,\, \Gamen.
\end{align}

To solve the nonlinear boundary value problem \eqref{nlbvp-perturbation-full}--\eqref{bc4-Y-pert} by the method of iteration, we define the following iteration sets for positive constants $\delta_{\rm e}$, $\delta_{\rm p}$, $\delta_{\rm v}$ and $L$ to be fixed later:
\begin{align}
\label{iterset1-full}
&\mcl{J}^{\rm ent}_{\delta_{\rm e}, L}:=\Bigl\{Y\in H^4(\Om_L):\,\,\|Y\|_{H^4(\Om_L)}\le \delta_{\rm e},\quad Y=S_{\rm en}-S_0\,\,\tx{on $\Gamen$},\\
&\phantom{\mcl{J}^{\rm ent}_{\delta_1, L}:=\{Y\in H^4(\Om_L):\,\,}
\der_{x_2}^kY=0\;
\tx{for $k=1,3$ on $\Lambda_{L}$}\Bigr\},\notag\\
\label{iterset2-full}
&\mcl{J}^{\rm pot}_{\delta_{\rm p}, L}:= \Bigl\{(\psi, \Psi)\in [H^4(\Om_L)]^2: \|(\psi, \Psi)\|_{H^4(\Om_L)}\le \delta_{\rm p},\\ &\phantom{\mcl{J}^{\rm pot}_{\delta_2, L}:=\{(\psi, \Psi)\in [H^4(\Om_L)]^2: }\der_{x_2}^k\psi=\der_{x_2}^k\Psi=0\;
\tx{for $k=1,3$ on $\Lambda_{L}$}\Bigr\},\notag\\
\label{iterset3-full}
&\mcl{J}^{\rm vort}_{\delta_{\rm v}, L}:=\Bigl\{\phi\in H^5(\Om_L):\,\, \|\phi\|_{H^5(\Omega_L)}\leq \delta_{\rm v},\,\,
\der_{x_2}^k\phi=0\;
\tx{for $k=0,2,4$ on $\Lambda_{L}$}\Bigr\},
\end{align}
where the boundary conditions in \eqref{iterset1-full}--\eqref{iterset3-full} are satisfied in the sense of trace.

For fixed functions $\til{Y}\in \mcl{J}^{\rm ent}_{\delta_{\rm e}, L}$, $(\tpsi, \tPsi)\in \mcl{J}^{\rm pot}_{\delta_{\rm p}, L}$, and $\til{\phi}\in \mcl{J}^{\rm vort}_{\delta_{\rm v}, L}$,
denote
\begin{align}
\label{definition-til-as-full}
&\til{a}_{j2}:=a_{j2}(\rx, \tPsi, \nabla\tpsi, \nabla\til{\phi})\quad\tx{for $j=1,2$},\\
\label{definition-til-f1-full}
&\til{f}_1:=f_1(\rx, \til{Y}, \nabla\til{Y}, \tPsi, \nabla\tPsi, \nabla\tpsi, \nabla\til{\phi}, \nabla(\nabla^{\perp}\til{\phi})),\\
\label{definition-til-f2-full}
&\til{f}_2:=f_2(\rx, \til{Y}, \tPsi, \nabla\tpsi, \nabla \til{\phi}),\quad \til{f}_3:=f_3(\rx, \til{Y}, \til{Y}_{x_2}, \tPsi, \nabla\tpsi, \nabla\til{\phi})
\end{align}
where  $(a_{12}, a_{22}, f_1, f_2, f_3)$ are given by Definition \ref{definition-nlbvp-pert-full}. By lengthy but direct computations, one can derive various properties of $(\til a_{12}, \til a_{22}, \til f_1, \til f_2, \til f_3)$ as in the following lemma.
\begin{lemma}
\label{lemma-coeff-full}
For each $\epsilon\in(0, T_{\max})$, there exist positive constants $\delta_2(\epsilon)\in(0, R_{\epsilon}]$ for $R_{\epsilon}$ from Lemma \ref{lemma-coeff-full-new2018} , $\kappa_0>0$, $\kappa_1>0$ and $C>0$ depending only on $(\gam, J_0, S_0, \rho_0, E_0, \epsilon)$ so that if $\max\{\delta_{\rm e},\delta_{\rm p}+\delta_{\rm v}\}\le \delta_2(\epsilon)$ and $L\in(0, T_{\max}-\epsilon]$ hold, then
$(\til{a}_{12}, \til{a}_{22}, \til{f}_1, \til{f}_2, \til{f_3})$ given by \eqref{definition-til-as-full}--\eqref{definition-til-f2-full} for $(\til{Y},\tpsi, \tPsi, \til{\phi})\in \mcl{J}^{\rm ent}_{\delta_{\rm e}, L}\times \mcl{J}^{\rm pot}_{\delta_{\rm p}, L} \times  \mcl{J}^{\rm vort}_{\delta_{\rm v}, L}$
satisfy the following properties:
\begin{itemize}
\item[(a)] for any $\rx\in\ol{\Om}_L$, we have
\begin{equation*}
\begin{split}
&\bar{u}(\rx)+\til{\psi}_{x_1}+\til{\phi}_{x_2}\ge \kappa_0,\\
&{\bf M}(\rx, \til{Y}, \til{\Psi}, \nabla\til{\psi}, \nabla\til{\phi})\cdot {\bf e_1}\ge \kappa_0,\\
&\beta(\rx, \til{\Psi}, \nabla\til{\psi}, \nabla\til{\phi})\le -\kappa_0;
\end{split}
\end{equation*}

\item[(b)]
$\displaystyle{\|\til{a}_{12}\|_{H^3(\Om_L)}
+\|\til{a}_{22}-\ba_{22}\|_{H^3(\Om_L)}\le C(\delta_{\rm p}+\delta_{\rm v})}$;

\item[(c)]
$\displaystyle{\til{a}_{22}\ge \kappa_1}$ in $\ol{\Om_L}$ for the constant $\kappa_1>0$;

\item[(d)] $\displaystyle{\til{a}_{12}=0}$ and $\displaystyle{\der_{x_2}\til{a}_{22}=0}$ hold on $\Lambda_L$;

\item[(e)] $\displaystyle{\|\til f_1\|_{H^3(\Om_L)}\le C \left(\|\til Y\|_{H^4(\Om_L)}
    +\|(\tpsi,\tPsi)\|^2_{H^4(\Om_L)}+\|\til{\phi}\|_{H^5(\Om_L)}\right)}$, and $\displaystyle{\der_{x_2}\til{f}_1=0}$ on $\Lambda_L$;

\item[(f)] $\displaystyle{\|\til{f}_2\|_{H^2(\Om_L)}\le C\left(\|(\til Y,\tpsi,\tPsi)\|^2_{H^4(\Om_L)}
    +\|\til{\phi}\|_{H^5(\Om_L)}+\|b-b_0\|_{C^2(\ol{\Om_L})}\right)}$, and $\displaystyle{\der_{x_2}\til{f}_2=0}$ on $\Lambda_L$ in the sense of trace;
\item[(g)] $\displaystyle{\|\til f_3\|_{H^3(\Om_L)}\le C\|\til Y\|_{H^4(\Om_L)}}$, and $\displaystyle{\der_{x_2}^k\til{f}_3=0}$ on $\Lambda_L$ for $k=0,2$ in the sense of trace.
\end{itemize}
\end{lemma}

Suppose that $\max\{\delta_{\rm e}, \delta_{\rm p}+\delta_{\rm v}\}\le \delta_2(\epsilon)$ and $L\in(0, T_{\max}-\epsilon]$ hold. For fixed $\til{U}:=(\tpsi, \tPsi, \til{\phi})\in \mcl{J}^{\rm pot}_{\delta_{\rm p},L}\times \mcl{J}^{\rm vort}_{\delta_{\rm v},L}$, we define a linear differential operator $\mcl{L}_1^{\til{U}}$ associated with $\til{U}$ by
\begin{equation*}
  \mcl{L}_1^{\til{U}}(\psi, \Psi):=\psi_{x_1x_1}+2\til{a}_{12}\psi_{x_1x_2}
-\til{a}_{22}\psi_{x_2x_2}+\bara_1({\rx})\psi_{x_1}
+\bb_1({\rx})\Psi_{x_1}+\bb_2({\rx})\Psi,
\end{equation*}
where $(\til{a}_{12}, \til{a}_{22})$ are given by \eqref{definition-til-as-full}. The definition of $\mcl{L}_1^{\til{U}}$ is similar to \eqref{def-L1} except that the coefficients $(\til{a}_{12}, \til{a}_{22})$  depend additionally on $\nabla\til{\phi}$ in this case.

For each $(\til{Y}, \til{\psi}, \til{\Psi}, \til{\phi})\in \mcl{J}^{\rm ent}_{\delta_{\rm e},L}\times \mcl{J}^{\rm pot}_{\delta_{\rm p},L}\times \mcl{J}^{\rm vort}_{\delta_{\rm v},L}$, we introduce a linear boundary value problem for $(\psi, \Psi, \phi)$ associated with $(\til{Y}, \til{\psi}, \til{\Psi}, \til{\phi})$ as follows:
\begin{align}
\label{lbvp1-iter-eqn-full}
  &\begin{cases}
  \mcl{L}_1^{\til{U}}(\psi, \Psi)=\til{f}_1\\
  \mcl{L}_2(\psi, \Psi)=\til{f}_2
  \end{cases}\quad\tx{in $\Om_L$},\\
  \label{lbvp1-iter-bc1-psi-full}
&\psi_{x_1}=u_{\rm en}-u_0-\til{\phi}_{x_2}\,\,\text{and}\,\, \psi=\varphi_{\rm en}
\,\,\tx{on}\,\, \Gamen,\quad
\psi_{x_2}=0\,\, \tx{on}\,\, \Lambda_L,\\
\label{lbvp1-iter-bc2-Psi-full}
&\Psi_{x_1}=E_{\rm en}-E_0\,\,\tx{on}\,\, \Gamen,\quad
\Psi_{x_2}=0\,\,\tx{on}\,\, \Lambda_L,\quad
\Psi=\Phi_{\rm ex}-\Phi_0(L,\cdot)\,\,\tx{on}\,\, \Gamex,\\
\label{lbvp1-iter-phi-full}
&\begin{cases}
\Delta\phi=\til{f}_3\quad &\mbox{in $\Om_L$},\\
\phi_{x_1}=0\,\, \tx{on}\,\, \Gamen,\quad \phi=0\,\, \tx{on}\,\, \der\Om_L\setminus \Gamen,
\end{cases}
\end{align}
where $(\til{f}_1, \til{f}_2, \til{f}_3)$ are given by \eqref{definition-til-f1-full}--\eqref{definition-til-f2-full}.

By using Lemma \ref{lemma-coeff-full}, one can make simple adjustments in the proofs of Propositions \ref{lemma-H1-estimate} and  \ref{proposition-wp-lbvp} to achieve the well-posedness of the linear boundary value problem \eqref{lbvp1-iter-eqn-full}--\eqref{lbvp1-iter-bc2-Psi-full}. Furthermore, the unique existence of a solution to \eqref{lbvp1-iter-phi-full} in $H^5$ follows from Lemma \ref{lemma-coeff-full}(g) and standard elliptic theory. So we have the following lemma as an extension of Proposition \ref{proposition-wp-lbvp}.
\begin{lemma}
\label{lemma-wp-lbvp-full}
For any fixed constant $\epsilon_0\in(0, T_{\max})$, there exists a constant $\bar{L}\in(0, T_{\max}-\epsilon_0]$  depending only on $(\gam, S_0, J_0, \rho_0, E_0, \epsilon_0)$ so that whenever $\max\{\delta_{\rm e}, \delta_{\rm p}+\delta_{\rm v}\}\le \frac{\delta_2(\epsilon_0)}{2}$ and $L\in (0, \bar{L}]$, the linear boundary value problem \eqref{lbvp1-iter-eqn-full}--\eqref{lbvp1-iter-phi-full} associated with $(\til{Y}, \tpsi, \tPsi, \til{\phi})\in \mcl{J}^{\rm ent}_{\delta_{\rm e},L} \times
\mcl{J}^{\rm pot}_{\delta_{\rm p},L} \times \mcl{J}^{\rm vort}_{\delta_{\rm v},L}$ has a unique solution $(\psi, \Psi, \phi)\in [H^4(\Om_L)]^2\times H^5(\Om_L)$. And, the solution satisfies the estimate
    \begin{equation}
    \label{estimate-sol-lbvp-full}
      \begin{split}
      &\|(\psi,\Psi)\|_{H^4(\Om_L)}\le C^*
  \Bigl(\|\til{Y}\|_{H^4(\Om_L)}+\|(\tpsi,\tPsi)\|_{H^4(\Om_L)}^2
    +\|\til{\phi}\|_{H^5(\Om_L)}+\sigma_{\rm v}\Bigr),\\
&\|\phi\|_{H^5(\Om_L)}\le C^*\|\til{Y}\|_{H^4(\Om_L)}
      \end{split}
    \end{equation}
    for
    \begin{equation}
    \label{def-sigma-v-new2018}
    \sigma_{\rm v}=\sigma(b,u_{\rm en}, v_{\rm en}, E_{\rm en},\Phi_{\rm ex}, S_{\rm en})
    \end{equation}
    given by \eqref{estimate-bc-full}. In \eqref{estimate-sol-lbvp-full}, the estimate constant $C>0$
can be chosen depending only on $(\gam, J_0, S_0, \rho_0, E_0, \epsilon_0, L)$. Furthermore, the solution $(\psi, \Psi, \phi)$ satisfies the compatibility conditions
 \begin{equation}
 \label{comp-cond-sol-lbvp-full}
 \begin{split}
   &\der_{x_2}^k\psi=\der_{x_2}^k\Psi=\der_{x_2}^{k+1}\phi =0\;
\tx{for $k=1,3$}\tx{ on $\Lambda_L$}
\end{split}
 \end{equation}
 in the sense of trace.
\end{lemma}
The proof of Lemma \ref{lemma-wp-lbvp-full} can be given by adjusting the proof of Proposition \ref{proposition-wp-lbvp} so we skip it.

Now, we outline the iteration scheme that we will apply to prove Theorem \ref{theorem-1} in Section \ref{subsec-proof-theorem1-a}.

\textbf{Step 1.} Assume that $\max\{\delta_{\rm e}, \delta_{\rm p}+\delta_{\rm v}\}\le \frac{\delta_2(\epsilon_0)}{2}$ and $L\in (0, \bar{L}]$ hold. Fix a function $\til{Y}\in \mcl{J}^{\rm ent}_{\delta_{\rm e},L}$. For each $(\tpsi, \tPsi, \til{\phi})$, let $(\psi, \Psi, \phi)$ be the unique solution to the linear boundary value problem \eqref{lbvp1-iter-eqn-full}--\eqref{lbvp1-iter-phi-full} associated with $(\til{Y}, \tpsi, \tPsi, \til{\phi}) \in \mcl{J}^{\rm ent}_{\delta_{\rm e},L}\times \mcl{J}^{\rm pot}_{\delta_{\rm p},L}\times \mcl{J}^{\rm vort}_{\delta_{\rm v},L}$. Hence one can define a mapping $\mathfrak{F}^{\til{Y}}_1: \mcl{J}^{\rm pot}_{\delta_{\rm p},L}\times \mcl{J}^{\rm vort}_{\delta_{\rm v},L}\rightarrow [H^4(\Om_L)]^2\times H^5(\Om_L)$ by
\begin{equation}
\label{definition-mfrakF-Ytil}
  \mathfrak{F}^{\til{Y}}_1(\tpsi, \tPsi, \til{\phi})=(\psi, \Psi, \phi).
\end{equation}
We first express $(\delta_{\rm v}, \delta_{\rm p})$ in terms  of $(\delta_{\rm e}, \sigma_{\rm v})$, and find two positive constants
$\delta^*$ and $\sigma_*$ so that if $\delta_{\rm e}\le \delta^*$ and $\sigma_{\rm  v}\le \sigma_*$, then, for each $\til{Y}\in \mcl{J}^{\rm ent}_{\delta_{\rm e},L}$, the iteration mapping $\mathfrak{F}^{\til{Y}}_1$ has a unique fixed point in $\mcl{J}^{\rm pot}_{\delta_{\rm p},L}\times \mcl{J}^{\rm vort}_{\delta_{\rm v},L}$.

{\textbf{Step 2.}} For each $\til{Y}\in \mcl{J}^{\rm ent}_{\delta_{\rm e},L}$, let $(\psi^{(\til Y)}, \Psi^{(\til Y)}, \phi^{(\til Y)})\in \mcl{J}^{\rm pot}_{\delta_{\rm p},L}\times \mcl{J}^{\rm vort}_{\delta_{\rm v},L}$ be the fixed point of the iteration mapping $\mathfrak{F}^{\til{Y}}_1$. In Appendix \ref{appendixtran}, we show that the linear transport equation ${\bf M}(\rx, \til{Y},  \Psi^{(\til Y)}, \nabla\psi^{(\til Y)}, \nabla\phi^{(\til Y)})\cdot \nabla Y=0$ with the boundary condition \eqref{bc4-Y-pert} has a unique solution $Y\in H^4(\Om_L)$. Then, we define another iteration mapping $\mathfrak{F}_2: \mcl{J}^{\rm ent}_{\delta_{\rm e},L}\rightarrow H^4(\Om_L)$ by
\begin{equation*}
\mathfrak{F}_2(\til{Y})=Y.
\end{equation*}
Next, we express $\delta_{\rm e}$ in terms of $\sigma_{\rm  v}$, and prove that $\mathfrak{F}_2$ has a fixed point $Y_*\in \mcl{J}^{\rm ent}_{\delta_{\rm e},L}$ as long as $\sigma_{\rm v}$ is sufficiently small. Therefore, $(\psi^{(Y_*)}, \Psi^{(Y_*)}, \phi^{(Y_*)}, Y_*)$ solves the nonlinear boundary value problem \eqref{nlbvp-perturbation-full}--\eqref{bc4-Y-pert}.
Finally, by the observation made right after Definition \ref{definition-nlbvp-pert-full}, this proves the existence of a solution to the nonlinear boundary value problem \eqref{hdecomp-system}--\eqref{eqn4-hdecomp} and \eqref{bc-vphi}--\eqref{bc-T}.

{\textbf{Step 3.}} Finally, we prove that if $\sigma_{\rm v}$ is sufficiently small, then $\mathfrak{F}_2$ has a unique fixed point in $\mcl{J}^{\rm ent}_{\delta_{\rm e},L}$. This proves the uniqueness of a solution to the nonlinear boundary value problem \eqref{hdecomp-system}--\eqref{eqn4-hdecomp} and \eqref{bc-vphi}--\eqref{bc-T} so that Theorem \ref{theorem-1} is proved.

\subsubsection{Proof of Theorem \ref{problem-1}}
\label{subsec-proof-theorem1-a}
We give the details of the proof for Theorem \ref{problem-1}  outlined above.

{\bf Step 1. The unique existence of a fixed point of $\mathfrak{F}_1^{\til{Y}}$.} We fix $\epsilon_0\in(0, T_{\max})$, and $L\in(0, \bar{L}]$ for $\bar{L}$ from Lemma \ref{lemma-wp-lbvp-full}.  Assume that
\begin{equation}
\label{fund-assmp-new2018}
\max\{\delta_{\rm e},\delta_{\rm p}+\delta_{\rm v}\}\le \delta_2(\epsilon_0)
\end{equation}
for $\delta_2(\epsilon_0)$ Lemma \ref{lemma-coeff-full}.  Fix $\til{Y}\in \mcl{J}^{\rm ent}_{\delta_{\rm e},L}$, and let $\mathfrak{F}^{\til{Y}}_1: \mcl{J}^{\rm pot}_{\delta_{\rm p},L}\times \mcl{J}^{\rm vort}_{\delta_{\rm v},L}\rightarrow [H^4(\Om_L)]^2\times H^5(\Om_L)$ be given by \eqref{definition-mfrakF-Ytil}.

We recall that $\sigma_{\rm v}=\sigma(b,u_{\rm en}, v_{\rm en}, E_{\rm en},\Phi_{\rm ex}, S_{\rm en})$ is given by \eqref{estimate-bc-full}. We choose the constants $(\delta_{\rm p}, \delta_{\rm v})$ in \eqref{iterset2-full}--\eqref{iterset3-full} in the following forms
\begin{equation}\label{constant-chocie1-full}
\begin{cases}
\delta_{\rm p}=m_1\delta_{\rm e}+m_2\sigma_{\rm v},\\
  \delta_{\rm v}=m_3\delta_{\rm e}
  \end{cases}
\end{equation}
for constants $(m_1, m_2, m_3)$ to be specified later, where $\delta_{\rm e}$ from \eqref{iterset1-full} is to be determined finally in Step 3.

Fix $(\tpsi, \tPsi, \til{\phi})\in  \mcl{J}^{\rm pot}_{\delta_{\rm p},L}\times \mcl{J}^{\rm vort}_{\delta_{\rm v},L}$, and set
\begin{equation*}
(\psi, \Psi, \phi):=\mathfrak{F}^{\til{Y}}_1(\tpsi, \tPsi, \til{\phi}).
\end{equation*}
For the constant $C^*$ from the estimate \eqref{estimate-sol-lbvp-full} given in Lemma \ref{lemma-wp-lbvp-full}, we choose $(m_1, m_2, m_3)$ as
\begin{equation}
\label{m3-choice-new2018}
m_1=12C^*,\quad m_2=4C^*,\quad m_3=2C^*.
\end{equation}
Under these choices, if $(\delta_{\rm e}, \sigma_{\rm v})$ satisfies
\begin{equation}\label{condition-for-dpot-full}
3\delta_{\rm e}+\sigma_{\rm v}\le \frac{1}{12(C^*)^2},
\end{equation}
then it follows from Lemma \ref{lemma-wp-lbvp-full} and \eqref{constant-chocie1-full} that
\begin{equation}
\label{apriori-fixedpt-full}
\begin{split}
  &\|(\psi,\Psi)\|_{H^4(\Om_L)}\le \frac 12 \delta_{\rm p}\quad \text{and}\quad \|\phi\|_{H^5(\Om_L)}\le \frac 12 \delta_{\rm v},
  \end{split}
\end{equation}
and, this implies that the mapping $\mathfrak{F}_1^{\til{Y}}$ maps $\mcl{J}^{\rm pot}_{\delta_{\rm p},L}\times \mcl{J}^{\rm vort}_{\delta_{\rm v},L}$ into itself for any $\til{Y}\in \mcl{J}^{\rm ent}_{\delta_{\rm e},L}$ provided that the condition \eqref{condition-for-dpot-full} holds. Then, one can repeat and adjust the argument in Step 1 of the proof of Theorem \ref{theorem-2-potential} given in Section \ref{subsection-thm-pf1}  to conclude that, for each $\til{Y}\in \mcl{J}^{\rm ent}_{\delta_{\rm e},L}$, the iteration mapping  $\mathfrak{F}_1^{\til{Y}}$  has a fixed point in $\mcl{J}^{\rm pot}_{\delta_{\rm p},L}\times \mcl{J}^{\rm vort}_{\delta_{\rm v},L}$.

We continue to assume that the condition \eqref{condition-for-dpot-full} holds. Later, we will choose $\delta_{\rm e}$ and an upper bound of $\sigma_{\rm  v}$, which becomes $\hat{\sigma}_{\rm bc} $ in Theorem \ref{theorem-1} so that \eqref{condition-for-dpot-full} holds.
For a fixed $\til{Y}\in \mcl{J}^{\rm ent}_{\delta_{\rm e},L}$, let $(\psi^{(i)}, \Psi^{(i)}, \phi^{(i)})$ ($i=1, 2$) be two fixed points of
$\mathfrak{F}_1^{\til{Y}}$ in $\mcl{J}^{\rm pot}_{\delta_{\rm p}, L}\times \mcl{J}^{\rm vort}_{\delta_{\rm v}, L}$, and set
\begin{equation*}
(\check{\psi}, \check{\Psi}, \check{\phi}):=(\psi^{(1)}, \Psi^{(1)}, \phi^{(1)})-(\psi^{(2)}, \Psi^{(2)}, \phi^{(2)}).
\end{equation*}
Similar to Step 2 of the proof of Theorem \ref{theorem-2-potential}, we use Proposition \ref{lemma-H1-estimate}, Definition  \ref{definition-nlbvp-pert-full}, and Lemma \ref{lemma-coeff-full-new2018} to get
\begin{equation}
\label{contraction-1}
  \begin{split}
  &\|(\check{\psi},\check{\Psi})\|_{H^1(\Om_L)}
  \le C^{\natural}\Bigl(\|\check{\phi}\|_{H^2(\Om_L)}  +\bigl(\delta_{\rm e}+\sigma_{\rm v}\bigr) \|(\check{\psi},\check{\Psi})\|_{H^1(\Om_L)}\Bigr),\\
  &\|\check{\phi}\|_{H^2(\Om_L)}\le C^{\natural}\delta_{\rm e}\Bigl(\|(\check{\psi},\check{\Psi})\|_{H^1(\Om_L)}+
  \|\check{\phi}\|_{H^2(\Om_L)}\Bigr)
  \end{split}
\end{equation}
for a constant $C^{\natural}>0$ depending only on $(\gam, J_0, S_0, \rho_0, E_0, \epsilon_0, L)$.
If $(\delta_{\rm e}, \sigma_{\rm v})$ satisfies
\begin{equation}
\label{condition1-for-dent-sig}
 C^{\natural}\bigl(\delta_{\rm e}+\sigma_{\rm v}\bigr)\le \frac 12,
\end{equation}
then it follows from \eqref{contraction-1} that
\begin{equation}
\label{two-inequ-contr}
  \|(\check{\psi},\check{\Psi})\|_{H^1(\Om_L)}\le
   2C^{\natural}\|\check{\phi}\|_{H^2(\Om_L)}\quad\tx{and}\quad
   \|\check{\phi}\|_{H^2(\Om_L)}\le  2C^{\natural}\delta_{\rm e}\|(\check{\psi},\check{\Psi})\|_{H^1(\Om_L)}.
\end{equation}
If the constant $\delta_{\rm e}$ additionally satisfies the condition
\begin{equation}
\label{condition2-for-dent-sig}
  \delta_{\rm e}\le \frac{1}{8(C_{\natural})^2},
\end{equation}
then the estimate \eqref{two-inequ-contr} implies that
\begin{equation*}
{\|(\check{\psi},\check{\Psi})\|_{H^1(\Om_L)}=0},
 \end{equation*}
from which we obtain that $\check{\psi}=\check{\Psi}=\check{\phi}=0$ in $\Om_L$. Thus, we conclude that  $\mathfrak{F}_1^{\til{Y}}$ has a unique fixed point in $\mcl{J}^{\rm pot}_{\delta_{\rm p},L}\times \mcl{J}^{\rm vort}_{\delta_{\rm v},L}$ provided that the conditions \eqref{condition1-for-dent-sig} and \eqref{condition2-for-dent-sig} hold.

{\bf Step 2. The existence of a fixed point of $\mathfrak{F}_2$.} In the previous step, we have shown that, for each $\til{Y}\in \mcl{J}^{\rm ent}_{\delta_{\rm e}, L}$, $\mathfrak{F}^{\til{Y}}_1$ has a unique fixed point in $ \mcl{J}^{\rm pot}_{\delta_{\rm p},L}\times \mcl{J}^{\rm vort}_{\delta_{\rm v},L}$ provided that the conditions \eqref{condition1-for-dent-sig} and \eqref{condition2-for-dent-sig} hold. Note that
the fixed point $(\psi, \Psi, \phi)$ solves the nonlinear boundary value problem \eqref{nlbvp-perturbation-full} with the boundary conditions \eqref{bc1-psi-pert}--\eqref{bc3-phi-pert}, where $Y$ is replaced by $\til{Y}$ in \eqref{nlbvp-perturbation-full}. Since \eqref{nlbvp-perturbation-full} is derived by rewriting \eqref{hdecomp-system} and \eqref{eqn3-hdecomp} in terms of perturbations from the background solution $(\vphi_0, 0, \Phi_0, S_0)$, the vector valued function ${\bf M}(\rx, \til{Y}, \Psi, \nabla\psi, \nabla\phi)$, given by Definition \ref{definition-nlbvp-pert-full}(iii)  satisfies
\begin{equation}
\label{M-div-free}
  \rm{div}\,{\bf M}(\rx, \til{Y}, \Psi, \nabla\psi, \nabla\phi)=0\quad\tx{in $\ol{\Om_L}$}.
\end{equation}
Furthermore, the boundary conditions $\psi_{x_2}=\phi=0$ on $\Lambda_L$ imply that
\begin{equation}
\label{M-slip-on-wall}
  {\bf M}(\rx, \til{Y}, \Psi, \nabla\psi, \nabla\phi)\cdot {\bf n}_{\Lambda}=0\quad\tx{on $\Lambda_L$}
\end{equation}
for the inward unit normal vector ${\bf n}_{\Lambda}$ on $\Lambda_L$.

It follows from Lemma \ref{lemma-coeff-full}(a) and \eqref{fund-assmp-new2018} that one has
\begin{equation}
\label{lwrbd-M-pert-vec}
  {\bf M}(\rx, \til{Y}, \Psi, \nabla\psi, \nabla\phi)\cdot {\bf e}_1\ge \kappa_0\quad\tx{in $\ol{\Om_L}$}.
\end{equation}
Now let us consider the following  boundary value problem for $Y$:
\begin{equation}\label{bvp-for-Y}
 {\bf M}(\rx, \til{Y}, \Psi, \nabla\psi, \nabla\phi)\cdot\nabla Y=0\quad\tx{in $\Om_L$},\qquad Y=S_{\rm en}-S_0\quad\tx{on $\Gamen$}.
\end{equation}
The well-posedness of the boundary value problem \eqref{bvp-for-Y} for $Y$ is stated in the following lemma.
\begin{lemma}\label{lemma-wp-transport-eqn}
Assume that the conditions  \eqref{fund-assmp-new2018}, \eqref{condition-for-dpot-full}, \eqref{condition1-for-dent-sig}, and \eqref{condition2-for-dent-sig} hold. Then, the boundary value problem \eqref{bvp-for-Y} has a unique solution $Y\in H^4(\Om_L)$. And, the solution $Y$ satisfies the following properties:
\begin{itemize}
\item[(a)]
the solution $Y$ can be represented as
\begin{equation*}
  Y(\rx)=(S_{\rm en}-S_0)\circ \mathscr{L}_{\til{Y}}(\rx)\quad\tx{in $\Om_L$}
\end{equation*}
where the mapping $\mathscr{L}_{\til{Y}}:\ol{\Om}_L\rightarrow [-1,1]$ satisfies
    \begin{equation}
    \label{deriv-L-mapping}
      \nabla \mathscr{L}_{\til{Y}}(\rm x)
      =\frac{- \left({\bf M}(\rx, \til{Y}, \Psi, \nabla\psi, \nabla\phi)\right)^{\perp}}{{\bf M}(\rx', \til{Y}(\rx'), \Psi(\rx'), \nabla\psi(\rx'), \nabla\phi(\rx'))|_{\rx'=(0, \mathscr{L}_{\til{Y}}(\rx))}\cdot{\bf e}_1}\quad\tx{in}\quad \Om_L,
    \end{equation}
    furthermore, if a sequence $\{\til{Y}_j\}\subset \mcl{J}^{\rm ent}_{\delta_{\rm e},L}$ converges to $\til{Y}_{\infty}\in \mcl{J}^{\rm ent}_{\delta_{\rm e},L}$ in $H^3(\Om_L)$, then $\{\mathscr{L}_{\til{Y}_j}\}$ converges to $\mathscr{L}_{\til{Y}_{\infty}}$ in $H^4(\Om_L)$;

\item[(b)] there exists a constant $C_{**}>0$ depending only on $(\gam, J_0, S_0, \rho_0, E_0, \epsilon_0, L)$ to satisfy
    \begin{equation}
    \label{estimate-entropy-full}
      \|Y\|_{H^4(\Om_L)}\le C_{**}\|S_{\rm en}-S_0\|_{C^4(\ol{\Gamen})};
    \end{equation}

\item[(c)] the compatibility conditions
\begin{equation*}
  \der_{x_2}^k Y=0\quad\tx{on $\Lambda_L$ for $k=1,3$}
\end{equation*}
hold in the sense of traces.
\end{itemize}
\end{lemma}
This lemma is inspired by the study on the transport equation in \cite{BDX3}, and the  proof of this lemma is similar to the proof of \cite[Lemma 3.3]{BDX3}. For the convenience of the readers,  we provide a proof of Lemma \ref{lemma-wp-transport-eqn} in Appendix \ref{appendixtran}.

Now we choose $\delta_{\rm e}$ in \eqref{iterset1-full} as
\begin{equation}
\label{dent-choice-full}
  \delta_{\rm e}=2C_{**}\sigma_{\rm v}
\end{equation}
for $C_{**}$ from (b) of Lemma \ref{lemma-wp-transport-eqn}. Under this choice of $\delta_{\rm e}$, the solution $Y$ to the boundary value problem \eqref{bvp-for-Y} is contained in $\mcl{J}^{\rm ent}_{\delta_{\rm e},L}$ because $\|S_{\rm en}-S_0\|_{C^4(\ol{\Gamen})}\le \sigma_{\rm v}$ by the definition of $\sigma_{\rm v}$ given in \eqref{def-sigma-v-new2018}.

We define another iteration mapping $\mathfrak{F}_2:\mcl{J}^{\rm ent}_{\delta_{\rm e},L}\rightarrow H^4(\Om_L)$ by
\begin{equation*}
  \mathfrak{F}_2(\til{Y})=Y
\end{equation*}
where $Y$ is the solution to the boundary value problem \eqref{bvp-for-Y}. The mapping $\mathfrak{F}_2$ is well defined and maps $\mcl{J}^{\rm ent}_{\delta_{\rm e},L}$ into itself provided that the conditions \eqref{fund-assmp-new2018}, \eqref{condition-for-dpot-full}, \eqref{condition1-for-dent-sig}, and\eqref{condition2-for-dent-sig} hold under the choice of $\delta_{\rm e}$ by \eqref{dent-choice-full}.
Then, we can apply the Rellich's theorem and the Schauder fixed point theorem to conclude that $\mathfrak{F}_2$ has a fixed point in $\mcl{J}^{\rm ent}_{\delta_{\rm e},L}$.

Let $Y_*\in \mcl{J}^{\rm ent}_{\delta_{\rm e},L}$ be a fixed point of $\mathfrak{F}_2$, and let $(\psi_*, \Psi_*, \phi_*)\in \mcl{J}^{\rm pot}_{\delta_{\rm p},L}\times \mcl{J}^{\rm vort}_{\delta_{\rm v},L}$ be the unique fixed point of $\mathfrak{F}_1^{Y_*}$. Then, $(\psi_*, \Psi_*, \phi_*, Y_*)$ is a solution to the nonlinear boundary value problem \eqref{nlbvp-perturbation-full}--\eqref{bc4-Y-pert}.

Now, we choose $\sigma_{\rm bd}^{(1)}$ so that whenever $\sigma_{\rm v}\in (0, \sigma_{\rm bd}^{(1)}]$, the conditions  \eqref{fund-assmp-new2018}, \eqref{condition-for-dpot-full}, \eqref{condition1-for-dent-sig}, and \eqref{condition2-for-dent-sig} hold under the choices of $(\delta_{\rm p}, \delta_{\rm v}, \delta_{\rm e})$ given by \eqref{constant-chocie1-full}--\eqref{dent-choice-full}.
\begin{itemize}
\item[(i)] By   \eqref{constant-chocie1-full}--\eqref{m3-choice-new2018} and \eqref{dent-choice-full}, the condition \eqref{fund-assmp-new2018} holds if
\begin{equation*}
\sigma_{\rm v}\le \frac{\delta_2(\epsilon_0)}{2C_{**}(1+14C^*)+4C^*}=:\sigma_1
\end{equation*}
for the constants $(\delta_2(\epsilon_0), C^*, C_{**})$ from Lemma \ref{lemma-coeff-full}, \eqref{estimate-sol-lbvp-full}, and Lemma \ref{lemma-wp-transport-eqn}(b), respectively;

\item[(ii)] The condition  \eqref{condition-for-dpot-full} holds if
\begin{equation*}
\sigma_{\rm v}\le \frac{1}{24(C^*)^2(3+2C^*)}=:\sigma_2;
\end{equation*}

\item[(iii)] By  \eqref{dent-choice-full}, the condition  \eqref{condition1-for-dent-sig} holds if
\begin{equation*}
\sigma_{\rm v}\le \frac{1}{2C^{\natural}(2C_{**}+1)}=:\sigma_3
\end{equation*}
for $C^{\natural}$ from \eqref{contraction-1};

\item[(iv)] Finally, the condition \eqref{condition2-for-dent-sig} holds if
\begin{equation*}
\sigma_{\rm v}\le \frac{1}{16(C^{\natural})^2C_{**}}=:\sigma_4.
\end{equation*}
\end{itemize}

Therefore, we choose $\sigma_{\rm bd}^{(1)}$ as
 \begin{equation}
 \label{sigma-bd-step2}
\sigma_{\rm bd}^{(1)} =\min\left\{\sigma_k: k=1,2,3,4
  \right\},
\end{equation}
so that if
\begin{equation}\label{choice1-sigma-full}
\sigma_{\rm v} \leq \sigma_{\rm bd}^{(1)},
  \end{equation}
then all the conditions \eqref{fund-assmp-new2018}, \eqref{condition-for-dpot-full}, \eqref{condition1-for-dent-sig}, and \eqref{condition2-for-dent-sig} hold under the choices of $(\delta_{\rm p}, \delta_{\rm v}, \delta_{\rm e})$ given by \eqref{constant-chocie1-full}--\eqref{m3-choice-new2018}, and \eqref{dent-choice-full}.

This proves the existence of a solution to the nonlinear boundary value problem \eqref{hdecomp-system}--\eqref{eqn4-hdecomp} with boundary conditions \eqref{bc-vphi}--\eqref{bc-T} whenever \eqref{choice1-sigma-full} holds. The estimates and the compatibility conditions \eqref{comp-cond-nlbvp-full}--\eqref{apriori-nlbvp-full-holder} stated in Theorem \ref{theorem-1} can be directly verified by using Lemmas \ref{lemma-coeff-full}--\ref{lemma-wp-transport-eqn},  \eqref{constant-chocie1-full}, \eqref{m3-choice-new2018}, and \eqref{dent-choice-full}.

{\bf Step 3. Uniqueness.}
Under the assumption of \eqref{choice1-sigma-full},
let $(\vphi^{(i)}, \Phi^{(i)}, \phi^{(i)}, S^{(i)})$  ($i=1$, $2$) be two solutions to \eqref{hdecomp-system}--\eqref{eqn4-hdecomp} with boundary conditions \eqref{bc-vphi}--\eqref{bc-T}, and assume that both solutions satisfy \eqref{supersonicity-full} and \eqref{apriori-nlbvp-full}. Then, we set
\begin{equation*}
  (\check{\psi}, \check{\Psi}, \check{\phi}, \hat{Y}):=(\vphi^{(1)}, \Phi^{(1)}, \phi^{(1)}, S^{(1)})-(\vphi^{(2)}, \Phi^{(2)}, \phi^{(2)}, S^{(2)}).
\end{equation*}
Similar to \eqref{contraction-1}, it can be checked that
\begin{equation}
  \label{contraction-final1}
  \begin{split}
  &\|(\check{\psi},\check{\Psi})\|_{H^1(\Om_L)}
\le C\left(\sigma_{\rm v} \|(\check{\psi},\check{\Psi})\|_{H^1(\Om_L)}
+\|\check{\phi}\|_{H^2(\Om_L)}+\|\check{Y}\|_{H^1(\Om_L)}
\right),\\
&\|\check{\phi}\|_{H^2(\Om_L)}
\le C\left(\sigma_{\rm v} \bigr(\|(\check{\psi},\check{\Psi})\|_{H^1(\Om_L)}
+\|\check{\phi}\|_{H^2(\Om_L)}\bigr)+\|\check{Y}\|_{H^1(\Om_L)}\right)
  \end{split}
\end{equation}
for a constant $C>0$ depending only on $(\gam, J_0, S_0, \rho_0, E_0, \epsilon_0, L)$. Any estimate constant $C>0$ appearing hereafter may vary, but it is regarded to depend only on $(\gam, J_0, S_0, \rho_0, E_0, \epsilon_0, L)$ unless otherwise specified.

By Lemma \ref{lemma-wp-transport-eqn}(a), one has
\begin{equation*}
  \check{Y}=(S_{\rm en}-S_0)\circ \mathscr{L}_{1}
  -(S_{\rm en}-S_0)\circ \mathscr{L}_{2},
\end{equation*}
where each $\mathscr{L}_j$ denotes $\mathscr{L}_{Y^{(j)}}$ for $j=1,2$, so we get
\begin{equation}\label{estimate-Y-hat}
\|\check{Y}\|_{H^1(\Om_L)}\le C\sigma_{\rm v}\|\mathscr{L}_{1}-\mathscr{L}_{2}\|_{H^1(\Om_L)}.
\end{equation}
For each $j=1$ and $2$, let ${\bm{\mathcal{M}}}^{(j)}(\rx)=(\mcl{M}_1^{(j)}(\rx), \mcl{M}_2^{(j)}(\rx))$ denote ${\bf M}(\rx, Y^{(j)}, \Psi^{(j)}, \nabla\psi^{(j)}, \nabla\phi^{(j)})$, for ${\bf M}$ defined by Definition  \ref{definition-nlbvp-pert-full}(iii). And, we set
\begin{equation*}
w^{(j)}(x_1,x_2):=\int_{-1}^{x_2} \mcl{M}_1^{(j)}(x_1, y)\,dy.
\end{equation*}
Then one has
$w^{(1)}(0, \mathscr{L}_1(\rx))-w^{(2)}(0, \mathscr{L}_2(\rx))=w^{(1)}(\rx)-w^{(2)}(\rx)
$ in $\Om_L$. It follows from \eqref{lwrbd-M-pert-vec} and  mean-value theorem that one has
\begin{equation*}
  (\mathscr{L}_1-\mathscr{L}_2)(\rx)=
  \frac{(w^{(1)}-w^{(2)})(\rx)-(w^{(1)}-w^{(2)})(0,\mathscr{L}_2(\rx))}
  {\int_0^1 \mcl{M}_1^{(1)}(0, t\mathscr{L}_1(\rx)+(1-t)\mathscr{L}_2(\rx))\,dt}.
\end{equation*}
Direct computations using this expression along with \eqref{apriori-nlbvp-full},  \eqref{lwrbd-M-pert-vec} and \eqref{M1-at-entrance} yield
\begin{equation}\label{L2-estimate-diff-Lmapping}
  \|\mathscr{L}_1-\mathscr{L}_2\|_{L^2(\Om_L)}
  \le  C\left(\sigma_{\rm v} \|\mathscr{L}_1-\mathscr{L}_2\|_{L^2(\Om_L)}+
  \|(\check{Y},\check{\psi},\check{\Psi})\|_{H^1(\Om_L)}
  +\|\check{\phi}\|_{H^2(\Om_L)}
  \right).
\end{equation}
 Furthermore, combining \eqref{apriori-nlbvp-full}, Definition \ref{definition-nlbvp-pert-full}, \eqref{deriv-L-mapping},  \eqref{L2-estimate-diff-Lmapping}, and \eqref{M1-at-entrance} gives
\begin{equation}\label{L2-estimate-diff-Lmapping-deriv}
\|\nabla \mathscr{L}_1-\nabla \mathscr{L}_2\|_{L^2(\Om_L)}
  \le  C\left(\sigma_{\rm v} \|\mathscr{L}_1-\mathscr{L}_2\|_{L^2(\Om_L)}+
  \|(\check{Y},\check{\psi},\check{\Psi})\|_{H^1(\Om_L)} +\|\check{\phi}\|_{H^2(\Om_L)}
  \right).
\end{equation}
Therefore, one has
\begin{equation*}
\|\mathscr{L}_1-\mathscr{L}_2\|_{H^1(\Om_L)}
  \le  C^{\spadesuit}\left(\sigma_{\rm v} \|\mathscr{L}_1-\mathscr{L}_2\|_{L^2(\Om_L)}+
  \|(\check{Y},\hat{\psi},\check{\Psi})\|_{H^1(\Om_L)}
  +\|\check{\phi}\|_{H^2(\Om_L)}
  \right)
\end{equation*}
for some constant $C^{\spadesuit}>0$. So if
\begin{equation*}
\sigma_{\rm v}\le \frac{1}{2C^{\spadesuit}},
\end{equation*}
then we get
\begin{equation*}
\|\mathscr{L}_1-\mathscr{L}_2\|_{H^1(\Om_L)}
  \le  2C^{\spadesuit}\left(
  \|(\check{Y},\check{\psi},\check{\Psi})\|_{H^1(\Om_L)}
  +\|\check{\phi}\|_{H^2(\Om_L)}
  \right)
\end{equation*}

Next, we substitute the previous estimate into the right-hand side of  \eqref{estimate-Y-hat} to get
\begin{equation*}
\|\check{Y}\|_{H^1(\Om_L)}\le C^{\dagger}\sigma_{\rm v}\left(
  \|\check{Y}\|_{H^1(\Om_L)}\|+\|(\check{\psi},\check{\Psi})\|_{H^1(\Om_L)}
  +\|\check{\phi}\|_{H^2(\Om_L)}
  \right)
\end{equation*}
for some constant $C^{\dagger}>0$. Therefore, if $\sigma_{\rm v}$ additionally satisfies
\begin{equation*}
\sigma_{\rm v}\le \frac{1}{2C^{\dagger}},
\end{equation*}
then we get
\begin{equation*}
\|\check{Y}\|_{H^1(\Om_L)}\le C\sigma_{\rm v}\left(
  \|(\check{\psi},\check{\Psi})\|_{H^1(\Om_L)}
  +\|\check{\phi}\|_{H^2(\Om_L)}
  \right)
\end{equation*}
Now, we substitute this estimate into the right-hand sides of the the estimates given in \eqref{contraction-final1} to obtain that
\begin{equation*}
 \|(\check{\psi},\check{\Psi})\|_{H^1(\Om_L)}
  +\|\check{\phi}\|_{H^2(\Om_L)}\le C^{\clubsuit}\sigma_{\rm v}\left( \|(\check{\psi},\check{\Psi})\|_{H^1(\Om_L)}
  +\|\check{\phi}\|_{H^2(\Om_L)}\right)
\end{equation*}
for some constant $C^{\clubsuit}>0$. So if
\begin{equation*}
\sigma_{\rm v}\le \frac{1}{2(C^{\spadesuit}+C^{\dagger}+C^{\clubsuit})},
\end{equation*}
then we finally conclude that
\begin{equation*}
(\vphi^{(1)}, \Phi^{(1)}, \phi^{(1)}, S^{(1)})-(\vphi^{(2)}, \Phi^{(2)}, \phi^{(2)}, S^{(2)})
=(\check{\psi}, \check{\Psi}, \check{\phi}, \check{Y})\equiv 0\quad\tx{in $\Om_L$}.
\end{equation*}

The proof of Theorem \ref{theorem-1} is completed by choosing $ \hat \sigma_{\rm bd}$ as
\begin{equation*}
 \hat \sigma_{\rm bd}=\min\left\{\sigma_{\rm bd}^{(1)}, \frac{1}{2(C^{\spadesuit}+C^{\dagger}+C^{\clubsuit})}\right\},
\end{equation*}
where $\sigma_{\rm bd}^{(1)}$ is given in \eqref{sigma-bd-step2}.
\hfill $\square$

\appendix

\section{Proof of Lemma \ref{lemma-H4-galerkin}}
\label{appendix-1}
In this appendix, we prove Lemma \ref{lemma-H4-galerkin}.
\begin{proof}[Proof of Lemma \ref{lemma-H4-galerkin}]  The proof is divided into 3 steps.
For a fixed $n\in \mathbb{N}$,  let $(\hat{\mfrak{f}}_1^{(n)}, \hat{\mfrak{f}}_2^{(n)})$ be given by \eqref{smooth-approx-in-lbvp}. To simplify notations, we will write $(\hat{\mfrak{f}}_1^{(n)}, \hat{\mfrak{f}}_2^{(n)})$ as $(\hat{\mfrak{f}}_1, \hat{\mfrak{f}}_2)$ hereafter. Note that $(\hat{\mfrak{f}}_1, \hat{\mfrak{f}}_2)$ are smooth in $\ol{\Om_L}$. For the rest of proof, we also fix $m\in \mathbb{N}$, and set
\begin{equation}
\label{definition-f2m-galerkin}
  \hat{\mathfrak{f}}_{l,m}:=\sum_{j=0}^m\langle \hat{\mathfrak{f}}_l, \eta_j \rangle \eta_j\quad\tx{in $\Om_L$ for $l=1,2$},
\end{equation}
where $\{\eta_j\}_{j=0}^\infty$ is the orthonormal basis given in \eqref{eigenfunction}.
And, let $(V_m, W_m)$ given in the form \eqref{def-projection-n} be the solution to \eqref{galerkin-he-system}--\eqref{galerkin-bc-others}.

{\textbf{Step 1. $H^2$ estimate for $W_m$.}}
It follows from the definition of $\mcl{L}_2$ given by \eqref{pfeqn-poisson-modi}, \eqref{def-coefficients3} and \eqref{galerkin-he-system} that
\begin{equation}
\label{galerking-form-global}
  \langle \mcl{L}_2(V_m, W_m)-\hat{\mathfrak{f}}_{2,m}, \eta_k \rangle =0 \quad\tx{for all $k\in \mathbb{Z}_+$, $0<x_1<L$}.
\end{equation}
By  \eqref{galerkin-bc-others}--\eqref{galerkin-slip-bc} and \eqref{galerking-form-global}, $W_m$ becomes a classical solution to the elliptic boundary value problem:
\begin{equation}
\label{lbvp-elliptic-wm}
  \begin{split}
  &\Delta W_m-\bar{h}_1W_m=\hat{\mathfrak{f}}_{2,m}+\bar{h}_2\der_{x_1}V_m \quad\tx{in $\Om_L$},\\
  &\der_{x_1}W_m=0\,\,\tx{on $\Gamen$},\quad
  \der_{x_2}W_m=0\,\,\tx{on $\Lambda_L$},\quad
  W_m=0\,\,\tx{on $\Gam_L$}.
  \end{split}
\end{equation}
Applying \cite[Theorems 8.8 and 8.12]{GilbargTrudinger} and the method of reflection to \eqref{lbvp-elliptic-wm} yields
\begin{equation*}
  \|W_m\|_{H^2(\Om_L)}\le C(\|\hat{\mathfrak{f}}_2\|_{L^2(\Om_L)}
  +\|W_m\|_{L^2(\Om_L)}+\|V_m\|_{H^1(\Om_L)})
\end{equation*}
for a constant $C>0$ depending only on $(\gam, S_0, J_0, \rho_0, E_0, \epsilon_0, L)$.
We combine this estimate with \eqref{estimate-H1-galerkin} to get
\begin{equation}
\label{estimate-H2-for-wm}
      \|W_m\|_{H^2(\Om_L)}\le C\Bigl(\|(\hat{\mathfrak{f}}_1\|_{L^2(\Om_L)}+\|\hat{\mathfrak{f}}_2\|_{L^2(\Om_L)}
      +\|g_1\|_{C^0(\ol{\Gamen})}\Bigr)
    \end{equation}
    for a constant $C>0$ depending only on $(\gam, S_0, J_0, \rho_0, E_0, \epsilon_0, L)$.

{\textbf{Step 2. $H^2$ estimate for $V_m$.}} We divide the proof for the $H^2$-estimate for $V_m$ into three parts.

{\it{Part 1. Energy estimate.}} The first equation in \eqref{galerkin-he-system} can be written as
\begin{equation}
\label{equation-galerkin-1}
  \langle \mcl{L}_1^{(n)}(V_m, W_m),\eta_k\rangle
  =\langle   \hat{\mathfrak{f}}_{1,m}, \eta_k \rangle\quad\tx{for $0<x_1<L$, and $k=0,1,\cdots, m$.}
\end{equation}
We define a linear hyperbolic differential operator $\mcl{L}_{\rm{hyp}}$ by
\begin{equation*}
         \mcl{L}_{\rm{hyp}}(V)
         := V_{x_1x_1}+2\til{a}_{12}^{(n)}V_{x_1x_2}
-\til{a}_{22}^{(n)}V_{x_2x_2}+\bara_1({\rx})V_{x_1},
     \end{equation*}
and use this definition to rewrite \eqref{equation-galerkin-1} as
\begin{equation}
\label{equation-hyperbolic-for-vm}
  \langle \mcl{L}_{\rm{hyp}}(V_m), \eta_k \rangle
  =\langle \hat{\mathfrak{f}}_{1,m}-\bb_1 \der_{x_1}W_m-\bb_2W_m, \eta_k \rangle \quad\tx{for $0<x_1<L$, and $k=0,1,\cdots, m$.}
\end{equation}
   Let us set
   \begin{equation*}
     q_m:=\der_{x_1}V_m\quad\tx{in $\Om_L$.}
   \end{equation*}
We first differentiate \eqref{equation-hyperbolic-for-vm} with respect to $x_1$, then multiply the resultant equation by $\vartheta''_k$ for each $k=0,1,\cdots,m$, and add up the results over $k=0$ to $m$, finally integrate the summation with respect to $x_1$ on the interval $[0, t]$ for $t$ varying in the interval $[0,L]$ to get
\begin{equation}
\label{energy-estimate-H2-galerkin}
  \begin{split}
  \int_{\Om_t}
 \mcl{L}_{\rm{hyp}}(q_m)
  \der_{x_1}q_m
  \,d\rx
  =&\int_{\Om_t}(
  \der_{x_1}f_{1,m}-\der_{x_1}(\bb_1\der_{x_1}W_m+\bb_2 W_m)
  -\der_{x_1}\ba_1q_m) \der_{x_1}q_m\,d\rx\\
  &+\int_{\Om_t} (-\der_{x_1}\til{a}_{12}^{(n)}\der_{x_2}q_m
 + \der_{x_1}\til{a}_{22}^{(n)}\der_{x_2x_2}V_m) \der_{x_1}q_m\,d\rx
  \end{split}
\end{equation}
for $\Om_t:=\{\rx=(x_1,x_2): 0<x_1<t,\,\, -1<x_2<1\}$.
 Using \eqref{galerkin-coeff} gives
\begin{equation*}
\begin{split}
{\tx{LHS of \eqref{energy-estimate-H2-galerkin}}}=&\frac 12 \left(\int_{\Gam_t}-\int_{\Gamen}\right) \left[ (\der_{x_1}q_m)^2+\til{a}_{22}^{(n)}(\der_{x_2}q_m)^2\right]\,dx_2\\
  &+\int_{\Om_t}(\ba_1-\der_{x_2}\til{a}_{12}^{(n)})(\der_{x_1}q_m)^2
  -\frac{\der_{x_1}\til{a}_{12}^{(n)}}{2}(\der_{x_2}q_m)^2
  +\der_{x_2}\til{a}_{22}^{(n)}\der_{x_1}q_m\der_{x_2}q_m
  \,d\rx
  \end{split}
\end{equation*}
for $\Gam_t=\{(t,x_2)\in\R^2: -1<x_2<1\}$. By Lemmas \ref{lemma-bgdcoeff} and \ref{lemma-coeff-regularity2}, \eqref{galerkin-coeff}, Morrey's inequality and Cauchy-Schwarz inequality, there exist positive constants $\lambda$, $\mu$, and $C$ depending only on $(\gam, J_0, S_0, \rho_0, E_0, \epsilon_0)$ to satisfy
\begin{equation*}
  \int_{\Om_t}
  \mcl{L}_{\rm{hyp}}(q_m)
  \der_{x_1}q_m\,d\rx
  \ge \lambda\int_{\Gam_t} |\nabla_{\rx}q_m|^2\,dx_2
  -\mu\int_{\Gamen} |\nabla_{\rx}q_m|^2\,dx_2
  -C\int_{\Om_t} |\nabla_{\rx}q_m|^2\,d\rx.
\end{equation*}
Substituting this inequality into the left-hand side of \eqref{energy-estimate-H2-galerkin} and applying \eqref{estimate-H2-for-wm} and Cauchy-Schwarz inequality yield
\begin{equation}
\label{gronwall-expression1}
\begin{split}
  \int_{\Gam_t} |\nabla_{\rx}q_m|^2\,dx_2\le
  \frac{\mu}{\lambda}\int_{\Gamen} |\nabla_{\rx}q_m|^2\,dx_2 +C\Big(&\int_{\Om_t}|\nabla_{\rx}q_m|^2+(\der_{x_2}^2V_m)^2\,d\rx\\
 & \quad  +(\|\hat{\mathfrak{f}}_1\|_{H^1(\Om_L)}
      +\|\hat{\mathfrak{f}}_2\|_{L^2(\Om_L)}
      +\|g_1\|_{C^0(\ol{\Gamen})})^2\Big)
      \end{split}
\end{equation}
for some constant $C>0$ depending only on $(\gam, J_0, S_0, \rho_0, E_0, \epsilon_0)$.
Next, we estimate $\int_{\Gamen} |\nabla_{\rx}q_m|^2\,dx_2$ and $\int_{\Om_t}(\der_{x_2}^2V_m)^2\,d\rx$, separately.

{\it{Part 2. Estimate of $\int_{\Gamen} |\nabla_{\rx}q_m|^2\,dx_2$.}}
We differentiate the boundary condition
$\der_{x_1}V_m=\sum_{j=0}^m \langle g_1, \eta_j\rangle \eta_j$
on $\Gamen$ with respect to $x_2$ to get
\begin{equation}
\label{H2-galerkin1}
  \int_{\Gamen} (\der_{x_2}q_m)^2\,dx_2
  =\int_{\Gamen}\Bigl(\sum_{j=0}^m\langle g_1, \eta_j\rangle \eta_j'(x_2)\Bigr)^2\,dx_2.
\end{equation}
By \eqref{eigenfunction}, for $0\le j,k\le m$, one has
\begin{equation}
\label{H2-galerkin2}
 \int_{\Gamen}\langle g_1, \eta_j\rangle \langle g_1, \eta_k\rangle \eta_j'\eta_k' dx_2=\begin{cases}
 \langle g_1, \eta_j\rangle^2(j\pi)^2 \quad&\mbox{for $j=k$},\\
 0\quad&\mbox{otherwise}.
 \end{cases}
\end{equation}
For any $j\in \mathbb{N}$, since $-\eta_j''=(j\pi)^2\eta_j$, integration by parts with using $\eta_j'(\pm 1)=0$ yields
\begin{equation}
\label{H2-galerkin3}
  \langle g_1, \eta_j\rangle=\frac{1}{(j\pi)^2}\langle g_1, -\eta_j''\rangle =\frac{1}{j\pi}\langle g_1', \frac{\eta_j'}{j\pi}\rangle.
\end{equation}
Note that the set $\{\frac{\eta_j'}{j\pi}\}_{j=1}^\infty$ forms an orthonormal basis in $L^2([-1,1])$. Therefore, we conclude from \eqref{H2-galerkin1}--\eqref{H2-galerkin3} that
\begin{equation}
\label{H2-galerkin4}
    \int_{\Gamen} (\der_{x_2}q_m)^2\,dx_2\le \int_{\Gamen}|g_1'|^2\,dx_2\le 2\|g_1\|^2_{C^1(\ol{\Gamen})}.
\end{equation}

For each $k=0, 1,\cdots, m$, multiplying \eqref{equation-hyperbolic-for-vm} by $\vartheta_k''$, summing up over $k=0$ to $m$, and integrating the result over $\Gamen$ with respect to $x_2$  give
\begin{equation}
\label{integral-qm-x1}
  \int_{\Gamen}(\der_{x_1}q_m)^2 +2\til{a}_{12}^{(n)}\der_{x_2}q_m\der_{x_1}q_m+\ba_1q_m\der_{x_1}q_m\,dx_2
  =\int_{\Gamen} (\hat{f}_{1,m}-\bb_2 W_m)\der_{x_1}q_m\,dx_2
\end{equation}
because $V_m=\der_{x_1}W_m=0$ on $\Gamen$ due to \eqref{galerkin-bc-others}. It follows from \eqref{galerkin-bc-others}, \eqref{estimate-H1-galerkin}, \eqref{estimate-H2-for-wm}, \eqref{H2-galerkin2}, Cauchy-Schwarz inequality and trace inequality, and \eqref{integral-qm-x1} that
\begin{equation*}
  \int_{\Gamen}(\der_{x_1}q_m)^2\,dx_2
  \le C\Bigl(\|\hat{f}_1\|_{H^1(\Om_L)}
  +\|\hat{f}_2\|_{L^2(\Om_L)}+\|g_1\|_{C^1(\ol{\Gamen})}\Bigr)^2,
\end{equation*}
where the constant $C>0$ depends only on $(\gam, J_0, S_0, \rho_0, E_0, \epsilon_0, L)$.
We combine this integral estimate with \eqref{H2-galerkin4} to finally get
\begin{equation}
\label{integral-qm-ent}
  \int_{\Gamen}|\nabla_{\rx}q_m|^2\,dx_2
  \le C\Bigl(\|\hat{f}_1^{(n)}\|_{H^1(\Om_L)}
  +\|\hat{f}_2^{(n)}\|_{L^2(\Om_L)}+\|g_1\|_{C^1(\ol{\Gamen})}\Bigr)^2,
\end{equation}
where the constant $C>0$ depends only on $(\gam, J_0, S_0, \rho_0, E_0, \epsilon_0, L)$.

{\it{Part 3. Estimate of $\int_{\Om_t}(\der_{x_2}^2V_m)^2\,d\rx$.}}
We rewrite \eqref{equation-hyperbolic-for-vm} as
\begin{equation}
\label{galerkin-v-2nd-derib}
  \langle \til{a}_{22}^{(n)}\der_{x_2}^2V_m, \eta_k\rangle
  =\langle \der_{x_1}q_m+2\til{a}_{12}^{(n)}\der_{x_2}q_m+\ba_1q_m
  +\bb_1\der_{x_1}W_m+\bb_2 W_m-\hat{\mathfrak{f}}_{1,m},\eta_k\rangle
\end{equation}
for $x_1\in (0,L)$, $k=0,1,\cdots, m$. Since $\eta_k''=-(k\pi)^2\eta_k$ for each $k\in \mathbb{Z}^+$, it follows from \eqref{galerkin-v-2nd-derib} that
\begin{equation*}
  \int_0^t\sum_{k=0}^m\vartheta_k\langle \til{a}_{22}^{(n)}\der_{x_2}^2V_m, \eta_k''\rangle dx_1
  =\int_0^t\sum_{k=0}^m \vartheta_k\langle \der_{x_1}q_m+2\til{a}_{12}^{(n)}\der_{x_2}q_m+\ba_1q_m
  +\bb_1\der_{x_1}W_m+\bb_2 W_m-\hat{\mathfrak{f}}_{1,m},\eta_k''\rangle dx_1,
\end{equation*}
which is the same as
\begin{equation}
\label{H2-estimate-v-2nd-deriv-expr}
\begin{split}
  &\int_{\Om_t}\til{a}_{22}^{(n)}(\der_{x_2}^2V_m)^2\,d\rx\\
  =\,\,&\int_{\Om_t}(\der_{x_1}q_m+2\til{a}_{12}^{(n)}\der_{x_2}q_m+\ba_1q_m
  +\bb_1\der_{x_1}W_m+\bb_2 W_m-\hat{\mathfrak{f}}_{1,m})\der_{x_2}^2V_m\,d\rx.
  \end{split}
\end{equation}
 Combining Lemmas \ref{lemma-bgdcoeff} and \ref{lemma-coeff-regularity2}, \eqref{galerkin-coeff}, Morrey's inequality and Cauchy-Schwarz inequality, \eqref{estimate-H2-for-wm}, and  \eqref{H2-estimate-v-2nd-deriv-expr} yields that
\begin{equation}\label{H2-estimate-v-2nd-deriv-final}
  \int_{\Om_t} (\der_{x_2}^2V_m)^2\,d\rx
  \le C \left(\int_{\Om_t} |\nabla_{\rx}q_m|^2\,d\rx+(\|\hat{\mathfrak{f}}_1\|_{L^2(\Om_L)}+\|\hat{\mathfrak{f}}_2\|_{L^2(\Om_L)}
      +\|g_1\|_{C^0(\ol{\Gamen})})^2\right)
\end{equation}
for some constant $C>0$ depending only on $(\gam, J_0, S_0, \rho_0, E_0, \epsilon_0,L)$.

For notational convenience, we define
\begin{equation*}
  \mathscr{Z}_m(t):=\int_{\Om_t} |\nabla_{\rx}q_m|^2\,d\rx\quad \text{for}\,\, t\in [0, L]
\end{equation*}
and
\begin{equation}
\label{definition-omega-nh}
  \mathscr{E}(\hat{f}_1, \hat{f}_2, g_1):=
  \Bigl(\|\hat{f}_1\|_{H^1(\Om_L)}+\|\hat{f}_2\|_{L^2(\Om_L)}
  +\|g_1\|_{C^1(\ol{\Gamen})}\Bigr)^2.
\end{equation}
It follows from \eqref{gronwall-expression1}, \eqref{integral-qm-ent}, and \eqref{H2-estimate-v-2nd-deriv-final} that  $ \mathscr{Z}_m$ satisfies a differential inequality
\begin{equation}
  \label{H2-diff-ineq}
  \mathscr{Z}_m'(t)\le \alp \mathscr{Z}_m(t)+\beta \mathscr{E}(\hat{f}_1, \hat{f}_2, g_1)\quad\tx{for $0<t<L$},
\end{equation}
where the constants $\alp$ and  $\beta$ depend only on $(\gam, J_0, S_0, \rho_0, \epsilon_0, L)$.
Applying Gronwall's inequality to \eqref{H2-diff-ineq} gives
\begin{equation}\label{H2-estimate-vm-1}
  \mathscr{Z}_m(L)\le C \mathscr{E}(\hat{f}_1, \hat{f}_2, g_1).
\end{equation}
Finally, the estimate \eqref{H2-estimate-vm-1}, together with \eqref{H2-estimate-v-2nd-deriv-expr}, yields
\begin{equation}\label{H2-estimate-vm-final}
  \|V_m\|_{H^2(\Om_L)}\le C\Bigl(\|\hat{f}_1\|_{H^1(\Om_L)}+\|\hat{f}_2\|_{L^2(\Om_L)}
  +\|g_1\|_{C^1(\ol{\Gamen})}\Bigr).
\end{equation}
In \eqref{H2-estimate-vm-1}--\eqref{H2-estimate-vm-final}, the constants $C$ may vary, but they depend only on $(\gam, J_0, S_0, \rho_0, \epsilon_0, L)$.

{\textbf{Step 3. Estimate for higher order weak derivatives of $(V_m, W_m)$.}} In order to complete a priori $H^4$ estimates of $(V_m, W_m)$ in $\Om_L$, we continue the bootstrap argument. All the details can be given by employing the ideas in Steps 1 and 2, but they are much more lengthy and technical. So, in this step, we only describe main differences in establishing the estimate for higher order weak derivatives of $(V_m, W_m)$ in $\Om_L$.

Note that the extension of $\hat{\mathfrak{f}}_{2,m}+\bar{h}_2\der_{x_1}V_m$ given by even reflection about $\Lambda_L$ is $H^1$ across $\Lambda_L$ without any additional compatibility condition. Therefore, back to \eqref{lbvp-elliptic-wm}, we apply Lemma \ref{lemma-H1-estimate}, \eqref{H2-estimate-vm-final} and the method of reflection to obtain that
\begin{equation}
\label{estimate-H3-for-wm}
  \|W_m\|_{H^3(\Om_L)}\le C\Bigl(\|\hat{\mathfrak{f}}_2\|_{H^1(\Om_L)}+\sqrt{\mathscr{E}(\hat{\mfrak f}_1, \hat{\mfrak f}_2, g_1)}\Bigr),
\end{equation}
where the constant $C>0$ depends only on $(\gam, S_0, J_0, \rho_0, E_0, \epsilon_0, L)$ and $\mathscr{E}(\hat{\mfrak f}_1, \hat{\mfrak f}_2, g_1)$ is given in \eqref{definition-omega-nh}.

For a priori $H^3$ estimate of $V_m$, we adapt the argument in Step 2. The main difference is that the compatibility condition
\begin{equation}
\label{compatibility-cond-g1}
  \frac{dg_1}{dx_2}(\pm 1)=0
\end{equation}
derived from \eqref{compatibility-condition-BC1} is used. For example, when we estimate
\begin{equation*}
  \int_{\Gamen} (\der_{x_2}^3V_m )^2\,dx_2=\int_{\Gamen} \Bigl(\sum_{j=0}^m \langle g_1, \eta_j\rangle \eta_j''\Bigr)^2\,dx_2=
  \int_{\Gamen} \Bigl(\sum_{j=0}^m \langle g_1, \eta_j''\rangle \eta_j \Bigr)^2\,dx_2,
\end{equation*}
to derive a differential inequality similar to \eqref{H2-diff-ineq}, we obtain from \eqref{compatibility-cond-g1}  that
\begin{equation*}
  \langle g_1, \eta_j'' \rangle =\langle g_1'', \eta_j \rangle,
\end{equation*}
from which one has
\begin{equation*}
  \int_{\Gamen} (\der_{x_2}^3V_m )^2\,dx_2 \le 2\Bigl(\|g_1\|_{C^2(\ol{\Gamen})}\Bigr)^2.
\end{equation*}
Furthermore, it can be directly checked from \eqref{def-f1}, \eqref{comp-cond1}, \eqref{distance-approx-origina} and \eqref{smooth-approx-in-lbvp} that
\begin{equation*}
  \der_{x_2}\hat{\mathfrak{f}}_{1}
 =0\quad\tx{on $\Lambda_L$}.
\end{equation*}
We use this to derive that
\begin{equation*}
  \|\hat{\mathfrak{f}}_{1,m}\|_{H^2(\Om_L)}
  \le \|\hat{\mathfrak{f}}_1\|_{H^2(\Om_L)}\quad\tx{for all $m\in \mathbb{Z}^+$.}
\end{equation*}
With the aid of  \eqref{H2-estimate-vm-final}-- \eqref{compatibility-cond-g1}, it follows from lengthy but straightforward computations that
\begin{equation}
  \label{estimate-H3-vm-final}
  \|V_m\|_{H^3(\Om_L)}\le C\Bigl(\|\hat{f}_1\|_{H^2(\Om_L)}
  +\|\hat{f}_2\|_{H^1(\Om_L)}
  +\|g_1\|_{C^2(\ol{\Gamen})}\Bigr),
\end{equation}
where the constant $C>0$ depending only on $(\gam, S_0, J_0, \rho_0, E_0, \epsilon_0, L)$.

Similar to $H^3$ estimate of $W_m$ given in \eqref{estimate-H3-for-wm}, a priori $H^4$-estimate of $W_m$ can be obtained by applying \eqref{estimate-H3-vm-final}, \cite[Theorems 8.8 and 8.12]{GilbargTrudinger} and the method of reflection to \eqref{lbvp-elliptic-wm} because the compatibility condition $\der_{x_2}(\hat{\mathfrak{f}}_{2,m}+\bar{h}_2\der_{x_1}V_m)=0$ holds on $\Lambda_{L}$. Furthermore, it can be directly checked from \eqref{comp-cond1}, \eqref{def-f2}, \eqref{distance-approx-origina} and \eqref{smooth-approx-in-lbvp}  that one has
\begin{equation*}
  \der_{x_2}\hat{\mathfrak{f}}_{2}=0\quad\tx{on $\Lambda_L$}.
\end{equation*}
Thus we obtain that
\begin{equation*}
  \|\hat{\mathfrak{f}}_{2,m}\|_{H^2(\Om_L)}
  \le \|\hat{\mathfrak{f}}_2\|_{H^2(\Om_L)}\quad\tx{for all $m\in \mathbb{Z}^+$.}
\end{equation*}
Therefore, we have
\begin{equation}
\label{estimate-H4-for-wm}
  \|W_m\|_{H^4(\Om_L)}\le
  C\Bigl(\|\hat{\mathfrak{f}}_1\|_{H^2(\Om_L)} +\hat{\mathfrak{f}}_2\|_{H^2(\Om_L)}
+\|g_1\|_{C^2(\ol{\Gamen})}\Bigr).
\end{equation}
Finally, adapting the argument in Step 2 and using the estimate \eqref{estimate-H4-for-wm} yield
\begin{equation}\label{estimate-H4-for-psi}
  \|V_m\|_{H^4(\Om_L)}\le C\Bigl(\|\hat{\mathfrak{f}}_1\|_{H^3(\Om_L)}+
\|\hat{\mathfrak{f}}_2\|_{H^2(\Om_L)}
+\|g_1\|_{C^3(\ol{\Gamen})} \Bigr ).
\end{equation}
This finishes the proof of Lemma \ref{lemma-H4-galerkin}.
\end{proof}

\section{Proof of Lemma \ref{lemma-wp-transport-eqn}}\label{appendixtran}
\begin{proof}[Proof of Lemma \ref{lemma-wp-transport-eqn}]
For a fixed $\til{Y}\in \mcl{J}^{\rm ent}_{\delta_{\rm e}, L}$, denote ${\bm{\mathcal{M}}}(\rx)=(\mathcal{M}_1, \mathcal{M}_2)(\rx)$ by
\begin{equation*}
  {\bm{\mathcal{M}}}(\rx):={\bf M}(\rx, \til{Y}, \Psi, \nabla\psi, \nabla\phi)
\end{equation*}
for  ${\bf M}(\rx, \til{Y}, \Psi, \nabla\psi, \nabla\phi)$ defined by Definition  \ref{definition-nlbvp-pert-full}(iii).
We define a function $w:\ol{\Om_L}\rightarrow \R$ by
\begin{equation*}
  w(x_1,x_2):=\int_{-1}^{x_2} \mathcal{M}_1(x_1, y)\,dy.
\end{equation*}
By adjusting the proof of \cite[Lemma 3.3]{BDX3} with using \eqref{compatibility-condition-full}, \eqref{M-div-free}, \eqref{M-slip-on-wall} and Lemma \ref{lemma-coeff-full}(a) , we obtain the following properties:
\begin{itemize}
\item[(i)] $\displaystyle{\nabla^{\perp} w=\bm{\mathcal{M}}}$ in $\ol{\Om_L}$;

\item[(ii)] $\displaystyle{0=w(0,-1)\le w(\rx)\le w(0,1)}$ in $\ol{\Om_L}$;

\item[(iii)] The function $w_0:=w(0,\cdot):[-1,1]\rightarrow [0, w(0,1)]$ is strictly increasing, and its inverse $w_0^{-1}: [0, w(0,1)] \rightarrow [-1,1]$ is well defined;

\item[(iv)] We define a Lagrangian coordinate mapping $\mathscr{L}_{\til{Y}}:\ol{\Om_L}\rightarrow [-1,1]$ by
\begin{equation}
\label{expression-L-mapping}
  \mathscr{L}_{\til{Y}}(\rx):=w_0^{-1}\circ w(\rx).
\end{equation}
From the definition of $\mathscr{L}_{\til{Y}}$, \eqref{deriv-L-mapping} stated in Lemma \ref{lemma-wp-transport-eqn}(a)   can be directly checked. And, the function $Y$ given by
\begin{equation*}
  Y=(S_{\rm en}-S_0)\circ \mathscr{L}_{\til{Y}}
\end{equation*}
solves the boundary value problem \eqref{bvp-for-Y}.
\end{itemize}
In order to complete the proof of Lemma \ref{lemma-wp-transport-eqn}, it suffices to show that there exists a constant $C_{**}>0$ depending only on $(\gam, J_0, S_0, \rho_0, E_0, \epsilon_0, L)$ to satisfy
\begin{equation}
\label{estimate-Lagrangian-map}
  \|\mathscr{L}_{\til{Y}}\|_{H^4(\Om_L)}\le C_{**}\quad\tx{for all $\til{Y}\in \mcl{J}^{\rm ent}_{\delta_{\rm e},L}$}.
\end{equation}
Once \eqref{estimate-Lagrangian-map} is verified, then the rest of Lemma \ref{lemma-wp-transport-eqn} can be easily proved by direct computations. In particular, the continuity of $\mathscr{L}_{\til{Y}}$ with respect to $\til{Y}$ in $H^3(\Om_L)$ stated in Lemma \ref{lemma-wp-transport-eqn}(a) follows from the smooth dependence of ${\bf M}$ on $(\rx, \til{Y},\Psi, \nabla\psi, \nabla\phi)$, and the continuous dependence of the fixed point  $(\psi, \Psi, \phi)\in \mcl{J}^{\rm pot}_{\delta_{\rm p},L}\times \mcl{J}^{\rm vort}_{\delta_{\rm v},L}$ of the iteration mapping $\mathfrak{F}_1^{\til Y}$ on $\til{Y}\in \mcl{J}^{\rm ent}_{\delta_{\rm e},L}$.

The rest of the proof devotes to verify \eqref{estimate-Lagrangian-map}. Due to the smooth dependence of ${\bf M}$ on $(\rx, \til{Y}, \Psi, \nabla\psi, \nabla\phi)$, there exists a constant $C>0$ depending only on $(\gam, J_0, S_0, \rho_0, E_0, \epsilon_0, L)$ to satisfy
\begin{equation}
\label{estimate-vector-M}
  \|{\bf M}(\cdot, \til{Y}, \Psi, \nabla\psi, \nabla\phi)\|_{H^3(\Om_L)}\le C\quad\tx{for all $\til{Y}\in \mcl{J}^{\rm ent}_{\delta_{\rm e},L}$}.
\end{equation}
 It follows from \eqref{rho-hdecomp}, \eqref{bc1-psi-pert}, \eqref{bc3-phi-pert}, \eqref{iterset1-full} , and Definition \ref{definition-nlbvp-pert-full}(iii) that one has
\begin{equation}
\label{M1-at-entrance}
  {\bf M}(\rx, \til{Y},\Psi, \nabla\psi, \nabla\phi)\cdot{\bf e}_1=
  \left(\frac{\gam-1}{\gam S_{\rm en}}\left(\Phi_0+\Psi-\frac 12(u_{\rm en}^2+v_{\rm en}^2)\right)\right)^{\frac{1}{\gam-1}}u_{\rm en}\quad \tx{on $\Gamen$}.
\end{equation}
Note that  the fixed point  $(\psi, \Psi, \phi)\in \mcl{J}^{\rm pot}_{\delta_{\rm p},L}\times \mcl{J}^{\rm vort}_{\delta_{\rm v},L}$ of the iteration mapping $\mathfrak{F}_1^{\til Y}$ solves the nonlinear boundary value problem \eqref{nlbvp-perturbation-full} with boundary conditions \eqref{bc1-psi-pert}--\eqref{bc3-phi-pert} with $Y$ being replaced by $\til{Y}$. So $\Psi$ can be considered as a solution to the linear boundary value problem:
\begin{equation}\label{lbvp-Psi-H5}
\begin{split}
&\Delta\Psi=F_2\quad\tx{in $\Om_L$}, \\
&\Psi_{x_1}=E_{\rm en}-E_0\,\,\tx{on $\Gamen$},\quad \Psi_{x_2}=0\,\,\tx{on $\Lambda_L$},
\end{split}
\end{equation}
where
\begin{equation*}
F_2:=\bar{h}_1\Psi+\bar{h}_2\psi_{x_1}+f_2(\rx, \til{Y}, \Psi, \nabla\psi, \nabla\phi)
 \end{equation*}
with $(\bar h_1, \bar h_2)$ given by \eqref{def-coefficients3}, and $f_2$ given by Definition \ref{definition-nlbvp-pert-full}(v).
Using \eqref{constant-chocie1-full} and \eqref{apriori-fixedpt-full} yields
\begin{equation}
\label{estimate-F2-trans}
  \|F_2\|_{H^3(\Om_L)}\le C(\delta_{\rm e}+\sigma_{\rm v})
\end{equation}
for $\sigma_{\rm v}$ given by \eqref{def-sigma-v-new2018}.
From \eqref{bc1-psi-pert}, \eqref{bc2-Psi-pert} and Definition \ref{definition-nlbvp-pert-full}, it can be directly checked that
\begin{equation}
\label{F2-compatibility}
\der_{x_2}F_2=0\quad\tx{on $\Lambda_L$}.
\end{equation}
Applying the Morrey's inequality to $F_2$ to obtain from \eqref{estimate-F2-trans} and \eqref{F2-compatibility} that
\begin{equation*}
\|F_2\|_{C^{1,\frac 12}(\ol{\Om_L})}\le C(\delta_{\rm e}+\sigma_{\rm v}).
\end{equation*}
In addition, one has $\der_{x_2}(E_{\rm en}-E_0)(\pm 1)=0$ due to \eqref{compatibility-condition-full}. Then, by the standard Schauder estimates and the method of reflection, we obtain that
\begin{equation}
\label{holder-est-Psi-full}
\|\Psi\|_{C^{3, \frac 12}(\ol{\Om_L\cap \{x_1<\frac L2\}})}\le
C(\delta_{\rm e}+\sigma_{\rm v}).
\end{equation}
The estimate constants $C$ appeared so far vary, but  they all depend only on $(\gam, J_0, S_0, \rho_0, E_0, \epsilon_0, L)$.
It follows from \eqref{M1-at-entrance} and \eqref{holder-est-Psi-full} that $ {\bf M}(\rx, \til{Y},\Psi, \nabla\psi, \nabla\phi)\cdot{\bf e}_1\in C^3(\ol{\Gamen})$.
Hence straightforward computations with using \eqref{lwrbd-M-pert-vec}, \eqref{deriv-L-mapping}, \eqref{estimate-vector-M}, \eqref{M1-at-entrance}, \eqref{holder-est-Psi-full}, the chain rule and the Sobolev inequality show that there exists a constant $C_{**}>0$ depending only on $(\gam, J_0, S_0, \rho_0, E_0, \epsilon_0, L)$ to satisfy the estimate \eqref{estimate-Lagrangian-map}. This finishes the proof of Lemma \ref{lemma-wp-transport-eqn}.
\end{proof}

\medskip

{\bf Acknowledgement.}
The research of Myoungjean Bae was supported in part by  Samsung Science and Technology Foundation
under Project Number SSTF-BA1502-02. The research of Ben Duan was supported in part by NSFC No.11871133, No.11671412, the Fundamental Research Funds for the Central Universities grant DUT18RC(3)000 and the High-level innovative and entrepreneurial talents support plan in Dalian grant 2017RQ041.
The research of Chunjing Xie was supported in part by  NSFC grants and 11631008, 11422105, and 11511140276,  and Young Changjiang Scholar of Ministry of Education in China. The authors would like to thank the hospitalities and support for many visits in Korea Institute for Advanced Study, Pohang University of Science and Technology, and Shanghai Jiao Tong University.

\medskip

Conflict of Interest: The authors declare that they have no conflict of interest.


\begin{thebibliography}{99}




\bibitem{Anderson-Monc}
L. Anderson, and V. Moncrief,
{\it Elliptic hyperbolic systems and the Einstein equations},
Ann. H. Poincare, {\textbf 4} (2003), 1--34.

\bibitem{MarkPhase}
U. M. Ascher, P. A. Markowich, P. Pietra, and C.
Schmeiser, {\it A phase plane analysis of transonic solutions for
the hydrodynamic semiconductor model}, Math. Models Methods Appl.
Sci., {\bf 1} (1991),  347--376.

\bibitem{BDX}
M. Bae, B. Duan, and C. J. Xie, {\it Existence and stability of multidimensional steady potential flows for Euler-Poisson equations},  Arch. Ration. Mech. Anal., {\bf 220} (2016)  155--191.

\bibitem{BDX3}
M. Bae, B. Duan, and C. J. Xie, {\it Subsonic solutions for steady Euler-Poisson system
in two-dimensional nozzles}, SIAM J. Math. Anal., {\textbf{46}} (2014), 3455--3480.

\bibitem{BDX2}
M. Bae, B. Duan, and C. J. Xie, {\it   Two dimensional subsonic flows with self-gravitation in bounded domain}, Math. Models Methods Appl.
Sci.,  {\textbf 25} (2015), 2721--2747.



\bibitem{Shu}
D. P. Chen, R. S. Eisenberg, J. W. Jerome, and C. W. Shu, {\it A
hydrodynamic model of temperature change in open ionic channels},
Biophys J., {\bf 69} (1995), 2304--2322.






\bibitem{DeMark1d}
P. Degond and  P. A. Markowich, {\it On a one-dimensional
steady-state hydrodynamic model for semiconductors}, Appl. Math.
Lett., {\bf 3} (1990),  25--29.


\bibitem{DeMark3d}
P. Degond and  P. A. Markowich, {\it A steady state potential flow
model for semiconductors}, Ann. Mat. Pura Appl. (4)  {\bf 165}
(1993), 87--98.


\bibitem{Fried}
K. O.  Friedrichs, {\it Symmetric positive linear differential equations}, Comm. Pure Appl. Math. {\bf 11} (1958), 333--418.

\bibitem{Gamba1d}
I. M. Gamba, {\it Stationary transonic solutions of
a one-dimensional hydrodynamic model for semiconductors}, Comm.
Partial Differential Equations, {\bf 17} (1992),  553--577.


\bibitem{GambaMorawetz}
I. M. Gamba and C. S. Morawetz, {\it A viscous
approximation for a {$2$}-{D} steady semiconductor or transonic gas
dynamic flow: existence theorem for potential flow}, Comm. Pure
Appl. Math., {\bf 49} (1996),  999--1049.


\bibitem{GilbargTrudinger}
D. Gilbarg and N. Trudinger,
{\it  Elliptic Partial Differential Equations of Second Order.}  2nd Ed. Springer-Verlag: Berlin.

\bibitem{Guo99}
Y. Guo, {\it Smooth irrotational flows in the large to the Euler-Poisson system in $\R^{3+1}$},   Comm. Math. Phys., {\bf 195} (1998), 249--265.

\bibitem{GHZ}
Y. Guo, L. J.  Han, J. J. Zhang, {\it  Absence of shocks for one dimensional Euler-Poisson system},  Arch. Ration. Mech. Anal., {\bf 223} (2017),  1057--1121.

\bibitem{Guo}
Y. Guo and W. Strauss, {\it Stability of semiconductor states with insulating and contact boundary conditions},
 Arch. Ration. Mech. Anal., {\bf 179} (2006),  1--30.


\bibitem{Ha-L}
Q. Han and F. Lin,
{\it Elliptic partial differential equations.}  Courant Institute of Math. Sci., NYU.

\bibitem{Huang}
F. M. Huang, R. H. Pan, and H. M. Yu,  {\it Large time behavior of Euler-Poisson system for semiconductor},
Sci. China Ser. A, {\bf 51} (2008),  965--972.

\bibitem{IP}
A. D. Ionescu and B. Pausader,{\it The Euler-Poisson system in 2D: global stability of the constant equilibrium solution}, Int. Math. Res. Not., (2013),  761--826.

\bibitem{LW}
D. Li and Y. F. Wu, {\it The Cauchy problem for the two dimensional Euler-Poisson system}, J. Eur. Math. Soc., {\bf  16} (2014),  2211--2266.

\bibitem{LiMark}
H. L. Li,   P. Markowich, and M. Mei,  {\it Asymptotic behavior of subsonic entropy solutions of the isentropic Euler-Poisson equations},  Quart. Appl. Math.,  {\bf 60} (2002),  773--796.

\bibitem{LXY1}
J. Li, Z. P. Xin, and H. C. Yin, {\it On transonic shocks in a
nozzle with variable end pressures},  Comm. Math. Phys., {\bf 291} (2009),  111--150.



\bibitem{LuoNX}
T. Luo, R. Natalini, and Z. P. Xin, {\it Large time behavior of the solutions to a hydrodynamic model for semiconductors},  SIAM J. Appl. Math., \textbf{59} (1999), 810--830.

\bibitem{LRXX}
T. Luo, J. Rauch,  C. J. Xie, and Z. P. Xin, {\it Stability of transonic shock solutions for one-dimensional Euler-Poisson equations}, Arch. Ration. Mech. Anal., {\bf 202} (2011), 787--827.

\bibitem{LuoXin}
T. Luo and Z. P. Xin, {\it Transonic shock solutions for a
system of Euler-Poisson equations}, Comm. Math. Sci., {\bf 10} (2012),  419--462.



\bibitem{MarkZAMP}
P. A. Markowich, {\it On steady state {E}uler-{P}oisson models
for semiconductors}, Z. Angew. Math. Phys., {\bf 42} (1991),
389--407.


\bibitem{MarkRSbook}
P. A. Markowich, C. A. Ringhofer, and C. Schmeiser, {\it
Semiconductor equations}, Springer-Verlag, Vienna, 1990.


\bibitem{Peng}
Y. J. Peng and I. Violet,  {\it Example of supersonic
solutions to a steady state {E}uler-{P}oisson system}, Appl. Math. Lett., {\bf 19} (2006),   1335--1340.

\bibitem{Rauch}
J. Rauch and F. Massey, {\it Differentiability of solutions to hyperbolic initial-boundary value problems}, Trans. Amer. Math. Soc., {\bf189} (1974), 303--318.

\bibitem{RosiniStability}
M. D. Rosini, {\it Stability of transonic strong shock
waves for the one-dimensional hydrodynamic model for
semiconductors}, J. Differential Equations, {\bf 199} (2004),
326--351.


\bibitem{RosiniPhase}
M. D. Rosini, {\it A phase analysis of transonic solutions
for the hydrodynamic semiconductor model}, Quart. Appl. Math., {\bf
63} (2005),  251--268.





\bibitem{WEP}
S. K.  Weng, {\it On steady subsonic flows for Euler-Poisson models},  SIAM J. Math. Anal., {\bf 46}  (2014),   757--779.




\bibitem{Yeh}
L. M. Yeh, {\it On a steady state Euler-Poisson model for semiconductors}, Comm. Partial Differential Equations, {\bf 21} (1996), 1007--1034.


\bibitem{Hattori}
C. Zhu and H. Hattori, {\it Asymptotic behavior of the solution to a nonisentropic hydrodynamic model of semiconductors}, J. Differential Equations, {\bf 144} (1998),  353--389.


\end{thebibliography}
\end{document}